\documentclass[11pt,reqno]{amsart}

\usepackage[margin=27mm]{geometry}
\usepackage[T1]{fontenc}
\usepackage{lmodern}
\usepackage{microtype}
\usepackage{mathtools}
\usepackage{amssymb}
\usepackage{mathrsfs}
\usepackage{booktabs}
\usepackage{array}
\usepackage{longtable}
\usepackage{enumitem}
\usepackage{graphicx}
\usepackage[section]{placeins}
\usepackage[numbers,sort&compress]{natbib}
\usepackage[colorlinks,linkcolor=blue,citecolor=blue,urlcolor=blue]{hyperref}
\usepackage[nameinlink,noabbrev,capitalise]{cleveref}

\providecommand{\doi}[1]{doi:~\href{https://doi.org/#1}{\nolinkurl{#1}}}
\numberwithin{equation}{section}
\makeatletter
\def\LT@makecaption#1#2#3{\LT@mcol\LT@cols c{\hbox to\z@{\hss\parbox[t]\LTcapwidth{\@makecaption{#2}{#3}\vskip\belowcaptionskip}\hss}}}
\makeatother
\setlist[enumerate]{label=(\alph*),leftmargin=*,itemsep=2pt}
\theoremstyle{plain}
\newtheorem{theorem}{Theorem}[section]
\newtheorem{lemma}[theorem]{Lemma}
\newtheorem{proposition}[theorem]{Proposition}

\theoremstyle{definition}
\newtheorem{definition}[theorem]{Definition}
\theoremstyle{remark}
\newtheorem{remark}[theorem]{Remark}

\newcommand{\R}{\mathbb{R}}
\newcommand{\C}{\mathbb{C}}
\newcommand{\Z}{\mathbb{Z}}
\newcommand{\Tor}{\mathbb{T}}
\newcommand{\one}{\mathbf{1}}
\newcommand{\T}{\mathsf{T}}
\newcommand{\cp}{\operatorname{cap}}

\newcommand{\bn}{\mathbf{n}}
\newcommand{\bt}{\mathbf{t}}
\newcommand{\Jm}{J_{-}}

\newcommand{\Ih}{I_h}

\newcommand{\Lam}{\Lambda}

\newcommand{\Lh}{\Lambda_h}
\newcommand{\Lhp}{\Lambda_h^{\dagger}}
\newcommand{\SL}{\mathcal{V}}
\newcommand{\Pz}{P_0}
\newcommand{\avg}[1]{\langle #1\rangle_\Gamma}
\newcommand{\mV}{\mathsf{V}}
\newcommand{\mK}{\mathsf{K}}
\newcommand{\mW}{\mathsf{W}}
\newcommand{\mB}{\mathsf{B}}
\newcommand{\mDs}{\mathsf{D}_s}
\newcommand{\mG}{\mathsf{G}}
\newcommand{\Mv}{M_v}

\newcommand{\dnorm}[1]{\lVert #1\rVert_*}

\newcommand{\Om}{\Omega}

\newcommand{\bXh}{b_{X,h}}
\newcommand{\Nd}{\mathcal{N}_{r_0}}

\newcommand{\epsh}{\varepsilon_h}

\newcommand{\dd}{\mathrm{d}}
\newcommand{\iu}{\mathrm{i}}
\newcommand{\eu}{\mathrm{e}}

\title[Shape perturbation of the SGBEM Dirichlet-to-Neumann matrix]{Shape
 Perturbation of the Galerkin Boundary Element
 Dirichlet-to-Neumann Matrix on Moving Polygons}
\author{Tokuhiro Eto}
\address{Universit\'{e} Claude Bernard Lyon 1, CNRS, Centrale Lyon,
 INSA Lyon, Universit\'{e} Jean Monnet, ICJ UMR5208,
 69622 Villeurbanne, France}
\email{eto@math.univ-lyon1.fr}
\date{}

\begin{document}

\begin{abstract}
For nonlocal interface evolution problems such as the Mullins--Sekerka flow,
certain parametric schemes discretize only the interface and approximate the Dirichlet-to-Neumann map by boundary elements.
Their error analysis requires bounds on how this discrete operator changes as the polygon moves.
In that analysis, the polygons are only Lipschitz and
close to a smooth curve only in a Lipschitz norm, the perturbations are
controlled only in the energy norm, and the kernels are singular on the
diagonal. We study the symmetric Galerkin boundary element matrix on the
polygons in a Lipschitz neighborhood of a smooth reference curve. Its first
variation, tested against discrete Neumann solutions, is bounded by the
energy norm of the perturbation alone. The second-order remainder of the
associated mass form is bounded by the product of the energy norm and the
Lipschitz norm of the perturbation, so that only one of the two factors
needs a uniform bound. On the interpolating polygon, the matrix is
consistent with the continuous operator without the logarithmic loss that a
pointwise comparison of the kernels would suggest, and its shape derivative
is consistent with that of the continuous operator. All bounds are uniform in the mesh size and depend
only on the constants of the reference curve.
\end{abstract}


\maketitle

\section{Introduction}\label{sec:intro}

In the Mullins--Sekerka flow \cite{MS63}, an interface moves with a normal
velocity given by the jump of the normal derivative of a harmonic chemical
potential, whose boundary value is prescribed by the curvature. In the
parametric scheme proposed by the author \cite{E26-SGBEM}, the chemical potential is eliminated
by a symmetric Galerkin boundary element (SGBEM) approximation of the
Dirichlet-to-Neumann map, so that only the curve is discretized and no bulk
mesh is needed. In every step, the Dirichlet-to-Neumann map of the current
polygon $X$ enters through one matrix, the SGBEM Steklov--Poincar\'e matrix
$\Lh(X)$, which is assembled from piecewise constant densities and
piecewise linear traces on $X$. The polygon moves from step to step, and so
does the matrix.

This paper asks how the matrix $\Lh(X)$ changes when the polygon $X$ moves, in which
norm of the direction of motion the change can be measured, and whether the
resulting bounds are uniform in the mesh size $h$ and in the polygon $X$.
The question arises in the error analysis of such schemes. There, the polygons
are close to a smooth curve but are only Lipschitz. The error analysis
controls the direction $e$ of a perturbation only in $H^1$, that is, $e'$
only in $L^2$. Moreover, the kernels of the boundary integral operators are
singular on the diagonal. A bound that measures $e$ in $W^{1,\infty}$ is
easy to obtain and too weak for this purpose, since every factor
$\|e'\|_{L^\infty}$ costs an inverse estimate $h^{-1/2}$ when it is
converted to $\|e\|_{H^1}$.

We work near a fixed reference curve $x$ from a family $\mathscr K$ of smooth
embedded closed curves on which three constants are uniform: a chord
constant, bounds for finitely many derivatives, and the coercivity constant
of the single layer operator together with a bound for its inverse
(\cref{ass:K}). The capacity condition $\operatorname{cap}\Gamma<1$ is one
way to obtain the coercivity (\cref{rem:cap}). The polygons are those
that are close to the nodal interpolant $\Ih x$ in $W^{1,\infty}$
(\cref{def:nbhd}). On this class, we prove the following bounds, with
constants that depend only on the constants of $\mathscr K$.
\begin{enumerate}[label=(\roman*)]
\item The first variation of $\Lh(X)$, tested against the solutions
$z_f$ of the discrete Neumann problem with data $f$ (\cref{sec:setting}), is bounded by
$\|e'\|_{L^2}$. When the segment from $X$ to $X+e$ stays in the class,
the increment of the associated mass form $\bXh$ is bounded by
$\|e'\|_{L^2}$ times one $L^\infty$ norm of the data
(\cref{thm:bem-c1p}).
\item Under the same segment condition, the second-order Taylor remainder
of $\bXh$ is bounded by $\|e'\|_{L^\infty}\|e'\|_{L^2}$, with only one
factor in $L^\infty$ (\cref{thm:bem-c3p}).
\item The matrix $\Lh(X)$ is Lipschitz in the polygon with respect to
$\|X'-Y'\|_{L^\infty}$, with quadratic Taylor remainders and mixed second
derivatives (\cref{thm:bem-c1}).
\item On the interpolating polygon $\Ih x$, applied to the interpolant of
$g\in H^3$, the matrix is consistent with the continuum
Dirichlet-to-Neumann map of order $h^{3/2}\|g\|_{H^3}$ in $(H^{1/2})^*$
and $h^2\|g\|_{H^3}$ in $(H^1)^*$, without a
logarithmic loss, and its shape derivative is consistent with the shape
derivative of the continuum operator, of order $h\|e'\|_{L^\infty}\|g\|_{H^2}$
in $(H^1)^*$ and $h^{1/2}\|e'\|_{L^\infty}\|g\|_{H^2}$ in $(H^{1/2})^*$
(\cref{lem:bem-bweak,lem:bem-c4}). Here, the reference curve has four
bounded derivatives.
\end{enumerate}
The bounds (i) and (ii) measure the direction in the weak norm, with a
single factor $\|e'\|_{L^\infty}$ in (ii), as such an error analysis needs. The bound (iii) is in the strong norm $\|e'\|_{L^\infty}$, and the
proof of (iv) uses it.

The strategy of the proofs is the following. Written on a fixed parameter circle, the
SGBEM matrices depend on the polygon only through the logarithm of the chord
$X(\xi)-X(\eta)$, and the arclength weights cancel exactly
(\cref{lem:bem-reduction,lem:bem-log}). The shape derivatives of the
operators are therefore powers of the divided-difference quotient
$q_e=(e(\xi)-e(\eta))/(X(\xi)-X(\eta))$. The first variation is controlled
by a maximal function of $e'$, which is bounded in $L^2$ by the
Hardy--Littlewood theorem, and the second-order remainder of the mass form, through a Calder\'on
commutator \cite{CoifmanMeyer1978}, by the boundedness of the Cauchy
integral on Lipschitz graphs \cite{CMM1982} and its uniform extension to
chord-arc curves \cite{David1984}.
In the consistency estimate, the single layer operators of the polygon and
of the curve differ by $O(h^2)$ on $L^2$, while a bound of the kernel
difference in absolute value gives only $O(h^2|\log h|)$. The logarithm
disappears because the leading part of the kernel difference is a
commutator with the Cauchy integral.

The bounds above concern the SGBEM matrix $\Lh(X)$ and its associated mass
form $\bXh$ on moving piecewise linear polygons, with constants uniform in
$h$. Two lines of work are related to them. For the
continuum operators, shape differentiability of boundary integral operators is classical for smooth
boundaries \cite{Potthast1994,CostabelLeLouer2012}, and shape holomorphy of
the Calder\'on projector of the Laplacian in $\R^2$ is known for $C^2$
Jordan curves \cite{HenriquezSchwab2021}. Parametric shape holomorphy of
boundary integral operators in an $L^2$-based setting has been extended to
Lipschitz boundaries in $\R^3$ \cite{DolzHenriquez2024} and to multiple open
arcs \cite{PHJH2024}. A recent report \cite{HDPS2026SAM} proves shape
holomorphy of the Calder\'on projector of the Laplacian on Lipschitz
boundaries, between fractional Sobolev trace spaces, and notes that this is
the analytic foundation for high-dimensional approximation of the
transported operators and their Galerkin discretizations. For the
conductivity equation in a planar domain, the Dirichlet-to-Neumann map is
differentiable with respect to the vertices of a triangular inclusion
\cite{BeFrVe2017}.

The second line is the shape sensitivity analysis of boundary integral
formulations and boundary element methods, which has a long history \cite{Bonnet1997}, including
symmetric Galerkin formulations for cracks \cite{Bonnet2001}.
\section{Setting}\label{sec:setting}

\subsection{Curves, polygons and norms}
Let $\Tor=\R/\Z$ and let $d$ be the distance on $\Tor$, that is, $d(\xi,\eta) := \min_{k\in\mathbb Z}|\xi - (\eta + k)|$ for $\xi,\,\eta\in\Tor$.
We identify $\R^2$ with $\C$, and $\Jm$ denotes the rotation by $-90^\circ$, that is, $\Jm z=-\iu z$ for $z\in\C$.
A closed curve is parametrized by a map $Y:\Tor\to\R^2$. It is embedded if it is injective, smooth if it is $C^\infty$, and oriented counter-clockwise if its signed area
\[
\frac12\int_\Tor\mathrm{Im}\big(\overline{Y(\xi)}\,Y'(\xi)\big)\,\dd\xi
\]
is positive. For a Lipschitz curve with $|Y'|>0$, the unit tangent is $\bt=Y'/|Y'|$ and the unit normal is $\bn=\Jm\bt$, where $'=\dd/\dd\xi$ and $|\cdot|$ is the Euclidean norm.
If $Y$ is embedded and oriented counter-clockwise, then $\bn$ points out of the bounded domain enclosed by $Y(\Tor)$.
A closed Lipschitz curve $Y$ has chord constant $c>0$ if $|Y(\xi)-Y(\eta)|\ge c\,d(\xi,\eta)$ for all $\xi,\eta\in\Tor$. Such a curve is embedded and satisfies $|Y'|\ge c$ almost everywhere. If moreover $|Y'|\le M$ almost everywhere, then the arclength between two points of $Y(\Tor)$ is bounded by $M/c$ times their distance, so that $Y(\Tor)$ is a chord-arc curve in the usual sense.
For $m\ge1$ and a function $f:\Tor\to\R^k$ whose derivatives up to the order $m-1$ are Lipschitz, we set
$\|f\|_{W^{m,\infty}}=\sum_{i=0}^m\|f^{(i)}\|_{L^\infty(\Tor)}$, where $\cdot^{(i)} = \mathrm{d}^i/\mathrm{d}\xi^i$ for $i \ge 1$.
In particular, $\|f\|_{W^{1,\infty}}=\|f\|_{L^\infty(\Tor)}+\|f'\|_{L^\infty(\Tor)}$.
Two parametrized curves $x,y:\Tor\to\R^2$ are compared through $\|x-y\|_{W^{1,\infty}}$ as functions, not their images.

A polygon is identified with a continuous piecewise linear map $X:\Tor\to\R^2$ on the
grid $\xi_j=jh\,(0\le j\le N)$ with the uniform mesh size $h=1/N$.
It has vertices $X_j = X(\xi_j)$ and edges $\Gamma_p=[X_p,X_{p+1}]\,(0\le p\le N-1)$, with the
convention $X_0=X_N$, and the indices of vertices, edges and vertex values are taken modulo $N$. The edge $\Gamma_p$ is the image of the parameter interval $[\xi_p,\xi_{p+1}]$, which we call a \emph{panel}, and $\Ih$ denotes the nodal interpolation on this grid.
We identify a polygon, and a displacement $U$ of its vertices, with the
vector of its $N$ distinct vertices in $\R^{2N}$, so that $U=(U_1,\dots,U_N)$
with $U_j\in\R^2$. On $\Tor$, we use the fixed spaces $\mathcal P_1$ of
continuous piecewise linear functions (traces), with the hat functions
$\varphi_j$, composing the basis of $\mathcal P_1$,
and $\mathcal P_0$ of piecewise constant functions (densities), with the
indicators $\iota_p$ of the panels. All Sobolev norms are taken on $\Tor$ in the parameter. For $s\in\R$, the norm of $H^s=H^s(\Tor)$ is defined by
\[
\|u\|_{H^s}= \sqrt{\sum_{n\in\Z}(1+n^2)^s\,|\hat u(n)|^2},\qquad \hat u(n)=\int_\Tor \eu^{-2\pi\iu n\xi}u(\xi)\,\dd\xi ,
\]
and for $s<0$ the space $H^s$ is the completion of $L^2$ in this norm. We write $\langle u,v\rangle=\int_\Tor u\,v\,\dd\xi$ for real functions $u$ and $v$, and also for its extension to the duality between $H^{-s}$ and $H^s$. A vector of vertex values is measured through its piecewise linear interpolant.
For a load vector $y\in\R^N$, we use the dual norms of $(H^{1/2})^*$ and $(H^1)^*$ on $\mathcal P_1$,
\[
\dnorm{y}=\sup_{w\in\mathcal P_1\setminus\{0\}}\frac{y^\T w}{\|w\|_{H^{1/2}}},\qquad
\|y\|_{(H^1)^*}=\sup_{w\in\mathcal P_1\setminus\{0\}}\frac{y^\T w}{\|w\|_{H^1}} .
\]

\subsection{The SGBEM Dirichlet-to-Neumann matrix}
With
\[
\omega_j=\frac12\Jm(X_{j+1}-X_{j-1})\in\R^2 ,
\]
let
$\Om(X)=\operatorname{blockdiag}(\omega_1^\T,\dots,\omega_N^\T)\in\R^{N\times2N}$
be the block diagonal matrix with the $1\times2$ blocks $\omega_j^\T$, so that
$\big(\Om(X)U\big)_j=\omega_j\cdot U_j$ for $1\le j\le N$.
Roughly speaking, these values represent a normal velocity when $U$ corresponds to the nodal displacement.
The lumped vertex mass is $\Mv=\operatorname{diag}(\mu_j)$ with $\mu_j=\tfrac12(|\Gamma_{j-1}|+|\Gamma_j|)$.

The matrix $\Lh(X)$ is the SGBEM Steklov--Poincar\'e matrix introduced in \cite[(5.19)]{E26-SGBEM}.
It is assembled from the single-layer and the
double-layer operators of the Laplacian:
\[
\SL g(y)=-\frac1{2\pi}\int_\Gamma\log|y-y'|\,g(y')\,\dd s_{y'},\qquad
\mathcal K g(y)=-\frac1{2\pi}\int_\Gamma\partial_{\bn(y')}\log|y-y'|\,g(y')\,\dd s_{y'},
\]
where $\partial_{\bn(y')}$ is the derivative in $y'$ along the normal $\bn(y')$ of $\Gamma$,
through their Galerkin matrices on piecewise constant densities and piecewise linear traces on the polygon. On $\Gamma=X(\Tor)$, we write $\varphi_j$ and $\iota_p$ also for $\varphi_j\circ X^{-1}$ and $\iota_p\circ X^{-1}$, so that $\iota_p$ is the indicator of the edge $\Gamma_p$. The matrices are
\[
\begin{array}{ll}
\mV,\quad \mV_{pq}=\int_{\Gamma_p}\SL\iota_q\,\dd s & \text{single layer, densities to densities,}\\
\mK,\quad \mK_{pj}=\int_{\Gamma_p}\mathcal K\varphi_j\,\dd s & \text{double layer, traces to densities,}\\
\mB,\quad \mB_{pj}=\int_{\Gamma_p}\varphi_j\,\dd s & \text{mass between traces and densities,}\\
\mDs,\quad (\mDs u)_p=(u_{p+1}-u_p)/|\Gamma_p| & \text{arclength derivative, traces to densities,}\\
\mW=\mDs^\T\mV\mDs & \text{hypersingular (Maue's formula).}
\end{array}
\]
The Galerkin integrals are evaluated exactly.
With $Q_{\rm int}=\tfrac12\mB+\mK$, $Q_{\rm ext}=\tfrac12\mB-\mK$ and the
Euclidean deflation $P=I-\one\one^\T/N$, the matrix is represented as
\[
\Lh(X)=P\big(2\mW+Q_{\rm int}^\T\mV^{-1}Q_{\rm int}+Q_{\rm ext}^\T\mathsf C\,Q_{\rm ext}\big)P,
\]
where $\mathsf C$ is the bordered inverse of $\mV$ on densities of arclength
mean zero ($\int\varsigma\,\dd s=0$), that is, the upper left $N\times N$ block of the inverse of the bordered matrix
\[
\begin{pmatrix}\mV&\mB\one\\ (\mB\one)^\T&0\end{pmatrix},\qquad (\mB\one)_p=|\Gamma_p| .
\] The formula requires $\mV$ and the bordered system to be invertible. For the polygons near a reference curve considered
below, both matrices are invertible by the coercivity of \cref{lem:bem-Vreg}(a),
with inverses bounded uniformly in the norms of \eqref{eq:bem-inv}.
When the polygon is clear from the context, we write $\Lh$ for $\Lh(X)$. Since $\mV$ is symmetric, so are $\mW$, $\mV^{-1}$ and $\mathsf C$, and since $P$ is symmetric with $P\one=0$, the matrix $\Lh$ is symmetric and $\Lh\one=0$ at every polygon at which it is defined.

We write $\Lhp$ for the Moore--Penrose inverse of $\Lh$, that is, the unique matrix $B$ with $\Lh B\Lh=\Lh$ and $B\Lh B=B$ for which $\Lh B$ and $B\Lh$ are symmetric. Let $\one^\perp=\{y\in\R^N\mid\one^\T y=0\}$ and $\Pi_{\one}=\one\one^\T/N$, so that $P=I-\Pi_{\one}$. If $\Lh+\Pi_{\one}$ is invertible, then we have
\begin{equation}\label{eq:mp-formula}
\Lhp=\big(\Lh+\Pi_{\one}\big)^{-1}-\Pi_{\one},\qquad \Lhp P=\Lhp,\qquad \one^\T\Lhp=0,\qquad \Lh\Lhp=P .
\end{equation}
Indeed, let $M=\Lh+\Pi_{\one}$ and $B=M^{-1}-\Pi_{\one}$. Since $\Lh\Pi_{\one}=\Pi_{\one}\Lh=0$ and $\Pi_{\one}^2=\Pi_{\one}$, we have $M\Pi_{\one}=\Pi_{\one}M=\Pi_{\one}$, so $B\Pi_{\one}=\Pi_{\one}B=0$ and $\Lh B=B\Lh=P$, from which the four conditions that characterize $\Lhp$ and the other three identities in \eqref{eq:mp-formula} follow.
In particular, for every $y\in\R^N$, the vector $\Lhp y=\Lhp Py$ lies in $\one^\perp$ and solves $\Lh\Lhp y=Py$.
If $\ker\Lh=\operatorname{span}\{\one\}$, which holds for the polygons considered below by \cref{lem:bem-a}, then $\Lh+\Pi_{\one}$ is invertible. Indeed, if $(\Lh+\Pi_{\one})v=0$, then multiplying by $\one^\T$ and using $\one^\T\Lh=0$ gives $\one^\T v=0$, so $\Lh v=0$, and $v\in\operatorname{span}\{\one\}\cap\one^\perp=\{0\}$.

Let $X$ be a polygon at which $\Lh$ is differentiable and $\Lh(X)+\Pi_{\one}$ is invertible. Then, $\Lhp$ is differentiable at $X$, and for every direction $e$, we have
\begin{equation}\label{eq:mp-deriv}
D\Lhp(X)[e]=-\Lhp(X)\,D\Lh(X)[e]\,\Lhp(X).
\end{equation}
Indeed, $\Lh$ is continuous at $X$ and the invertible matrices form an open set, so \eqref{eq:mp-formula} holds at every polygon near $X$. Differentiating the inverse in \eqref{eq:mp-formula}, we obtain $D\Lhp(X)[e]=-M^{-1}HM^{-1}$ with $M=\Lh(X)+\Pi_{\one}$ and $H=D\Lh(X)[e]$, and since $H$ is symmetric with $H\one=0$, as $\Lh$ is at every polygon near $X$, the terms with $\Pi_{\one}$ in $M^{-1}=\Lhp(X)+\Pi_{\one}$ drop out.

Loads are written $\Om(X)f$, where $\Om(X)$ acts on the vertex values of a
piecewise linear vector field $f$. We set
\[
z_f=\Lhp(X)\Om(X)f,\qquad
\bXh(f,\varphi)=2\big(\Om(X)\varphi\big)^\T\Lhp(X)\big(\Om(X)f\big).
\]
The form $\bXh$ is the discrete mass of the parametric scheme of \cite{E26-SGBEM}.

\subsection{The continuum operators}
Let $\Lam_\Gamma=\Lam^{\rm int}_\Gamma+\Lam^{\rm ext}_\Gamma$
be the two-phase Dirichlet-to-Neumann map of a Jordan curve $\Gamma$ that bounds a Lipschitz domain, the exterior part
being the one for functions bounded at infinity.
Let $\avg{g}=\tfrac{1}{|\Gamma|}\int_\Gamma g\,\dd s$, and let $\Pz g=g-\avg g$ be the
$L^2(\dd s)$-orthogonal projection onto the arclength-mean-zero function space.

The operators of a curve are applied to functions of the parameter through the parametrization. For an embedded closed Lipschitz curve $Y:\Tor\to\R^2$ whose image $\Gamma_Y=\Gamma(Y)=Y(\Tor)$ is such a curve, and for a function $g$ on $\Tor$, we write
\[
\Lam_Yg=\Lam_{\Gamma_Y}g=\big(\Lam_{\Gamma_Y}(g\circ Y^{-1})\big)\circ Y .
\]
In particular, $\Lam_\Gamma g$ for the reference curve $x$ with $\Gamma=x(\Tor)$, and $\Lam_{x+\varepsilon e}g=\Lam_{\Gamma(x+\varepsilon e)}g$ for a perturbed curve, are functions on $\Tor$.

For a closed Lipschitz curve $Y:\Tor\to\R^2$ with a positive chord constant, in particular for an embedded polygon with edges of positive length or a smooth embedded closed curve with $|Y'|>0$, the parameter kernels are
\[
\mathcal L_Y\varsigma(\xi)=-\frac1{2\pi}\int_\Tor\log|Y(\xi)-Y(\eta)|\,\varsigma(\eta)\,\dd\eta,\qquad
\mathcal K_Ya(\xi)=-\frac1{2\pi}\,\mathrm{p.v.}\!\int_\Tor\mathrm{Im}\,\frac{Y'(\eta)}{Y(\eta)-Y(\xi)}\,a(\eta)\,\dd\eta.
\]
The logarithmic kernel is weakly singular, $\big|\log|Y(\xi)-Y(\eta)|\big|\le|\log d(\xi,\eta)|+C$ with $C$ depending only on $\|Y'\|_{L^\infty}$ and the chord constant, so $\mathcal L_Y\varsigma$ is an absolutely convergent integral for $\varsigma\in L^2$,
while the kernel of $\mathcal K_Y$ is of order $d(\xi,\eta)^{-1}$ on a general Lipschitz curve, and the principal value is the limit as $\delta\to0$ of the integrals over $\{\eta\in \Tor\mid d(\xi,\eta)>\delta\}$.

\subsection{The family of reference curves}

\begin{definition}[The hypotheses (K)]\label{ass:K}
Let $\mathscr K$ be a family of smooth embedded closed curves $x:\Tor\to\R^2$, oriented counter-clockwise, for which the following constants are uniform on $\mathscr K$.
\begin{enumerate}[label=(K\arabic*)]
\item\label{K1} The curve satisfies $|x'|\le M_*$ and has chord constant $c_*$.
\item\label{K2} For the order $r\ge0$ named in each result, $\|x\|_{W^{r+1,\infty}}\le C_r$. The order named is sufficient for the proof and is not claimed to be minimal.
\item\label{K3} The single layer is uniformly coercive on $L^2$ densities, $\langle\mathcal L_x\varsigma,\varsigma\rangle\ge\mu_1\|\varsigma\|_{H^{-1/2}}^2$ for all $\varsigma\in L^2$ with $\mu_1>0$. It maps $L^2$ boundedly onto $H^1$, and $\|\varsigma\|_{L^2}\le B_1\|\mathcal L_x\varsigma\|_{H^1}$ for all $\varsigma\in L^2$.
\end{enumerate}
We call $M_*$, $c_*$, $\mu_1$, $B_1$ and the $C_r$ the constants of $\mathscr K$.
\end{definition}

\begin{remark}[The capacity condition]\label{rem:cap}
Let $\Gamma=x(\Tor)$ for a smooth embedded closed curve $x$. The equilibrium density $w_{\rm eq}$ is the unique $w\in H^{-1/2}(\Gamma)$ with $\int_\Gamma w\,\dd s=1$ whose single layer potential $\SL w$ is constant on $\Gamma$, $\lambda_{\rm eq}(\Gamma)$ is this constant, and $\cp\Gamma=\eu^{-2\pi\lambda_{\rm eq}(\Gamma)}$ is the logarithmic capacity of $\Gamma$ \cite[\S6.6.1]{Steinbach2008}. If $\cp\Gamma<1$, then $\lambda_{\rm eq}(\Gamma)>0$, and the argument of the proof of \cite[Thm.~6.23]{Steinbach2008}, which combines this with the coercivity of $\SL$ on densities of mean zero \cite[Thm.~6.22]{Steinbach2008}, shows that $\SL$ is coercive on $H^{-1/2}(\Gamma)$.
\end{remark}

\subsection{Neighborhoods of a reference curve}
Let $C_L$ be the constant of \cref{lem:bem-cont} for the family of Lipschitz curves $Y$ with $|Y'|\le M_*+c_*/2$ and chord constant at least $c_*/2$. It depends only on $M_*$ and $c_*$. Let $C_I$ be the constant in $\|(\Ih x-x)'\|_{L^\infty}\le C_Ih$, which depends only on $C_1$.

\begin{definition}\label{def:nbhd}
Let $x\in\mathscr K$. Choose $r_0,h_0>0$ and put $r_1=r_0+C_Ih_0$ such that
\begin{equation}\label{eq:nbhd}
r_1\le\tfrac14c_*,\qquad
2C_Lr_1\le\min\big(\tfrac12\mu_1,(2B_1)^{-1}\big).
\end{equation}
Let $\mathfrak C$ be the set of Lipschitz curves $Y$ with $\|Y'-x'\|_{L^\infty}\le2r_1$, and let $\Nd$ be
the set of polygons $X$ with $\|X-\Ih x\|_{W^{1,\infty}}\le r_0$, for $h\le h_0$.
\end{definition}

Both conditions in \eqref{eq:nbhd} involve only the constants of $\mathscr K$, and the second gives the uniform bounds of \cref{lem:bem-Vreg}.

Since $\|(\Ih x-x)'\|_{L^\infty}\le C_Ih$, every $X\in\Nd$ satisfies $\|X'-x'\|_{L^\infty}\le r_1$, so $\Nd\subset\mathfrak C$. For every polygon $Z\in\mathfrak C$, we have $\ker\Lh(Z)=\operatorname{span}\{\one\}$ by \cref{lem:bem-a}, so $\Lh(Z)+\Pi_{\one}$ is invertible. The polygons $Z$ with $\|Z'-x'\|_{L^\infty}<2r_1$ form an open set contained in $\mathfrak C$, so $\Lh$ is defined near each of them. By \cref{lem:bem-log}, the matrices $\tilde V$ and $\tilde K$ of \cref{lem:bem-reduction} are twice differentiable functions of the vertices at $Z$, and by \cref{lem:bem-reduction} and \cref{lem:bem-Vreg}(a), $\Lh$ is a smooth function of $\tilde V$ and $\tilde K$ near $Z$. Hence, $\Lh$ is twice differentiable at $Z$. Thus, \eqref{eq:mp-deriv} holds at every polygon $Z$ with $\|Z'-x'\|_{L^\infty}<2r_1$, in particular on a neighborhood of every polygon of $\Nd$. Since $|\Delta(Y-x)|\le2r_1d(\xi,\eta)$ for $Y\in\mathfrak C$, where $\Delta Y(\xi,\eta)=Y(\xi)-Y(\eta)$, the curves $Y\in\mathfrak C$ have a uniform chord constant. Indeed, with $M_1=M_*+2r_1$ and $c_{\rm ca}=c_*-2r_1\ge c_*/2$, we have
\[
|Y'|\le M_1,\qquad |Y(\xi)-Y(\eta)|\ge c_{\rm ca}\,d(\xi,\eta)\qquad\text{for all }Y\in\mathfrak C,\ \xi,\eta\in\Tor .
\]
The set $\mathfrak C$ is convex and contains $\Nd$, the curve $x$, the segments between them, and $x+\varepsilon e$ for every Lipschitz direction $e$, that is, every Lipschitz map $e:\Tor\to\C$, with $\|\varepsilon e'\|_{L^\infty}\le2r_1$. By \cref{prop:bem-dtn}(i), the polygons of $\mathfrak C$, the curve $x$ and the curves $x+\varepsilon e$ with a piecewise linear direction $e$ and $|\varepsilon|\,\|e'\|_{L^\infty}\le2r_1$ are Jordan curves, oriented counter-clockwise, that bound Lipschitz domains. The matrix $\Lh(Y)$ is defined for $Y\in\mathfrak C$ by the formula of \cref{lem:bem-reduction}.

Throughout, $C$ and $c$ denote constants that depend only on the constants of $\mathscr K$ for the order $r$ named in the result at hand, and on $r_0$, $h_0$ and $r_1$, which are themselves fixed by these constants. They do not depend on $x\in\mathscr K$ otherwise, nor on the polygon, on $h$ or $N$, or on the functions and directions to which the estimates are applied. Every estimate below whose left side is linear in a function or a direction carries the norm of that function or direction as a factor.

\section{Shape perturbation of the boundary element operators}\label{sec:shape}

This section states the main results and the two structural facts behind them, \cref{lem:bem-reduction,lem:bem-log}. The proofs of the stated results are given in \cref{sec:app-bem}. Throughout, $x\in\mathscr K$, and $\Nd$, $\mathfrak C$, $r_0$, $h_0$, $r_1$, $M_1$ and $c_{\rm ca}$ are as in \cref{def:nbhd}.

\subsection{The parameter circle}

\begin{proposition}[The parameter circle]\label{lem:bem-reduction}
Let $X$ be a polygon with a positive chord constant, and let $\tilde V_X$ be the Galerkin matrix of $\mathcal L_X$ on $\mathcal P_0$, with the entries $(\tilde V_X)_{pq}=\langle\mathcal L_X\iota_q,\iota_p\rangle$. Assume that $\tilde V_X$ is invertible, and that the bordered matrix
\[
\begin{pmatrix}\tilde V_X&h\one\\ h\one^\T&0\end{pmatrix},
\]
which imposes parameter mean zero on the densities, is invertible. By the coercivity of \cref{lem:bem-Vreg}(a), this holds for every polygon $X\in\mathfrak C$, in particular for every $X\in\Nd$, and the inverses are bounded as in \eqref{eq:bem-inv}. Let $\tilde K_X$, with $(\tilde K_X)_{pj}=\langle\mathcal K_X\varphi_j,\iota_p\rangle$, be the Galerkin matrix of
$\mathcal K_X$ on the fixed spaces $\mathcal P_0$ and $\mathcal P_1$, let $\tilde B$, with $\tilde B_{pj}=\langle\varphi_j,\iota_p\rangle$, be the pairing of densities and
traces, and let $\mathsf D_\xi$, with $(\mathsf D_\xi u)_p=(u_{p+1}-u_p)/h$, be the parameter difference matrix from traces to
densities. Then, with $\tilde Q_{\rm int}=\tfrac12\tilde B+\tilde K_X$ and
$\tilde Q_{\rm ext}=\tfrac12\tilde B-\tilde K_X$, we have
\[
\Lh(X)=P\big(2\mathsf D_\xi^\T\tilde V_X\mathsf D_\xi+\tilde Q_{\rm int}^\T\tilde V_X^{-1}\tilde Q_{\rm int}+\tilde Q_{\rm ext}^\T\tilde C_X\tilde Q_{\rm ext}\big)P,
\]
where $\tilde C_X$ is the bordered inverse of $\tilde V_X$ on densities of
parameter mean zero, that is, the upper left $N\times N$ block of the inverse of the bordered matrix above. The matrices $\tilde B$, $\mathsf D_\xi$ and $P$ do not
depend on $X$.
\end{proposition}

Thus, the shape dependence of $\Lh$ is the dependence of the parameter
integrals in $\tilde V_X$ and $\tilde K_X$ on the chord
$\Delta X(\xi,\eta)=X(\xi)-X(\eta)$. Both kernels are derived from its
logarithm. The single layer is
\begin{equation}\label{eq:logform}
\mathcal L_X\varsigma(\xi)=-\frac1{2\pi}\int_\Tor\mathrm{Re}\log\Delta X(\xi,\eta)\,\varsigma(\eta)\,\dd\eta ,
\end{equation}
and the kernel of $\mathcal K_X$ is $-\frac1{2\pi}$ times the imaginary part of the logarithmic derivative in $\eta$ of the chord, $\partial_\eta\Delta X(\eta,\xi)/\Delta X(\eta,\xi)=X'(\eta)/(X(\eta)-X(\xi))$. For a
direction $e$, let
\[
q_e(\xi,\eta)=\frac{e(\xi)-e(\eta)}{X(\xi)-X(\eta)}
\]
be the divided-difference quotient of $e$ relative to $X$.

The logarithm makes the shape derivatives explicit. Moving the polygon to
$X+\varepsilon e$ multiplies the chord by $1+\varepsilon q_e$, so the first
and the second shape derivatives of $\log\Delta X$ are $q_e$ and $-q_e^2$.
For the single layer, this gives the kernels of the shape derivatives
directly. For the double layer, an integration by parts on each panel moves
the derivative in $\eta$ from the kernel to the trace.

\begin{lemma}[The log form]\label{lem:bem-log}
Let $X$ be a closed polygon with a positive chord constant, and let $e$ be a
piecewise linear direction. For $\varsigma\in L^1$ and a P1 function $a$, let
\begin{gather*}
D\mathcal L_X[e]\varsigma=-\frac1{2\pi}\int_\Tor\mathrm{Re}(q_e)\,\varsigma\,\dd\eta,\qquad
D^2\mathcal L_X[e,e]\varsigma=\frac1{2\pi}\int_\Tor\mathrm{Re}(q_e^2)\,\varsigma\,\dd\eta,\\
D\mathcal K_X[e]a=\frac1{2\pi}\int_\Tor\mathrm{Im}(q_e)\,a'\,\dd\eta,\qquad
D^2\mathcal K_X[e,e]a=-\frac1{2\pi}\int_\Tor\mathrm{Im}(q_e^2)\,a'\,\dd\eta .
\end{gather*}
Then, the Galerkin matrices $\tilde V_Y$ and $\tilde K_Y$ of \cref{lem:bem-reduction} are twice differentiable functions of the vertices of the polygon $Y$ at $Y=X$, and their first and second derivatives in the direction $e$ are the Galerkin matrices of these operators,
\[
\begin{aligned}
\big(D\tilde V_X[e]\big)_{pq}&=\langle D\mathcal L_X[e]\iota_q,\iota_p\rangle,&\qquad\big(D^2\tilde V_X[e,e]\big)_{pq}&=\langle D^2\mathcal L_X[e,e]\iota_q,\iota_p\rangle,\\
\big(D\tilde K_X[e]\big)_{pj}&=\langle D\mathcal K_X[e]\varphi_j,\iota_p\rangle,&\qquad\big(D^2\tilde K_X[e,e]\big)_{pj}&=\langle D^2\mathcal K_X[e,e]\varphi_j,\iota_p\rangle .
\end{aligned}
\]
\end{lemma}

A maximal function controls $q_e$. Up to a factor controlled by the chord constant, $q_e(\xi,\eta)$ is
the average of $e'$ over the parameter interval between $\eta$ and $\xi$, so it is
bounded by the largest such average around $\xi$, and this maximal function is
bounded in $L^2$ by $\|e'\|_{L^2}$. Let
\[
Me'(\xi)=\sup_{0<r\le1/2}\frac1{2r}\int_{\xi-r}^{\xi+r}|e'(\eta)|\,\dd\eta
\]
be the periodic Hardy--Littlewood maximal function of $e'$. Since
$|e(\xi)-e(\eta)|\le2\,d(\xi,\eta)\,Me'(\xi)$, and symmetrically in $\eta$, every polygon $X\in\mathfrak C$ and every piecewise linear direction $e$ satisfy
\begin{equation}\label{eq:maximal}
|q_e(\xi,\eta)|\le\frac2{c_{\rm ca}}\min\big(Me'(\xi),Me'(\eta)\big).
\end{equation}

\subsection{Shape perturbation in the weak norm}
In the following two theorems, $f$ and $\varphi$ are piecewise linear and
are measured by the norms of these functions on the parameter circle, and
$z_f$, $z_\varphi$ and $\bXh$ are as in \cref{sec:setting}.

\begin{theorem}[{Maximal-function bound of the first variation}]\label{thm:bem-c1p}
Assume \ref{K1}, \ref{K2} with $r=1$, and \ref{K3}. Let $X\in\Nd$ with $h\le h_0$, and let $e$, $f$ and $\varphi$ be piecewise linear. Then,
\begin{equation}\label{eq:C1}
\big|z_\varphi^\T D_X\Lh(X)[e]\,z_f\big|\le C\,\|e'\|_{L^2}\,\|f\|_{L^2}\,\|\varphi\|_{L^2}.
\end{equation}
If moreover the segment $[X,X+e]$ lies in $\Nd$, then
\[
\big|b_{X+e,h}(f,\varphi)-\bXh(f,\varphi)\big|\le C\,\|e'\|_{L^2}\big(\|f\|_{L^\infty}\|\varphi\|_{L^2}+\|f\|_{L^2}\|\varphi\|_{L^\infty}\big).
\]
\end{theorem}

By \eqref{eq:maximal}, the quotient $q_e$ is bounded pointwise by the maximal function $Me'$, so the first variation of $\Lh$ is bounded as a bilinear form on $H^1$ by a multiple of $\|e'\|_{L^2}$ (\cref{cor:bem-c1p-ab}). The discrete regularity of \cref{lem:bem-areg} bounds the vectors $z_f$ and $z_\varphi$ in $H^1$ by $C\|f\|_{L^2}$ and $C\|\varphi\|_{L^2}$. The bound for the increment follows from the mean value inequality along the segment $[X,X+e]$, and only this bound requires the segment to lie in $\Nd$. The neighborhood $\Nd$, defined in $W^{1,\infty}$, keeps a uniform chord constant for the polygons and the single layer uniformly coercive, and within it the direction $e$ appears in both bounds only through $\|e'\|_{L^2}$, so no inverse estimate converting a factor $\|e'\|_{L^\infty}$ into $\|e'\|_{L^2}$ is needed.

\begin{theorem}[{Second variation and the Calder\'on commutator}]\label{thm:bem-c3p}
Assume \ref{K1}, \ref{K2} with $r=1$, and \ref{K3}. Let $X$, $e$, $f$ and $\varphi$ be piecewise linear with $h\le h_0$, and let the segment
$[X,X+e]$ lie in $\Nd$. Then, the second-order Taylor remainder of
$X\mapsto\bXh(f,\varphi)$ in the direction $e$,
\[
\delta^2\bXh(f,\varphi)=b_{X+e,h}(f,\varphi)-\bXh(f,\varphi)-D_X\bXh(f,\varphi)[e],
\]
satisfies
\begin{equation}\label{eq:C3}
|\delta^2\bXh(f,\varphi)|\le C\,\|e'\|_{L^\infty}\,\|e'\|_{L^2}\,\|f\|_{L^\infty}\,\|\varphi\|_{L^2}.
\end{equation}
\end{theorem}

The terms with two derivatives
of the kernels are bounded as in \cref{thm:bem-c1p}. The cross terms
$\Lhp D\Lh[e]\Lhp D\Lh[e]\Lhp$ require that $D\Lh[e]w$ be the load of an
$L^2$ function for traces $w$, which is a statement one derivative higher.
Since each factor $\|e'\|_{L^\infty}$ in place of $\|e'\|_{L^2}$ costs
$h^{-1/2}$ by the inverse estimate, the $L^\infty$ factor must fall on one
of the two derivatives and the $L^2$ factor on the other.
Differentiating the logarithmic kernel once turns it into a Cauchy kernel,
whose truncations are bounded on chord-arc curves \cite[Thm.~6, (66)]{David1984}
and depend Lipschitz-continuously on the curve, uniformly in the truncation
(\cref{prop:bem-david,lem:bem-holo}). For the flat curve
$X(\xi)=\xi$, the first variation of the double layer is, up to a Hilbert
transform, the first commutator of Calder\'on, with kernel
$(e(\xi)-e(\eta))/(\xi-\eta)^2$. Our bound is its boundedness relative to a
chord-arc background, and for a flat background it is Calder\'on's theorem
\cite{Calderon1965} (see also \cite[Thm.~3]{CoifmanMeyer1978}).

\subsection{Consistency with the continuum operators}\label{sec:bem-toolbox}
This subsection compares the matrix with the continuum operators. \Cref{thm:bem-c1} is the Lipschitz bound in the strong norm $\|e'\|_{L^\infty}$ that the proof of \cref{lem:bem-c4} uses.

\begin{theorem}[{Lipschitz and mixed derivatives of $\Lh$}]\label{thm:bem-c1}
Assume \ref{K1} and \ref{K3}. The result depends on \ref{K2} only through the radius $r_1$ of $\mathfrak C$, which depends on the interpolation constant $C_I$. Let $X,Y\in\mathfrak C$, let $e$ and $z$ be Lipschitz directions, and let $v$ and $w$ be real Lipschitz functions on $\Tor$. Then, the function $\varepsilon\mapsto\mathfrak s_{X+\varepsilon_1e+\varepsilon_2z,h}(v,w)$ of $\varepsilon\in\R^2$, where $\mathfrak s_{X,h}$ is the Schur form \eqref{eq:bem-ext}, is twice continuously differentiable near $\varepsilon=0$. Let $D\mathfrak s_{X,h}[e](v,w)$ and $D^2\mathfrak s_{X,h}[e,z](v,w)$ be its derivatives $\partial_{\varepsilon_1}$ and $\partial_{\varepsilon_1}\partial_{\varepsilon_2}$ at $\varepsilon=0$. We have
\[
\begin{aligned}
|D\mathfrak s_{X,h}[e](v,w)|&\le C\|e'\|_{L^\infty}\|v\|_{H^{1/2}}\|w\|_{H^{1/2}},\\
|D^2\mathfrak s_{X,h}[e,z](v,w)|&\le C\|e'\|_{L^\infty}\|z'\|_{L^\infty}\|v\|_{H^{1/2}}\|w\|_{H^{1/2}},\\
|\mathfrak s_{X,h}(v,w)-\mathfrak s_{Y,h}(v,w)|&\le C\|X'-Y'\|_{L^\infty}\|v\|_{H^{1/2}}\|w\|_{H^{1/2}},\\
\big|\mathfrak s_{Y,h}(v,w)-\mathfrak s_{X,h}(v,w)-D\mathfrak s_{X,h}[Y-X](v,w)\big|&\le C\|X'-Y'\|_{L^\infty}^2\|v\|_{H^{1/2}}\|w\|_{H^{1/2}}.
\end{aligned}
\]
\end{theorem}

For a polygon $X\in\mathfrak C$ and P1 functions $v$ and $w$, we have $\mathfrak s_{X,h}(v,w)=w^\T\Lh(X)v$ (\cref{sec:bem-unif}). The theorem therefore bounds $\Lh$ and its shape derivatives as forms on $H^{1/2}$, and for curves that are not polygons and traces that are not P1 we understand $\Lh$ through $\mathfrak s_{X,h}$.

\Cref{lem:bem-bweak} is the consistency on the interpolating polygon. Its order is $\epsh=h^{3/2}$,
the order of the interpolation and Galerkin errors in $H^{1/2}$. In the weaker norm of $(H^1)^*$, the order is $h^2$ under \ref{K2} with $r=3$.

\begin{theorem}[{Consistency on the interpolating polygon}]\label{lem:bem-bweak}
Assume \ref{K1} and \ref{K3}. Let $X=\Ih x$ with $h\le h_0$, and let $\mathsf r=\Lh(X)\Ih g-\Mv\Ih(\Lam_\Gamma g)$. If \ref{K2} holds with $r=2$, then
\[
\dnorm{\mathsf r}\le Ch\|g\|_{H^2}\qquad(g\in H^2),
\]
and if \ref{K2} holds with $r=3$, then
\[
\dnorm{\mathsf r}\le C\epsh\|g\|_{H^3},\qquad\|\mathsf r\|_{(H^1)^*}\le Ch^2\|g\|_{H^3}\qquad(g\in H^3).
\]
\end{theorem}

We do not know whether the order $\epsh=h^{3/2}$ of the bound in $\dnorm{\cdot}$ under \ref{K2} with $r=3$ is optimal.

The comparison of the polygon with the curve, \cref{lem:bem-geo}, is a
boundary-approximation estimate of the type of \cite{Nedelec1976} and
\cite[Ch.~8]{SauterSchwab2011}. It has the order $\epsh$ without a
logarithmic loss, because the leading part of the difference of the
single-layer kernels is the commutator of $\theta=\Ih x-x$ with the Cauchy
integral, which is bounded on $L^2$ by $C\|\theta\|_{L^\infty}\le Ch^2$,
whereas a bound of the kernels in absolute value would give only
$h^2|\log h|$.

\Cref{lem:bem-c4} compares the shape derivative of $\Lh$ with that of the continuum map, tested with P1 functions. Let $e$ be a P1 direction, let $g\in H^2$, and let $x_\varepsilon=x+\varepsilon e$, which lies in $\mathfrak C$ for $|\varepsilon|\,\|e'\|_{L^\infty}\le2r_1$ and is piecewise smooth. Since $g$ is continuously differentiable, it is Lipschitz. By \cref{prop:bem-dtn}, we have
\begin{equation}\label{eq:bem-dtn-schur}
\int_\Tor(\Lam_{x_\varepsilon}g)\,\varphi_j\,|x_\varepsilon'|\,\dd\xi=\langle\mathcal S_{x_\varepsilon}g,\varphi_j\rangle ,
\end{equation}
where $\langle\mathcal S_Yv,w\rangle$ denotes the Schur form \eqref{eq:bem-ext} with the exact inverses $\mathcal L_Y^{-1}$ and $\mathcal L_{Y,0}^{-1}$ in place of the Galerkin inverses.

For a Lipschitz direction $e$, functions $\varsigma,\psi\in L^2$ and a real Lipschitz function $a$, the functions $\varepsilon\mapsto\langle\mathcal L_{x+\varepsilon e}\varsigma,\psi\rangle$ and $\varepsilon\mapsto\langle\mathcal K_{x+\varepsilon e}a,\psi\rangle$ are differentiable at $\varepsilon=0$ by \cref{lem:bem-cont}(iii), and we denote their derivatives at $\varepsilon=0$ by $D\mathcal L_x[e](\varsigma,\psi)$ and $D\mathcal K_x[e](a,\psi)$. For a P1 direction $e$, the derivative $D\mathcal L_x[e](\varsigma,\psi)$ is linear in $\varsigma$ and in $\psi$, and $D\mathcal K_x[e](a,\psi)$ is linear in $\psi$. \Cref{lem:bem-dS} shows that $\varepsilon\mapsto\langle\mathcal S_{x_\varepsilon}v,w\rangle$ is differentiable at $\varepsilon=0$ for real Lipschitz functions $v$ and $w$ and expresses its derivative through $D\mathcal L_x[e]$ and $D\mathcal K_x[e]$. With $v=g$ and $w=\varphi_j$, both sides of \eqref{eq:bem-dtn-schur} are therefore differentiable at $\varepsilon=0$.

\begin{theorem}[{Consistency of the shape derivative}]\label{lem:bem-c4}
Assume \ref{K1}, \ref{K2} with $r=2$, and \ref{K3}. Let $X=\Ih x$ with $h\le h_0$, let $g\in H^2$, and let $e$ be a piecewise linear direction. For each $j$, the map $\varepsilon\mapsto\int_\Tor(\Lam_{x+\varepsilon e}\,g)\,\varphi_j\,|x'+\varepsilon e'|\,\dd\xi$, where $\Lam_{x+\varepsilon e}$ is the two-phase Dirichlet-to-Neumann map of the Lipschitz curve $x+\varepsilon e$, is differentiable at $\varepsilon=0$ by \eqref{eq:bem-dtn-schur} and \cref{lem:bem-dS}. Let
\[
\mathsf r_j=\big(D\Lh(X)[e]\,\Ih g\big)_j-\frac{\dd}{\dd\varepsilon}\int_\Tor(\Lam_{x+\varepsilon e}\,g)\,\varphi_j\,|x'+\varepsilon e'|\,\dd\xi\Big|_{\varepsilon=0}.
\]
Then,
\[
\|\mathsf r\|_{(H^1)^*}\le Ch\,\|e'\|_{L^\infty}\|g\|_{H^2},\qquad
\dnorm{\mathsf r}\le Ch^{1/2}\,\|e'\|_{L^\infty}\|g\|_{H^2}.
\]
\end{theorem}

The proofs also use classical boundary element facts, such as the coercivity of the single layer operator, the layer representation of the Dirichlet-to-Neumann map and the discrete regularity of $\Lhp$, in versions whose constants depend only on the constants of $\mathscr K$. They are proved, or cited with their sources, in \cref{sec:bem-unif,sec:bem-a}.

\appendix
\section{Proofs of the boundary element estimates}\label{sec:app-bem}

The notation is that of \cref{sec:setting,sec:shape}, and $x\in\mathscr K$ is a fixed curve of the family of \cref{ass:K}. Every result below uses \ref{K1} and \ref{K3}. A result stated on $\Nd$ also uses \ref{K2} with $r=1$, through the constant $C_I$ fixed before \cref{def:nbhd}, and higher orders of \ref{K2} are named where they are used. The constants $C$ and $c$ are as in \cref{sec:setting}.

For the parameter difference quotient $\mathsf D_\xi:\mathcal P_1\to\mathcal P_0$ of \cref{lem:bem-reduction}, we have $\mDs=\mG_X^{-1}\mathsf D_\xi$ with $\mG_X=N\operatorname{diag}(|\Gamma_p|)$. For $X\in\mathfrak C$ the edge lengths satisfy $c_{\rm ca}h\le|\Gamma_p|\le M_1h$, so the grid is quasi-uniform in arclength, and the approximation properties and inverse estimates of $\mathcal P_0$ and $\mathcal P_1$ on the uniform parameter grid hold with constants independent of $h$.

\subsection{The parameter-circle reduction and the log form}

\begin{proof}[Proof of \cref{lem:bem-reduction}]
We write $\tilde V=\tilde V_X$ and $\tilde K=\tilde K_X$, with the entries of \cref{lem:bem-reduction}. Arclength measure is $|X'|\,\dd\xi$, and $|X'|=N|\Gamma_p|$ on panel $p$. In the single-layer matrix, both variables carry arclength measure. In the pairing matrix, only the test variable carries arclength measure, since there is no trial integral. In the double-layer matrix, only the test variable does as well, because the trial measure $|X'(\eta)|\,\dd\eta$ is part of the parameter kernel of $\mathcal K_X$. Maue's formula writes the hypersingular matrix as the single layer between arclength derivatives, and $\partial_s=|X'|^{-1}\partial_\xi$ on $\mathcal P_1$. Hence, we have
\[
\mV=\mG_X\tilde V\mG_X,\qquad\mB=\mG_X\tilde B,\qquad\mK=\mG_X\tilde K,\qquad\mW=\mathsf D_\xi^\T\tilde V\mathsf D_\xi .
\]
Put $\tilde Q_{\rm int}=\tfrac12\tilde B+\tilde K$ and $\tilde Q_{\rm ext}=\tfrac12\tilde B-\tilde K$. In the interior Schur complement, the factors $\mG_X$ cancel, since
\[
Q_{\rm int}^\T \mV^{-1}Q_{\rm int}=\tilde Q_{\rm int}^\T \mG_X(\mG_X\tilde V\mG_X)^{-1}\mG_X\tilde Q_{\rm int}=\tilde Q_{\rm int}^\T\tilde V^{-1}\tilde Q_{\rm int}.
\]
The bordering vector is $\mB\one=\mG_X\tilde B\one=h\,\mG_X\one$, so the bordered matrix of $\mV$ equals $\operatorname{diag}(\mG_X,1)$ times the bordered matrix of $\tilde V$ of \cref{lem:bem-reduction} times $\operatorname{diag}(\mG_X,1)$. Hence, $\mathsf C=\mG_X^{-1}\tilde C\mG_X^{-1}$, and $Q_{\rm ext}^\T\mathsf CQ_{\rm ext}=\tilde Q_{\rm ext}^\T\tilde C\tilde Q_{\rm ext}$. Together with the interior term and $\mW=\mathsf D_\xi^\T\tilde V\mathsf D_\xi$, we obtain
\[
\Lh(X)=P\big(2\mathsf D_\xi^\T\tilde V\mathsf D_\xi+\tilde Q_{\rm int}^\T\tilde V^{-1}\tilde Q_{\rm int}+\tilde Q_{\rm ext}^\T\tilde C\tilde Q_{\rm ext}\big)P,
\]
where $\tilde C$ is the bordered inverse of $\tilde V$ on $\{\varsigma\in\mathcal P_0\mid\int\varsigma\,\dd\xi=0\}$, and $P=I-\one\one^\T/N$ is the Euclidean deflation.
\end{proof}

\begin{proof}[Proof of \cref{lem:bem-log}]
Let $c>0$ be a chord constant of $X$, and let $\mathcal U$ be the neighborhood of $X$ in the space of vertices that consists of the polygons $Y$ with $\|Y'-X'\|_{L^\infty}<c/2$. Every $Y\in\mathcal U$ satisfies $|\Delta Y(\xi,\eta)|\ge\frac c2\,d(\xi,\eta)$. For directions $e$ and $z$, let $q^Y_e=\Delta e/\Delta Y$ and $q^Y_z=\Delta z/\Delta Y$. Since $|\Delta e(\xi,\eta)|\le\|e'\|_{L^\infty}d(\xi,\eta)$, we have $|q^Y_e|\le C\|e'\|_{L^\infty}$, and likewise for $z$, uniformly in $Y\in\mathcal U$.

For $\xi\ne\eta$, the chord $\Delta Y(\xi,\eta)$ is linear in the vertices of $Y$ and does not vanish on $\mathcal U$. Hence, $\log|\Delta Y(\xi,\eta)|$ is a smooth function of $Y\in\mathcal U$, and
\[
D_Y\log|\Delta Y|[e]=\mathrm{Re}\,q^Y_e,\qquad D_Y^2\log|\Delta Y|[e,z]=-\mathrm{Re}\big(q^Y_eq^Y_z\big).
\]
Both derivatives are bounded by $C\|e'\|_{L^\infty}$ and $C\|e'\|_{L^\infty}\|z'\|_{L^\infty}$, uniformly in $(\xi,\eta)$ and $Y\in\mathcal U$, and $\|e'\|_{L^\infty}$ is bounded by a multiple of the vertex norm of $e$. The entries $(\tilde V_Y)_{pq}=-\frac1{2\pi}\int_p\int_q\log|\Delta Y|\,\dd\eta\,\dd\xi$ have integrable integrands, so these bounds allow us to differentiate under the integral sign twice. Hence, $Y\mapsto\tilde V_Y$ is twice differentiable on $\mathcal U$, and the entries of $D\tilde V_X[e]$ and of $D^2\tilde V_X[e,e]$ are the integrals over $p\times q$ of $-\frac1{2\pi}\mathrm{Re}\,q_e$ and of $\frac1{2\pi}\mathrm{Re}(q_e^2)$, as stated.

For the double layer, let $u$ be a P1 function, fix a panel $p$ and a point $\xi$ in its interior, and let $Y\in\mathcal U$. For $\eta$ on the panel $p$, the vector $Y(\eta)-Y(\xi)$ is a real multiple of $Y'(\eta)$, so the kernel of $\mathcal K_Y$ vanishes there. On the panels $q\ne p$, which have positive distance from $\xi$, the kernel is bounded and smooth in $Y$. Hence,
\[
\mathcal K_Yu(\xi)=-\frac1{2\pi}\sum_{q\ne p}\int_{\xi_q}^{\xi_{q+1}}\mathrm{Im}\frac{Y'(\eta)}{Y(\eta)-Y(\xi)}\,u(\eta)\,\dd\eta .
\]
Since $Y(\eta)-Y(\xi)$ is linear in $Y$, the derivative of the kernel in the direction $e$ is
\[
D_Y\,\mathrm{Im}\frac{Y'(\eta)}{Y(\eta)-Y(\xi)}[e]=\mathrm{Im}\frac{e'(\eta)\big(Y(\eta)-Y(\xi)\big)-Y'(\eta)\big(e(\eta)-e(\xi)\big)}{\big(Y(\eta)-Y(\xi)\big)^2}=\partial_\eta\,\mathrm{Im}\,q^Y_e(\xi,\eta).
\]
On each panel $q\ne p$, the function $q^Y_e(\xi,\cdot)$ is smooth up to the endpoints. Differentiating under the integral sign and integrating by parts on these panels, we obtain
\[
D_Y\big(\mathcal K_Yu(\xi)\big)[e]=-\frac1{2\pi}\sum_{q\ne p}\Big(\big[\mathrm{Im}\,q^Y_e(\xi,\cdot)\,u\big]_{\xi_q}^{\xi_{q+1}}-\int_{\xi_q}^{\xi_{q+1}}\mathrm{Im}\,q^Y_e(\xi,\eta)\,u'(\eta)\,\dd\eta\Big).
\]
Since $q^Y_e(\xi,\cdot)$ and $u$ are continuous at the vertices, the boundary terms telescope around the polygon to $-\big[\mathrm{Im}\,q^Y_e(\xi,\cdot)\,u\big]_{\xi_p}^{\xi_{p+1}}$. On the panel $p$, the quotient $q^Y_e(\xi,\eta)=(e_{p+1}-e_p)/(Y_{p+1}-Y_p)$ does not depend on $\eta$, and $u$ is linear. Hence, this sum equals $-\int_{\xi_p}^{\xi_{p+1}}\mathrm{Im}\,q^Y_e\,u'\,\dd\eta$, and we obtain
\[
D_Y\big(\mathcal K_Yu(\xi)\big)[e]=\frac1{2\pi}\int_\Tor\mathrm{Im}\,q^Y_e(\xi,\eta)\,u'(\eta)\,\dd\eta .
\]
The integrand is bounded by $C\|e'\|_{L^\infty}|u'(\eta)|$, and it is smooth in $Y$, with the derivative $D_Y\,\mathrm{Im}\,q^Y_e[z]=-\mathrm{Im}(q^Y_eq^Y_z)$ bounded by $C\|e'\|_{L^\infty}\|z'\|_{L^\infty}$, uniformly in $\xi$ and $Y\in\mathcal U$. Since $(\tilde K_Yu)_p=\int_{\xi_p}^{\xi_{p+1}}\mathcal K_Yu\,\dd\xi$, differentiation under the integral sign in $\xi$, as for the single layer, gives the derivatives of $\tilde K_Y$ at $X$, with the entries as stated.
\end{proof}

\subsection{Maximal-function bounds}

The maximal functions are bounded on $L^2$ by the maximal theorem \cite[Ch.~I]{Stein1970}, which we use on $\R$ and for the exponent $2$.

\begin{proposition}[Maximal theorem]\label{prop:bem-hl}
For $f\in L^2(\R;\C)$ and $t\in\R$, let
\[
\mathcal M_\R f(t)=\sup_{r>0}\frac1{2r}\int_{t-r}^{t+r}|f(s)|\,\dd s .
\]
There is a constant $A$ such that $\|\mathcal M_\R f\|_{L^2(\R)}\le A\|f\|_{L^2(\R)}$ for every $f\in L^2(\R;\C)$.
\end{proposition}

Every interval in the definition of $Me'(\xi)$, $0\le\xi<1$, lies in $[-\tfrac12,\tfrac32]$. \Cref{prop:bem-hl}, applied to the periodic extension of $e'$ restricted to $[-\tfrac12,\tfrac32]$, therefore gives $\|Me'\|_{L^2}\le C\|e'\|_{L^2}$. Let $X\in\mathfrak C$ be a polygon and $e$ a piecewise linear direction. By the kernels of \cref{lem:bem-log} and \eqref{eq:maximal}, we have $|D\mathcal L_X[e]\varsigma(\xi)|\le C\,Me'(\xi)\|\varsigma\|_{L^1}$ and $|D\mathcal K_X[e]u(\xi)|\le C\,Me'(\xi)\|u'\|_{L^1}$, and the Cauchy--Schwarz inequality gives
\begin{equation}\label{eq:bem-tri}
|\langle D\mathcal L_X[e]\varsigma,\psi\rangle|\le C\|e'\|_{L^2}\|\varsigma\|_{L^1}\|\psi\|_{L^2},\qquad
|\langle D\mathcal K_X[e]u,\psi\rangle|\le C\|e'\|_{L^2}\|u'\|_{L^1}\|\psi\|_{L^2}.
\end{equation}

\subsection{The Cauchy integral}\label{sec:bem-cauchy}

For $M\ge c>0$, we write $\mathfrak F(M,c)$ for the set of closed Lipschitz curves $Y:\Tor\to\C$ with $|Y'|\le M$ almost everywhere and chord constant $c$, so that $|Y'|\ge c$ almost everywhere by \cref{sec:setting}. The class $\mathfrak F(M_1,c_{\rm ca})$ contains $\mathfrak C$. For $X\in\mathfrak F(M,c)$, $\delta>0$ and $f\in L^2(\Tor;\C)$, the truncated Cauchy integral is
\[
\mathcal C^\delta_Xf(\xi)=\frac1{\pi\iu}\int_{d(\xi,\eta)>\delta}\frac{f(\eta)X'(\eta)}{X(\eta)-X(\xi)}\,\dd\eta .
\]
On the region of integration the kernel is bounded by $M|f(\eta)|/(c\delta)$, so the integral converges absolutely. Since $d(\xi,\eta)\le\tfrac12$, the truncations with $\delta\ge\tfrac12$ vanish. Let $\mathcal K^\delta_Xa$ be defined as $\mathcal K_Xa$ in \cref{sec:setting}, with the principal value replaced by the integral over $\{\eta\in\Tor\mid d(\xi,\eta)>\delta\}$. Since $\operatorname{Re}(w/(\pi\iu))=\operatorname{Im}w/\pi$, we have, for real $a$,
\begin{equation}\label{eq:bem-KC}
\mathcal K^\delta_Xa=-\tfrac12\operatorname{Re}\mathcal C^\delta_Xa .
\end{equation}

The truncations are bounded uniformly by David's theorem \cite[Thm.~6, (66)]{David1984}, which we state for the kernel $1/z$ and the exponent $2$, with the function $\bar\mu$ of \cite[(10)]{David1984}.

\begin{proposition}[David's theorem]\label{prop:bem-david}
There are constants $C$ and $\beta$ with the following property. Let $\Gamma$ be a rectifiable curve with arclength measure $\mu$, and let $\bar\mu(z)$ be the supremum of $\mu(D)/r(D)$ over the discs $D$ containing $z$, where $r(D)$ is the radius of $D$. Then, the maximal chordal truncation
\[
T^*_\mu g(z)=\sup_{\varepsilon>0}\Big|\int_{|z-w|\ge\varepsilon}\frac{g(w)}{z-w}\,\dd\mu(w)\Big|
\]
satisfies $\|T^*_\mu g\|_{L^2(\mu)}\le C\|(1+\bar\mu)^\beta g\|_{L^2(\mu)}$ for every continuous $g$ with compact support.
\end{proposition}

We transfer \cref{prop:bem-david} to the parameter truncations. Let $Y\in\mathfrak F(M,c)$ and $\Gamma=Y(\Tor)$. The chord constant makes $Y$ injective on a period, and the parameters of the points of $\Gamma$ in a disc of radius $r$ lie within distance $2r/c$ of any one of them. Hence, the disc has measure bounded by $4Mr/c$, and $\bar\mu\le4M/c$. For fixed $\varepsilon>0$, the chordal truncation at level $\varepsilon$ is therefore bounded on $L^2(\mu)$ by $C(1+4M/c)^\beta$, first on continuous functions with compact support and then, since its kernel is bounded by $1/\varepsilon$ and $\mu(\Gamma)$ is finite, on $L^2(\mu)$ by density. For $f\in L^2$, let $g$ be the function on $\Gamma$ with $g(Y(\eta))=f(\eta)Y'(\eta)/|Y'(\eta)|$. Since $\dd\mu=|Y'|\,\dd\eta$, the truncation of $g$ at level $\varepsilon$, evaluated at $Y(\xi)$, is $-\pi\iu$ times the chordal truncation
\[
\frac1{\pi\iu}\int_{|Y(\eta)-Y(\xi)|\ge\varepsilon}\frac{f(\eta)Y'(\eta)}{Y(\eta)-Y(\xi)}\,\dd\eta .
\]
Since $c\le|Y'|\le M$ almost everywhere, the norms of $L^2(\mu)$ and $L^2(\dd\xi)$ are equivalent, so the chordal truncations are bounded on $L^2(\dd\xi)$ by a constant that depends only on $M$ and $c$. Since the parameter truncations with $\delta\ge\tfrac12$ vanish, let $\delta<\tfrac12$. The chord constant shows that $\{d(\xi,\eta)>\delta\}$ lies in $\{|Y(\eta)-Y(\xi)|\ge c\delta\}$. The parameter truncation at $\delta$ and the chordal truncation at the level $c\delta$ therefore differ by the integral over a part of $\{d(\xi,\eta)\le\delta\}$ on which the kernel is bounded by $M|f(\eta)|/(c\delta)$, and the difference is bounded by $C(M/c)\,\mathcal Mf(\xi)$. Here $\mathcal Mf$ is the periodic Hardy--Littlewood maximal function of $f$, defined as $Me'$ in \cref{sec:shape}, and \cref{prop:bem-hl} bounds it on $L^2$ as for $Me'$. Hence, there is a constant $C_D$, depending only on $M$ and $c$, such that
\begin{equation}\label{eq:bem-david}
\|\mathcal C^\delta_Yf\|_{L^2}\le C_D\|f\|_{L^2}\qquad(Y\in\mathfrak F(M,c),\ \delta>0,\ f\in L^2).
\end{equation}

The truncated operators depend on the curve in a Lipschitz way, uniformly in the truncation. We prove this for the Cauchy integrals, which depend holomorphically on the curve $X+z(Y-X)$ for complex $z$, and differentiate the logarithmic kernel before using it, because $\log|\Delta X+z\Delta(Y-X)|$ is not holomorphic in $z$. The proof of the following lemma uses the bound \eqref{eq:bem-david}, which rests on the boundedness of the Cauchy integral on Lipschitz graphs \cite{CMM1982} and its extension to chord-arc curves \cite{David1984}.

\begin{lemma}[Lipschitz dependence of the truncated Cauchy integrals]\label{lem:bem-holo}
Let $M\ge c>0$. There is a constant $C_\Delta$, depending only on $M$ and $c$, such that for $X,Y\in\mathfrak F(M,c)$, $\delta>0$ and $f\in L^2$,
\[
\|\mathcal C^\delta_Xf-\mathcal C^\delta_Yf\|_{L^2}\le C_\Delta\|X'-Y'\|_{L^\infty}\|f\|_{L^2}.
\]
\end{lemma}

\begin{proof}
Let $s=\|X'-Y'\|_{L^\infty}$, and let $C_1$ and $C_2$ be the constants of \eqref{eq:bem-david} for $\mathfrak F(M,c)$ and for $\mathfrak F(M+c/2,c/2)$. Since the derivative of $Y-X$ is bounded by $s$, we have
\begin{equation}\label{eq:bem-dXY}
|(Y-X)(\xi)-(Y-X)(\eta)|\le s\,d(\xi,\eta).
\end{equation}
If $s=0$, then $Y-X$ is constant, and the truncations of $X$ and $Y$ coincide. If $s>c/4$, then \eqref{eq:bem-david} gives $\|\mathcal C^\delta_Xf-\mathcal C^\delta_Yf\|_{L^2}\le2C_1\|f\|_{L^2}\le(8C_1/c)\,s\|f\|_{L^2}$.

Let $0<s\le c/4$, fix $\delta$ and $f$, and put $u=\mathcal C^\delta_Yf-\mathcal C^\delta_Xf$. For $z\in\C$, let $W_z=X+\frac{cz}{4s}(Y-X)$, so that $W_0=X$ and $W_{z_0}=Y$ for $z_0=4s/c\in(0,1]$. By \eqref{eq:bem-dXY} and the chord constant of $X$, the chord constant of $W_z$ is at least $c\,(1-|z|/4)$. We put
\[
\Phi(z)=\int_\Tor\overline{u(\xi)}\,\mathcal C^\delta_{W_z}f(\xi)\,\dd\xi .
\]
For $|z|<3$, the chord constant of $W_z$ is at least $c/4$, so $|W_z(\eta)-W_z(\xi)|\ge c\delta/4$ on $\{d(\xi,\eta)>\delta\}$. There, the integrand $\overline{u(\xi)}f(\eta)W_z'(\eta)/(W_z(\eta)-W_z(\xi))$ is holomorphic in $z$, and its derivative in $z$ is bounded by $C|u(\xi)||f(\eta)|/\delta$, which is integrable on $\Tor\times\Tor$. Hence, $\Phi$ is holomorphic in $|z|<3$ by differentiation under the integral sign.

For $|z|\le2$, we have $|W_z'|\le M+c/2$ almost everywhere and the chord constant of $W_z$ is at least $c/2$, so $W_z\in\mathfrak F(M+c/2,c/2)$, and the Cauchy--Schwarz inequality and \eqref{eq:bem-david} give $|\Phi(z)|\le C_2\|u\|_{L^2}\|f\|_{L^2}$. For $|w|\le1$, Cauchy's estimate on the circle of radius $1$ around $w$ gives $|\Phi'(w)|\le C_2\|u\|_{L^2}\|f\|_{L^2}$. Since $\Phi(z_0)-\Phi(0)=\|u\|_{L^2}^2$ and $|z_0|\le1$, the mean value inequality on the unit disc gives
\[
\|u\|_{L^2}^2\le C_2\|u\|_{L^2}\|f\|_{L^2}\,\frac{4s}c .
\]
Hence, $\|u\|_{L^2}\le(4C_2/c)\,s\|f\|_{L^2}$, and the lemma holds with $C_\Delta=\max(4C_2,8C_1)/c$.
\end{proof}

The double layer is the limit of its truncations on Lipschitz functions. Let $X\in\mathfrak F(M,c)$, let $a$ be a real Lipschitz function, and let $\xi$ be a point at which $X$ is differentiable, so that $|X'(\xi)|\ge c$. The function $\operatorname{Im}\big(X'(\eta)/(X(\eta)-X(\xi))\big)\big(a(\eta)-a(\xi)\big)$ is bounded by $M\|a'\|_{L^\infty}/c$, so its truncated integrals converge as $\delta\to0$. Let $J_\delta(\xi)$ be the integral of $X'(\eta)/(X(\eta)-X(\xi))$ over the arc $\{\eta\in\Tor\mid d(\xi,\eta)>\delta\}$. It is the increment of a continuous logarithm of $\eta\mapsto X(\eta)-X(\xi)$ over this arc, so it depends continuously on $\delta$, and its exponential $(X(\xi-\delta)-X(\xi))/(X(\xi+\delta)-X(\xi))$ tends to $-1$. Hence, $\operatorname{Im}J_\delta(\xi)$ tends to $(2m(\xi)+1)\pi$ for an integer $m(\xi)$, which does not depend on $a$. Since $X$ is differentiable almost everywhere, the truncations $\mathcal K^\delta_Xa$ converge almost everywhere, and the limit is
\begin{equation}\label{eq:bem-pv}
\mathcal K_Xa(\xi)=-\frac{2m(\xi)+1}2\,a(\xi)-\frac1{2\pi}\int_\Tor\operatorname{Im}\frac{X'(\eta)}{X(\eta)-X(\xi)}\,\big(a(\eta)-a(\xi)\big)\,\dd\eta .
\end{equation}
By \eqref{eq:bem-KC}, \eqref{eq:bem-david} and Fatou's lemma, we obtain $\|\mathcal K_Xa\|_{L^2}\le\tfrac12C_D\|a\|_{L^2}$.

The derivatives of the layer operators are computed on the truncations.

\begin{lemma}[Integration by parts on the truncations]\label{lem:bem-trunc}
Let $M\ge c>0$, let $X\in\mathfrak F(M,c)$, and let $f$ and $\varphi$ be real periodic $C^1$ functions. For $0<\delta<\tfrac12$, we have
\begin{equation}\label{eq:bem-dsC}
\int_\Tor\mathcal C^\delta_Xf\,\varphi'\,\dd\xi=-\int_\Tor\varphi\,X'\,\mathcal C^\delta_X(f'/X')\,\dd\xi+\frac1{\pi\iu}R_\delta,
\end{equation}
where the remainder $R_\delta$, given in the proof, tends to $0$ as $\delta\to0$. For real $f\in L^2$, let $\mathcal D^\delta_Xf=\tfrac12\operatorname{Re}\big[\iu X'\,\mathcal C^\delta_X(f/X')\big]$. Then, we have
\begin{equation}\label{eq:bem-dsL}
\lim_{\delta\to0}\int_\Tor\mathcal D^\delta_Xf\,\varphi\,\dd\xi=-\int_\Tor\mathcal L_Xf\,\varphi'\,\dd\xi .
\end{equation}
\end{lemma}
\begin{proof}
Let $k(\xi,\eta)=(X(\eta)-X(\xi))^{-1}$. On $\{d(\xi,\eta)>\delta\}$, we have $|k|\le1/(c\delta)$, so Fubini's theorem applies to all the integrals below, and $\partial_\xi k=X'(\xi)k^2$ and $\partial_\eta k=-X'(\eta)k^2$ almost everywhere. We integrate by parts in $\xi$ for fixed $\eta$, exchange the order of integration, and integrate by parts in $\eta$ for fixed $\xi$, using $f(\eta)X'(\eta)k^2=-f(\eta)\partial_\eta k$, each time over the arc $\{d(\xi,\eta)>\delta\}$. Collecting the boundary terms of both steps, we obtain \eqref{eq:bem-dsC} with
\[
\begin{aligned}
R_\delta&=\int_\Tor X'(\xi)\Big(\frac{N_\delta^-(\xi)}{X(\xi-\delta)-X(\xi)}-\frac{N_\delta^+(\xi)}{X(\xi+\delta)-X(\xi)}\Big)\,\dd\xi,\\
N_\delta^\pm(\xi)&=\varphi(\xi)\big(f(\xi\pm\delta)-f(\xi)\big)-f(\xi)\big(\varphi(\xi\pm\delta)-\varphi(\xi)\big).
\end{aligned}
\]
Let $B$ bound $f$, $f'$, $\varphi$ and $\varphi'$. Since $|N_\delta^\pm|\le2B^2\delta$ and $|X(\xi\pm\delta)-X(\xi)|\ge c\delta$, the integrand of $R_\delta$ is bounded by $4MB^2/c$. At every point $\xi$ at which $X$ is differentiable, both quotients tend to $(\varphi f'-f\varphi')(\xi)/X'(\xi)$, since $X'(\xi)\ne0$. Hence, the integrand of $R_\delta$ tends to $0$ almost everywhere, and $R_\delta\to0$ by dominated convergence.

For \eqref{eq:bem-dsL}, let $L(\xi,\eta)=-\frac1{2\pi}\log|X(\xi)-X(\eta)|$ be the kernel of $\mathcal L_X$. Since $\int_\Tor|L(\xi,\eta)|\,\dd\xi$ is bounded uniformly in $\eta$, Fubini's theorem and dominated convergence show that $\int_\Tor\mathcal L_Xf\,\varphi'\,\dd\xi$ is the limit as $\delta\to0$ of the integrals of $f(\eta)L(\xi,\eta)\varphi'(\xi)$ over $\{d(\xi,\eta)>\delta\}$. For fixed $\eta$, the function $\xi\mapsto L(\xi,\eta)$ is Lipschitz on the arc $\{d(\xi,\eta)>\delta\}$, with $\partial_\xi L=-\frac1{2\pi}\operatorname{Re}\big(X'(\xi)/\Delta X\big)$ almost everywhere, and integration by parts gives
\[
\begin{aligned}
\int_{d(\xi,\eta)>\delta}L\varphi'\,\dd\xi&=b_\delta(\eta)-\int_{d(\xi,\eta)>\delta}\partial_\xi L\,\varphi\,\dd\xi,\\
b_\delta(\eta)&=L(\eta-\delta,\eta)\big(\varphi(\eta-\delta)-\varphi(\eta+\delta)\big)+\frac{\varphi(\eta+\delta)}{2\pi}\log\frac{|X(\eta+\delta)-X(\eta)|}{|X(\eta-\delta)-X(\eta)|}.
\end{aligned}
\]
The term with $\varphi(\eta-\delta)-\varphi(\eta+\delta)$ is bounded by $C\delta(1+|\log\delta|)\|\varphi'\|_{L^\infty}$, and the quotient in the logarithm lies in $[c/M,M/c]$ and tends to $1$ at every point of differentiability of $X$. Hence, $b_\delta$ is bounded uniformly in $\delta$ and tends to $0$ almost everywhere, so $\int f\,b_\delta\,\dd\eta\to0$ by dominated convergence. Since $\int_{d(\xi,\eta)>\delta}f(\eta)\,\partial_\xi L\,\dd\eta=\mathcal D^\delta_Xf(\xi)$, multiplying the identity above by $f(\eta)$, integrating in $\eta$, exchanging the order of integration in the integral of $\partial_\xi L\,\varphi$ and letting $\delta\to0$ gives \eqref{eq:bem-dsL}.
\end{proof}

\begin{lemma}[First derivatives of the layer operators]\label{lem:bem-H1}
Let $M\ge c>0$, let $X\in\mathfrak F(M,c)$, and let $C_D$ be the constant of \eqref{eq:bem-david}. For real $f\in L^2$, the function $\mathcal L_Xf$ lies in $H^1$, and
\begin{equation}\label{eq:bem-dsLH1}
\|\partial_\xi\mathcal L_Xf\|_{L^2}\le\frac{MC_D}{2c}\|f\|_{L^2}.
\end{equation}
For a real periodic $C^1$ function $f$, the function $\mathcal K_Xf$ lies in $H^1$, and
\begin{equation}\label{eq:bem-dsKH1}
\|\partial_\xi\mathcal K_Xf\|_{L^2}\le\frac{MC_D}{2c}\|f'\|_{L^2}.
\end{equation}
\end{lemma}

\begin{proof}
We use the following criterion. Let $u\in L^2$, and let $E_n\in L^2$ satisfy $\|E_n\|_{L^2}\le A$ and $\int E_n\varphi\,\dd\xi\to-\int u\varphi'\,\dd\xi$ for every real periodic $C^1$ function $\varphi$. Then, $u\in H^1$ and $\|u'\|_{L^2}\le A$. Indeed, a subsequence of $(E_n)$ converges weakly in $L^2$ to some $E$ with $\|E\|_{L^2}\le A$, and $\int E\varphi\,\dd\xi=-\int u\varphi'\,\dd\xi$ shows that $E$ is the weak derivative of $u$.

For real $f\in L^2$, the function $\mathcal L_Xf$ lies in $L^2$ by Schur's test, since its kernel is symmetric and integrable uniformly in $\xi$. By \eqref{eq:bem-david} and $|X'|\ge c$, we have $\|\mathcal D^\delta_Xf\|_{L^2}\le\tfrac12MC_D\|f/X'\|_{L^2}\le(MC_D/2c)\|f\|_{L^2}$. By \eqref{eq:bem-dsL} and the criterion, applied to $E_n=\mathcal D^{\delta_n}_Xf$ for a sequence $\delta_n\to0$, we obtain $\mathcal L_Xf\in H^1$ and \eqref{eq:bem-dsLH1}.

For a real periodic $C^1$ function $f$, we multiply \eqref{eq:bem-dsC} by $-\tfrac12$ and take real parts. By \eqref{eq:bem-KC}, we obtain
\[
\int_\Tor\mathcal K^\delta_Xf\,\varphi'\,\dd\xi=\int_\Tor E^\delta\varphi\,\dd\xi-\tfrac12\operatorname{Re}\frac{R_\delta}{\pi\iu},\qquad E^\delta=\tfrac12\operatorname{Re}\big[X'\,\mathcal C^\delta_X(f'/X')\big],
\]
and \eqref{eq:bem-david} gives $\|E^\delta\|_{L^2}\le(MC_D/2c)\|f'\|_{L^2}$. The truncations $\mathcal K^\delta_Xf$ converge to $\mathcal K_Xf$ almost everywhere by \eqref{eq:bem-pv}, and they are bounded in $L^2$ uniformly in $\delta$ by \eqref{eq:bem-KC} and \eqref{eq:bem-david}, so they converge in $L^1$ by Vitali's theorem. Hence, the left side tends to $\int\mathcal K_Xf\,\varphi'\,\dd\xi$, and $R_\delta\to0$ by \cref{lem:bem-trunc}. The criterion, applied to $E_n=-E^{\delta_n}$ for a sequence $\delta_n\to0$, gives $\mathcal K_Xf\in H^1$ and \eqref{eq:bem-dsKH1}.
\end{proof}

The bounds on $H^{1/2}$ are obtained from those on $L^2$ and $H^1$ by interpolation. Let $\psi_0=1$, and for $n\ge1$ let $\psi_n(\xi)=\sqrt2\cos2\pi n\xi$ and $\psi_{-n}(\xi)=\sqrt2\sin2\pi n\xi$. These functions form an orthonormal basis of real $L^2$, and with $w_n=(1+n^2)^{1/2}$ and $u_n=\langle u,\psi_n\rangle$ we have $\|u\|_{H^s}^2=\sum_nw_n^{2s}u_n^2$ for real $u\in L^2$. A trigonometric polynomial is a finite linear combination of the $\psi_n$, and the Fourier partial sums of $u$ are $\sum_{|n|\le N}u_n\psi_n$.

We use the exact interpolation of the $K$-method \cite[Thm.~2.2(i)]{ChandlerWildeHewettMoiola2015} and the corrected form \cite[Cor.~3.2$'$]{ChandlerWildeHewettMoiola2022} of \cite[Cor.~3.2]{ChandlerWildeHewettMoiola2015}. By the definition of the norms in \cref{sec:setting}, the map $u\mapsto(\hat u(n))_{n\in\Z}$ is, for every $s\in\R$, a unitary isomorphism of the complex space $H^s$ onto the space of sequences that are square summable with the weight $(1+n^2)^s$. For $s_0<s_1$, the space $H^{s_1}$ is continuously embedded in $H^{s_0}$, so $(H^{s_0},H^{s_1})$ is a compatible pair in the sense of \cite[\S2]{ChandlerWildeHewettMoiola2022}, and the map is injective on the sum $H^{s_0}$ of the two spaces, as \cite[Cor.~3.2$'$]{ChandlerWildeHewettMoiola2022} requires. Hence, this corollary identifies the space $K_{1/2,2}(H^{s_0},H^{s_1})$ of the $K$-method with $H^{(s_0+s_1)/2}$, with equality of norms. Since the $K$-method is exact of exponent $\frac12$ \cite[Thm.~2.2(i)]{ChandlerWildeHewettMoiola2015}, a linear map from $H^{s_0}$ to $H^{t_0}$ with norm bounded by $B_0$ whose restriction to $H^{s_1}$ maps into $H^{t_1}$ with norm bounded by $B_1$, where $s_0<s_1$ and $t_0<t_1$, has norm bounded by $(B_0B_1)^{1/2}$ from $H^{(s_0+s_1)/2}$ to $H^{(t_0+t_1)/2}$. A real linear map $T$ extends to complex functions by $T(u+\iu v)=Tu+\iu Tv$, and the extension has the same norms as $T$ between the spaces $H^s$, since $\|u+\iu v\|_{H^s}^2=\|u\|_{H^s}^2+\|v\|_{H^s}^2$ for real $u$ and $v$.

\begin{proposition}[Interpolation]\label{lem:bem-fs}
Let $A,A_0,A_1\ge0$.
\begin{enumerate}[label=(\alph*)]
\item Let $T$ be a linear map from real $L^2$ to real $L^2$ such that $\langle Tu,v\rangle=\langle u,Tv\rangle$ and $\|Tu\|_{H^1}\le A\|u\|_{L^2}$ for all $u,v\in L^2$. Then, $\|Tu\|_{H^{1/2}}\le A\|u\|_{H^{-1/2}}$ for all $u\in L^2$.
\item Let $T$ be a linear map from the real trigonometric polynomials to real $L^2$ such that $\|Tu\|_{L^2}\le A_0\|u\|_{L^2}$ and $\|Tu\|_{H^1}\le A_1\|u\|_{H^1}$ for all trigonometric polynomials $u$. Then, $\|Tu\|_{H^{1/2}}\le(A_0A_1)^{1/2}\|u\|_{H^{1/2}}$ for all trigonometric polynomials $u$.
\end{enumerate}
\end{proposition}

Part (a) is the case $s_0=-1$, $s_1=0$, $t_0=0$ and $t_1=1$, because $T$ extends to a map from $H^{-1}$ to $L^2$ with norm bounded by $A$. Indeed, Parseval's identity gives $|\langle f,g\rangle|\le\|f\|_{H^{-1}}\|g\|_{H^1}$ for $f,g\in L^2$, so the symmetry of $T$ gives $|\langle Tu,v\rangle|=|\langle u,Tv\rangle|\le A\|u\|_{H^{-1}}\|v\|_{L^2}$, and $v=Tu$ gives $\|Tu\|_{L^2}\le A\|u\|_{H^{-1}}$ on the dense subspace $L^2$ of $H^{-1}$. Part (b) is the case $s_0=t_0=0$ and $s_1=t_1=1$. Indeed, $T$ extends by continuity to a map from $L^2$ to $L^2$ with norm bounded by $A_0$, and its restriction to $H^1$ has norm bounded by $A_1$ in $H^1$ by the following density argument with $V=L^2$ and $s=1$.

We use this density argument several times. Let $s\in\{\frac12,1\}$, let $V\subset L^2$ be a linear space that contains the trigonometric polynomials, and let $T:V\to L^2$ be linear with $\|Tu\|_{L^2}\le A_0\|u\|_{L^2}$ for $u\in V$ and $\|Tu\|_{H^s}\le A\|u\|_{H^s}$ for trigonometric polynomials $u$. Then, $Tu\in H^s$ and $\|Tu\|_{H^s}\le A\|u\|_{H^s}$ for every $u\in V\cap H^s$. Indeed, the Fourier partial sums $u_N$ of $u$ satisfy $\|u_N\|_{H^s}\le\|u\|_{H^s}$ and $u_N\to u$ in $L^2$, so $Tu_N\to Tu$ in $L^2$, and the Fourier coefficients of $Tu_N$ converge to those of $Tu$. Hence, Fatou's lemma for the series $\sum_nw_n^{2s}\langle Tu,\psi_n\rangle^2$ gives the claim.

The following lemma states classical results \cite{CoifmanMeyer1978,David1984} uniformly on the class $\mathfrak F(M,c)$, and its proof rests on \cref{lem:bem-holo}.

\begin{lemma}[Lipschitz dependence and first derivatives on chord-arc curves]\label{lem:bem-cont}
Let $M\ge c>0$. There is a constant $C$, depending only on $M$ and $c$, such that the following hold for all $X,Y\in\mathfrak F(M,c)$.
\begin{enumerate}[label=(\roman*)]
\item For $\varsigma\in L^2$, we have
\[
\|(\mathcal L_X-\mathcal L_Y)\varsigma\|_{H^1}\le C\|X'-Y'\|_{L^\infty}\|\varsigma\|_{L^2},\qquad
\|(\mathcal L_X-\mathcal L_Y)\varsigma\|_{H^{1/2}}\le C\|X'-Y'\|_{L^\infty}\|\varsigma\|_{H^{-1/2}}.
\]
\item For real Lipschitz functions $a$, we have
\[
\|(\mathcal K_X-\mathcal K_Y)a\|_{L^2}\le C\|X'-Y'\|_{L^\infty}\|a\|_{L^2},\qquad
\|(\mathcal K_X-\mathcal K_Y)a\|_{H^1}\le C\|X'-Y'\|_{L^\infty}\|a\|_{H^1}.
\]
\item Let $e$ and $z$ be Lipschitz directions, let $X_\varepsilon=X+\varepsilon_1e+\varepsilon_2z$ for $\varepsilon=(\varepsilon_1,\varepsilon_2)\in\R^2$, let $\varsigma,\psi\in L^2$, and let $a$ be a real Lipschitz function. Then, the functions
\[
F(\varepsilon)=\langle\mathcal L_{X_\varepsilon}\varsigma,\psi\rangle,\qquad G(\varepsilon)=\langle\mathcal K_{X_\varepsilon}a,\psi\rangle
\]
are twice continuously differentiable near $\varepsilon=0$. Their first derivatives at $\varepsilon=0$ are
\[
\begin{aligned}
\partial_{\varepsilon_1}F(0)&=-\frac1{2\pi}\int_\Tor\!\int_\Tor\operatorname{Re}q_e(\xi,\eta)\,\varsigma(\eta)\,\psi(\xi)\,\dd\eta\,\dd\xi,\\
\partial_{\varepsilon_1}G(0)&=\frac1{2\pi}\int_\Tor\!\int_\Tor\operatorname{Im}q_e(\xi,\eta)\,a'(\eta)\,\psi(\xi)\,\dd\eta\,\dd\xi,
\end{aligned}
\]
and they satisfy
\[
|\partial_{\varepsilon_1}F(0)|\le C\|e'\|_{L^\infty}\|\varsigma\|_{H^{-1/2}}\|\psi\|_{H^{-1/2}},\qquad
|\partial_{\varepsilon_1}G(0)|\le C\|e'\|_{L^\infty}\|a\|_{H^{1/2}}\|\psi\|_{H^{-1/2}}.
\]
The same holds for $\partial_{\varepsilon_2}$, with $z$ in place of $e$. The mixed second derivatives at $\varepsilon=0$ are
\[
\begin{aligned}
\partial_{\varepsilon_1}\partial_{\varepsilon_2}F(0)&=\frac1{2\pi}\int_\Tor\!\int_\Tor\operatorname{Re}\big(q_eq_z\big)(\xi,\eta)\,\varsigma(\eta)\,\psi(\xi)\,\dd\eta\,\dd\xi,\\
\partial_{\varepsilon_1}\partial_{\varepsilon_2}G(0)&=-\frac1{2\pi}\int_\Tor\!\int_\Tor\operatorname{Im}\big(q_eq_z\big)(\xi,\eta)\,a'(\eta)\,\psi(\xi)\,\dd\eta\,\dd\xi .
\end{aligned}
\]
\end{enumerate}
\end{lemma}

\begin{proof}
Let $s=\|X'-Y'\|_{L^\infty}$, let $C_D$ and $C_\Delta$ be the constants of \eqref{eq:bem-david} and \cref{lem:bem-holo} for $\mathfrak F(M,c)$, and let $T=\mathcal L_X-\mathcal L_Y$. All constants below depend only on $M$ and $c$, and $|\langle u,\psi\rangle|\le\|u\|_{H^{1/2}}\|\psi\|_{H^{-1/2}}$ by Parseval's identity.

For (i), the chords satisfy $|\Delta X|,|\Delta Y|\ge c\,d(\xi,\eta)$ and $\big||\Delta X|-|\Delta Y|\big|\le|\Delta(X-Y)|\le s\,d(\xi,\eta)$, so $\big|\log|\Delta X|-\log|\Delta Y|\big|\le s/c$. Hence, the kernel of $T$ is bounded by $s/(2\pi c)$, and $\|T\varsigma\|_{L^2}\le Cs\|\varsigma\|_{L^2}$. Let $\mathcal D^\delta_X$ be as in \eqref{eq:bem-dsL}. For $0<\delta<\frac12$, we have
\[
X'\mathcal C^\delta_X(\varsigma/X')-Y'\mathcal C^\delta_Y(\varsigma/Y')=(X'-Y')\,\mathcal C^\delta_X(\varsigma/X')+Y'(\mathcal C^\delta_X-\mathcal C^\delta_Y)(\varsigma/X')+Y'\,\mathcal C^\delta_Y\big(\varsigma(1/X'-1/Y')\big),
\]
and $|1/X'-1/Y'|\le s/c^2$ almost everywhere. By \eqref{eq:bem-david} and \cref{lem:bem-holo}, there is a constant $K_0$, depending only on $M$ and $c$, such that, uniformly in $\delta$,
\begin{equation}\label{eq:bem-Dlip}
\|\mathcal D^\delta_X\varsigma-\mathcal D^\delta_Y\varsigma\|_{L^2}\le K_0\,s\,\|\varsigma\|_{L^2}.
\end{equation}
By \eqref{eq:bem-dsL} for $X$ and for $Y$, the integrals $\int(\mathcal D^\delta_X\varsigma-\mathcal D^\delta_Y\varsigma)\varphi\,\dd\xi$ tend to $-\int T\varsigma\,\varphi'\,\dd\xi$ for every real periodic $C^1$ function $\varphi$. The criterion at the beginning of the proof of \cref{lem:bem-H1} gives $T\varsigma\in H^1$ and $\|\partial_\xi T\varsigma\|_{L^2}\le K_0s\|\varsigma\|_{L^2}$, which, together with the bound on $L^2$, proves the first bound of (i). The kernel of $T$ is real, symmetric and bounded, so $T$ is symmetric on $L^2$, and \cref{lem:bem-fs}(a) gives the second bound of (i).

For (ii), let $a$ be a real Lipschitz function. By \eqref{eq:bem-pv}, the truncations $\mathcal K^\delta_Xa$ and $\mathcal K^\delta_Ya$ converge almost everywhere to $\mathcal K_Xa$ and $\mathcal K_Ya$. By \eqref{eq:bem-KC} and \cref{lem:bem-holo}, we have $\|(\mathcal K^\delta_X-\mathcal K^\delta_Y)a\|_{L^2}\le\frac12C_\Delta s\|a\|_{L^2}$ for every $\delta$, and Fatou's lemma gives the first bound of (ii). Let first $f$ be a real periodic $C^1$ function. We subtract the identity for $\mathcal K^\delta_Yf$ in the proof of \cref{lem:bem-H1} from that for $\mathcal K^\delta_Xf$. The difference $E^\delta_X-E^\delta_Y$ of the functions $E^\delta$ there is half the real part of the decomposition above with $f'$ in place of $\varsigma$, so $\|E^\delta_X-E^\delta_Y\|_{L^2}\le K_0s\|f'\|_{L^2}$ uniformly in $\delta$. As in that proof, the remainders tend to $0$ and the left sides converge, so the criterion gives $\|\partial_\xi(\mathcal K_X-\mathcal K_Y)f\|_{L^2}\le K_0s\|f'\|_{L^2}$. With the first bound, we obtain
\[
\|(\mathcal K_X-\mathcal K_Y)f\|_{H^1}\le C\,s\,\|f\|_{H^1}.
\]
The density argument after \cref{lem:bem-fs}, applied in $H^1$ to $\mathcal K_X-\mathcal K_Y$ on the space of real Lipschitz functions, on which the first bound of (ii) holds, gives $(\mathcal K_X-\mathcal K_Y)a\in H^1$ and the second bound of (ii).

For (iii), we have $|\Delta(\varepsilon_1e+\varepsilon_2z)|\le\rho\,d(\xi,\eta)$ with $\rho=|\varepsilon_1|\|e'\|_{L^\infty}+|\varepsilon_2|\|z'\|_{L^\infty}$, so $q_\varepsilon=\varepsilon_1q_e+\varepsilon_2q_z$ satisfies $|q_\varepsilon|\le\frac12$ for $\rho\le c/2$. Let $\operatorname{Log}$ be the principal branch of the logarithm, which is holomorphic on $\{|w-1|<1\}$. For such $\varepsilon$, the curve $X_\varepsilon$ lies in $\mathfrak F(M+c/2,c/2)$, and $\log|\Delta X_\varepsilon|=\log|\Delta X|+\operatorname{Re}\operatorname{Log}(1+q_\varepsilon)$. Hence, we have
\[
F(\varepsilon)=F(0)-\frac1{2\pi}\int_\Tor\!\int_\Tor\operatorname{Re}\operatorname{Log}\big(1+q_\varepsilon(\xi,\eta)\big)\,\varsigma(\eta)\,\psi(\xi)\,\dd\eta\,\dd\xi .
\]
The function $\varepsilon\mapsto\operatorname{Re}\operatorname{Log}(1+q_\varepsilon(\xi,\eta))$ is smooth for $\rho<c/2$, with derivatives bounded uniformly in $(\xi,\eta)$, and $\varsigma(\eta)\psi(\xi)$ is integrable on $\Tor\times\Tor$. Hence, differentiation under the integral sign shows that $F$ is twice continuously differentiable near $0$, and the formulas for $\partial_{\varepsilon_1}F(0)$ and $\partial_{\varepsilon_1}\partial_{\varepsilon_2}F(0)$ follow from $\partial_{\varepsilon_1}\operatorname{Re}\operatorname{Log}(1+q_\varepsilon)=\operatorname{Re}q_e$ and $\partial_{\varepsilon_1}\partial_{\varepsilon_2}\operatorname{Re}\operatorname{Log}(1+q_\varepsilon)=-\operatorname{Re}(q_eq_z)$ at $\varepsilon=0$. For the bound, let $C'$ be the constant of (i) for $\mathfrak F(M+c/2,c/2)$. For small $t\ne0$, we obtain
\[
|F(t,0)-F(0)|=|\langle(\mathcal L_{X+te}-\mathcal L_X)\varsigma,\psi\rangle|\le C'|t|\,\|e'\|_{L^\infty}\|\varsigma\|_{H^{-1/2}}\|\psi\|_{H^{-1/2}},
\]
and dividing by $|t|$ and letting $t\to0$ gives the bound for $\partial_{\varepsilon_1}F(0)$.

For $G$, let $\rho<c/2$, and let $\xi$ be a point at which $X$, $e$ and $z$ are differentiable, so that $|X'(\xi)|\ge c$. Let $r_\varepsilon=\Delta X_\varepsilon/\Delta X=1+q_\varepsilon$. For $0<\delta<\frac12$, the function $\eta\mapsto\operatorname{Im}\operatorname{Log}r_\varepsilon(\xi,\eta)$ is Lipschitz on the arc $\{\eta\in\Tor\mid d(\xi,\eta)>\delta\}$, and its derivative is
\[
\operatorname{Im}\frac{X_\varepsilon'(\eta)}{X_\varepsilon(\eta)-X_\varepsilon(\xi)}-\operatorname{Im}\frac{X'(\eta)}{X(\eta)-X(\xi)}
\]
almost everywhere. Integrating by parts against $a$ on this arc, we obtain
\[
\mathcal K^\delta_{X_\varepsilon}a(\xi)-\mathcal K^\delta_Xa(\xi)=b_\delta(\xi)+\frac1{2\pi}\int_{d(\xi,\eta)>\delta}\operatorname{Im}\operatorname{Log}r_\varepsilon(\xi,\eta)\,a'(\eta)\,\dd\eta,
\]
where $b_\delta(\xi)=-\frac1{2\pi}\big[a(\xi-\delta)\operatorname{Im}\operatorname{Log}r_\varepsilon(\xi,\xi-\delta)-a(\xi+\delta)\operatorname{Im}\operatorname{Log}r_\varepsilon(\xi,\xi+\delta)\big]$. Both values $r_\varepsilon(\xi,\xi\pm\delta)$ tend to $X_\varepsilon'(\xi)/X'(\xi)$, so $b_\delta(\xi)\to0$, and the integrand is bounded by $\pi\|a'\|_{L^\infty}$. Since both truncations converge almost everywhere by \eqref{eq:bem-pv}, we obtain, for almost every $\xi$,
\[
\mathcal K_{X_\varepsilon}a(\xi)-\mathcal K_Xa(\xi)=\frac1{2\pi}\int_\Tor\operatorname{Im}\operatorname{Log}\big(1+q_\varepsilon(\xi,\eta)\big)\,a'(\eta)\,\dd\eta .
\]
Since $\mathcal K_Xa\in L^2$ by the bound after \eqref{eq:bem-pv}, multiplying by $\psi$ and integrating gives
\[
G(\varepsilon)=G(0)+\frac1{2\pi}\int_\Tor\!\int_\Tor\operatorname{Im}\operatorname{Log}\big(1+q_\varepsilon(\xi,\eta)\big)\,a'(\eta)\,\psi(\xi)\,\dd\eta\,\dd\xi .
\]
As for $F$, the function $G$ is twice continuously differentiable near $0$, and its first and mixed second derivatives are the stated ones, since $\partial_{\varepsilon_1}\operatorname{Im}\operatorname{Log}(1+q_\varepsilon)=\operatorname{Im}q_e$ and $\partial_{\varepsilon_1}\partial_{\varepsilon_2}\operatorname{Im}\operatorname{Log}(1+q_\varepsilon)=-\operatorname{Im}(q_eq_z)$ at $\varepsilon=0$. For the bound, we use $\operatorname{Im}q_e=\operatorname{Re}q_{-\iu e}$. Hence, $\partial_{\varepsilon_1}G(0)$ is $-\partial_{\varepsilon_1}F(0)$ for the direction $-\iu e$ in place of $e$ and the density $a'$ in place of $\varsigma$. Since $\|(-\iu e)'\|_{L^\infty}=\|e'\|_{L^\infty}$ and $\|a'\|_{H^{-1/2}}\le C\|a\|_{H^{1/2}}$, the bound for $F$ gives the bound for $G$. The derivatives in $\varepsilon_2$ are treated in the same way.
\end{proof}

On $\mathfrak C$ we apply \cref{lem:bem-holo,lem:bem-cont} with $M=M_1$ and $c=c_{\rm ca}$. The class $\mathfrak F(M_*+c_*/2,c_*/2)$, on which the constant $C_L$ fixed before \cref{def:nbhd} is taken, contains $\mathfrak C$ because $r_1\le c_*/4$.

\subsection{Uniform invertibility and the Schur complements}\label{sec:bem-unif}

\begin{lemma}[Uniform invertibility]\label{lem:bem-Vreg}
Assume \ref{K1} and \ref{K3}. There are $c,C>0$, depending only on $M_*$, $c_*$, $\mu_1$ and $B_1$, such that for all $X\in\mathfrak C$ the following hold.
\begin{enumerate}[label=(\alph*)]
\item $\langle\mathcal L_X\varsigma,\varsigma\rangle\ge c\|\varsigma\|_{H^{-1/2}}^2$ for all $\varsigma\in L^2$.
\item The operator $\mathcal L_X$ maps $L^2$ onto $H^1$, and for all $\varsigma\in L^2$ we have
\[
\begin{aligned}
\|\mathcal L_X\varsigma\|_{H^1}&\le C\|\varsigma\|_{L^2},&\qquad
\|\varsigma\|_{L^2}&\le C\|\mathcal L_X\varsigma\|_{H^1},\\
\|\mathcal L_X\varsigma\|_{H^{1/2}}&\le C\|\varsigma\|_{H^{-1/2}},&\qquad
\|\varsigma\|_{H^{-1/2}}&\le C\|\mathcal L_X\varsigma\|_{H^{1/2}}.
\end{aligned}
\]
\item For real Lipschitz functions $a$, we have
\[
\|\mathcal K_Xa\|_{L^2}\le C\|a\|_{L^2},\qquad
\|\mathcal K_Xa\|_{H^{1/2}}\le C\|a\|_{H^{1/2}},\qquad
\|\mathcal K_Xa\|_{H^1}\le C\|a\|_{H^1}.
\]
If $X$ is a polygon, these bounds hold for all $a$ in $L^2$, in $H^{1/2}$ and in $H^1$, respectively, and for $a\in H^1$ the densities $\varsigma_\pm\in L^2$ with $\mathcal L_X\varsigma_\pm=\tfrac12a\pm\mathcal K_Xa$ satisfy $\|\varsigma_\pm\|_{L^2}\le C\|a\|_{H^1}$.
\end{enumerate}
\end{lemma}

\begin{proof}
Let $X\in\mathfrak C$, and let $\kappa=2C_Lr_1$, so that $\kappa\le\min(\frac12\mu_1,(2B_1)^{-1})$ by \eqref{eq:nbhd}. The curves $X$ and $x$ lie in $\mathfrak F(M_*+c_*/2,c_*/2)$, and $\|X'-x'\|_{L^\infty}\le2r_1$. Hence, \cref{lem:bem-cont}(i) gives, for $\varsigma\in L^2$ and $D\varsigma=(\mathcal L_X-\mathcal L_x)\varsigma$,
\[
\|D\varsigma\|_{H^1}\le\kappa\|\varsigma\|_{L^2},\qquad\|D\varsigma\|_{H^{1/2}}\le\kappa\|\varsigma\|_{H^{-1/2}} .
\]
By \ref{K3} and $|\langle D\varsigma,\varsigma\rangle|\le\|D\varsigma\|_{H^{1/2}}\|\varsigma\|_{H^{-1/2}}$, we obtain $\langle\mathcal L_X\varsigma,\varsigma\rangle\ge(\mu_1-\kappa)\|\varsigma\|_{H^{-1/2}}^2\ge\frac12\mu_1\|\varsigma\|_{H^{-1/2}}^2$, which is (a).

For (b), we first prove the forward bounds. Let $M_2=M_*+c_*/2$ and $c_2=c_*/2$. Since $c_2\,d(\xi,\eta)\le|\Delta X|\le M_2\,d(\xi,\eta)$, the kernel of $\mathcal L_X$ is bounded by $C(1+|\log d(\xi,\eta)|)$, so the Cauchy--Schwarz inequality gives $|\mathcal L_X\varsigma(\xi)|\le C\|\varsigma\|_{L^2}$ for every $\xi$, and \eqref{eq:bem-dsLH1} gives $\|\mathcal L_X\varsigma\|_{H^1}\le C\|\varsigma\|_{L^2}$. The kernel is symmetric, so \cref{lem:bem-fs}(a) gives $\|\mathcal L_X\varsigma\|_{H^{1/2}}\le C\|\varsigma\|_{H^{-1/2}}$. Next, \ref{K3} gives $\|\varsigma\|_{L^2}\le B_1\|\mathcal L_x\varsigma\|_{H^1}\le B_1\|\mathcal L_X\varsigma\|_{H^1}+B_1\kappa\|\varsigma\|_{L^2}$, and $B_1\kappa\le\frac12$ gives $\|\varsigma\|_{L^2}\le2B_1\|\mathcal L_X\varsigma\|_{H^1}$. By (a), we have $\frac12\mu_1\|\varsigma\|_{H^{-1/2}}^2\le\|\mathcal L_X\varsigma\|_{H^{1/2}}\|\varsigma\|_{H^{-1/2}}$, so $\|\varsigma\|_{H^{-1/2}}\le(2/\mu_1)\|\mathcal L_X\varsigma\|_{H^{1/2}}$. To show that $\mathcal L_X$ maps $L^2$ onto $H^1$, let $g\in H^1$. Since $D\varsigma\in H^1$ by the bound of $D$, \ref{K3} defines a map $\Phi$ on $L^2$ by $\mathcal L_x\Phi\varsigma=g-D\varsigma$, and by the same bound it is a contraction with constant $B_1\kappa\le\frac12$. Its fixed point $\varsigma$ satisfies $\mathcal L_X\varsigma=\mathcal L_x\varsigma+D\varsigma=g$.

For (c), let $C_D$ be the constant of \eqref{eq:bem-david} for $\mathfrak F(M_2,c_2)$. The bound on $L^2$ is the bound after \eqref{eq:bem-pv}. By the definition of the norms, $\|u\|_{H^1}\le\|u\|_{L^2}+(2\pi)^{-1}\|u'\|_{L^2}$ and $\|u'\|_{L^2}\le2\pi\|u\|_{H^1}$ for $u\in H^1$. For a real periodic $C^1$ function $f$, these inequalities, the bound on $L^2$ and \cref{lem:bem-H1} give
\[
\|\mathcal K_Xf\|_{H^1}\le\tfrac12C_D\|f\|_{L^2}+\frac{M_2C_D}{4\pi c_2}\|f'\|_{L^2}\le\tfrac12C_D\big(1+M_2/c_2\big)\|f\|_{H^1}.
\]
Hence, \cref{lem:bem-fs}(b), applied to $\mathcal K_X$ on the trigonometric polynomials with $A_0=\frac12C_D$ and $A_1=\frac12C_D(1+M_2/c_2)$, gives $\|\mathcal K_Xf\|_{H^{1/2}}\le C\|f\|_{H^{1/2}}$ for trigonometric polynomials $f$. For a real Lipschitz function $a$, the bounds on $H^1$ and on $H^{1/2}$ follow from the density argument after \cref{lem:bem-fs}, applied with $s=1$ and $s=\frac12$ to $\mathcal K_X$ on the space of real Lipschitz functions, on which the bound after \eqref{eq:bem-pv} holds.

Let now $X$ be a polygon, let $\xi$ be a parameter that is not a vertex, and let $\ell_p$ be the line through the edge $\Gamma_p$. For $\eta$ on $p$, we have $|\operatorname{Im}(X'(\eta)/(X(\eta)-X(\xi)))|\le|X'(\eta)|/\operatorname{dist}(X(\xi),\ell_p)$ if $X(\xi)\notin\ell_p$, and the kernel vanishes if $X(\xi)\in\ell_p$. So the kernel of $\mathcal K_X$ is bounded in $\eta$, and for every integrable $a$ the truncations converge at $\xi$ by dominated convergence. Hence, \eqref{eq:bem-KC}, \eqref{eq:bem-david} and Fatou's lemma give $\|\mathcal K_Xa\|_{L^2}\le\frac12C_D\|a\|_{L^2}$ for all $a\in L^2$, and the density argument after \cref{lem:bem-fs}, now with $V=L^2$, gives the bounds on $H^1$ and on $H^{1/2}$. Finally, for $a\in H^1$, the function $\frac12a\pm\mathcal K_Xa$ lies in $H^1$ with norm bounded by $(\frac12+C)\|a\|_{H^1}$, and (b) gives $\varsigma_\pm$ with $\|\varsigma_\pm\|_{L^2}\le C\|a\|_{H^1}$.
\end{proof}

The Dirichlet-to-Neumann map of \cref{sec:setting} is related to the layer operators by the jump relations of the layer potentials on Lipschitz domains \cite[Lemma~4.1]{Costabel1988}, the representations of the Steklov--Poincar\'e operators \cite[(6.40)--(6.43), \S7.5]{Steinbach2008} and Maue's formula \cite[Thm.~6.15]{Steinbach2008}. We state the resulting representation in the parameter form as part (ii) of the following proposition. Part (i) is elementary, and we show it after the proposition, together with the reasons why the cited results apply.

\begin{proposition}[Layer representation of the Dirichlet-to-Neumann map]\label{prop:bem-dtn}
Let $Y$ be a polygon of $\mathfrak C$, the curve $x$, or a curve $x+\varepsilon e$ with a P1 direction $e$ and $|\varepsilon|\,\|e'\|_{L^\infty}\le2r_1$.
\begin{enumerate}[label=(\roman*)]
\item The image $\Gamma_Y=Y(\Tor)$ is a Jordan curve, oriented counter-clockwise, that bounds a Lipschitz domain, and $\bn=\Jm Y'/|Y'|$ is its outward normal.
\item If $Y$ is a polygon, let $g\in H^1$ and let $w$ be a P1 function. Otherwise, let $g\in H^2$ and let $w$ be Lipschitz. Let $\varsigma,\varsigma_{\rm int},\varsigma_{\rm ext}\in L^2$ and $c,c'\in\R$ satisfy
\[
\mathcal L_Y\varsigma=g+c,\quad\int_\Tor\varsigma\,\dd\xi=0,\qquad\mathcal L_Y\varsigma_{\rm int}=(\tfrac12+\mathcal K_Y)g,\qquad\mathcal L_Y\varsigma_{\rm ext}=(\tfrac12-\mathcal K_Y)g+c',\quad\int_\Tor\varsigma_{\rm ext}\,\dd\xi=0 .
\]
Then, we have
\begin{equation}\label{eq:bem-sp}
2\langle\mathcal L_Y\partial_\xi g,\partial_\xi w\rangle+\langle\varsigma_{\rm int},(\tfrac12+\mathcal K_Y)w\rangle+\langle\varsigma_{\rm ext},(\tfrac12-\mathcal K_Y)w\rangle=\int_\Tor\varsigma\,w\,\dd\xi .
\end{equation}
If $Y$ is not a polygon, both sides of \eqref{eq:bem-sp} are also equal to $\int_\Tor(\Lam_Yg)\,w\,|Y'|\,\dd\xi$.
\end{enumerate}
\end{proposition}

If $Y$ is $x$ or $x+\varepsilon e$ and $g\in H^2$, then $\int_\Tor(\Lam_Yg-\varsigma/|Y'|)\,w\,|Y'|\,\dd\xi=0$ for all Lipschitz $w$, which are dense in $L^2$, so $\Lam_Yg=\varsigma/|Y'|$ almost everywhere.

For (i), the curves on the segment from $x$ to $Y$ lie in the convex set $\mathfrak C$, so they are injective by their chord constant (\cref{sec:setting}). Since $\|Y'-x'\|_{L^\infty}\le\frac12c_*$ by \eqref{eq:nbhd}, every short arc of $\Gamma_Y$ is a Lipschitz graph over a tangent line of $x$, and the signed area, which depends continuously on the curve in $W^{1,\infty}$, stays positive along the segment. Part (ii) is the cited representation written in the parameter, where $\mathcal L_Y\varsigma$ is the single layer operator applied to the arclength density $\varsigma/|Y'|$. The curves of the proposition are smooth off finitely many corners, so $\sigma=\frac12$ off the corners by \cite[(6.11) and p.~123]{Steinbach2008}. Maue's formula is stated for piecewise smooth curves and for continuous functions that are $C^1$ on the smooth pieces, and for this reason (ii) is stated only for such curves. It extends to the functions $g$ and $w$ of (ii) by density in $H^{1/2}$, because both of its sides are continuous in the norm of $H^{1/2}$.

\begin{lemma}[Gauge solutions]\label{lem:bem-gauge}
Assume \ref{K1} and \ref{K3}, and let $Y\in\mathfrak C$. For every $g\in H^1$ there is a pair $(\varsigma,c)\in L^2\times\R$ with
\[
\mathcal L_Y\varsigma=g+c,\qquad\int_\Tor\varsigma\,\dd\xi=0 .
\]
The pair is unique, and $\|\varsigma\|_{L^2}+|c|\le C\|g\|_{H^1}$ with $C$ independent of $Y$.
\end{lemma}

\begin{proof}
Let $\varsigma_1=\mathcal L_Y^{-1}\one\in L^2$, which exists by \cref{lem:bem-Vreg}(b). Since $\|\one\|_{H^{1/2}}=\|\one\|_{H^1}=1$, \cref{lem:bem-Vreg}(a) and (b) give
\[
I_1=\int_\Tor\varsigma_1\,\dd\xi=\langle\mathcal L_Y\varsigma_1,\varsigma_1\rangle\ge c\|\varsigma_1\|_{H^{-1/2}}^2\ge c\,C^{-2},\qquad\|\varsigma_1\|_{L^2}\le C .
\]
Let $g\in H^1$ and $\rho=\mathcal L_Y^{-1}g$, so that $\|\rho\|_{L^2}\le C\|g\|_{H^1}$. By \cref{lem:bem-Vreg}(b), a pair $(\varsigma,c)$ is a solution if and only if $\varsigma=\rho+c\varsigma_1$ and $c=-I_1^{-1}\int\rho\,\dd\xi$. Hence, the solution exists and is unique, and we obtain
\[
|c|\le I_1^{-1}\|\rho\|_{L^2}\le C\|g\|_{H^1},\qquad\|\varsigma\|_{L^2}\le\|\rho\|_{L^2}+|c|\,\|\varsigma_1\|_{L^2}\le C\|g\|_{H^1} .
\]
\end{proof}

We write $\mathcal L_{Y,0}^{-1}g=\varsigma$ for the bordered inverse of \cref{lem:bem-gauge}, $g\in H^1$.

The jump constant of the double layer is given by Gauss's lemma, which we state in the parameter form.

\begin{proposition}[Gauss's lemma]\label{prop:bem-gauss}
Assume \ref{K1} and \ref{K3}, and let $X$ be a polygon of $\mathfrak C$. Then, $\mathcal K_X\one=-\tfrac12$ at every parameter that is not a vertex.
\end{proposition}

\begin{proof}
The polygon $X$ bounds a Lipschitz domain by \cref{prop:bem-dtn}(i), and the kernel of $\mathcal K_X$ is that of the double layer operator $K$ of \cite[(6.14)]{Steinbach2008}, written in the parameter. Let $\xi$ be a parameter that is not a vertex. The boundary is straight near $X(\xi)$, so $\sigma=\frac12$ at $X(\xi)$ by \cite[(6.11) and p.~123]{Steinbach2008}. The kernel of $\mathcal K_X$ is bounded at $\xi$, as in the proof of \cref{lem:bem-Vreg}(c), so the principal value of $K$, which is cut off in the image, and that of $\mathcal K_X$, which is cut off in the parameter, both equal the absolutely convergent integral.

The representation formula \cite[(6.1)]{Steinbach2008} with $u=\one$ shows that the double layer potential $W\one$ equals $-1$ in the interior. Hence, \cite[Lemma~6.11]{Steinbach2008} with $v=\one$ gives $-1=-1+\sigma+K\one$ at $X(\xi)$, that is, $K\one=-\sigma=-\frac12$. The lemma applies on the polygon, because for $v=\one$ the term of its proof that contains the differences $v(y)-v(X(\xi))$ vanishes, and the remaining terms use only the representation formula \cite[(6.1)]{Steinbach2008} with $u=\one$ on the Lipschitz domain and the definition \cite[(6.11)]{Steinbach2008} of $\sigma$. Since the two principal values coincide, we obtain $\mathcal K_X\one(\xi)=K\one(X(\xi))=-\frac12$.
\end{proof}

\begin{definition}[Galerkin inverses]\label{def:bem-galinv}
Let $X\in\mathfrak C$ and $g\in H^{1/2}$. The Galerkin inverse $\mathcal L_{X,h}^{-1}g\in\mathcal P_0$ is the density with $\langle\mathcal L_X\mathcal L_{X,h}^{-1}g,\chi\rangle=\langle g,\chi\rangle$ for all $\chi\in\mathcal P_0$. The bordered Galerkin inverse $\mathcal L_{X,h,0}^{-1}g=\varsigma$, which is the operator form of $\tilde C$, is the density $\varsigma\in\mathcal P_0$ with $\int\varsigma\,\dd\xi=0$ for which there is $\nu\in\R$ with $\langle\mathcal L_X\varsigma,\chi\rangle=\langle g,\chi\rangle+\nu\int\chi\,\dd\xi$ for all $\chi\in\mathcal P_0$. For $g\in H^1$, the exact inverses are $\mathcal L_X^{-1}g$ and the bordered inverse $\mathcal L_{X,0}^{-1}g$ of \cref{lem:bem-gauge}, and $\mathcal L_{X,h}^{-1}g$ and $\mathcal L_{X,h,0}^{-1}g$ are the Galerkin solutions corresponding to them.
\end{definition}

\begin{lemma}[Galerkin inverses]\label{lem:bem-galinv}
Assume \ref{K1} and \ref{K3}, and let $X\in\mathfrak C$. Then, the following hold with constants that do not depend on $X$ and $h$.
\begin{enumerate}[label=(\alph*)]
\item The Galerkin inverses of $g\in H^{1/2}$ exist and are unique, that is, the Galerkin system and the bordered system with its multiplier $\nu$ are uniquely solvable, and the exact inverses of $g\in H^1$ exist. The Galerkin inverses are symmetric, $\langle\mathcal L_{X,h}^{-1}a,b\rangle=\langle a,\mathcal L_{X,h}^{-1}b\rangle$, and likewise for the bordered one. Both pairings equal $\tilde b^\T\tilde V^{-1}\tilde a$, respectively $\tilde b^\T\tilde C\tilde a$, with $\tilde a_p=\langle a,\iota_p\rangle$ and $\tilde b_p=\langle b,\iota_p\rangle$.
\item Let $\varsigma$ be $\mathcal L_{X,h}^{-1}g$ or $\mathcal L_{X,h,0}^{-1}g$ with $g\in H^{1/2}$, or $\mathcal L_X^{-1}g$ or $\mathcal L_{X,0}^{-1}g$ with $g\in H^1$. Then, we have
\begin{equation}\label{eq:bem-inv}
\|\varsigma\|_{H^{-1/2}}\le C\|g\|_{H^{1/2}} .
\end{equation}
For the Galerkin inverses, we have moreover $\langle g,\varsigma\rangle=\langle\mathcal L_X\varsigma,\varsigma\rangle\ge0$.
\item Let $g\in H^1$, let $\varsigma$ be $\mathcal L_X^{-1}g$ or $\mathcal L_{X,0}^{-1}g$, and let $\varsigma_h$ be the corresponding Galerkin solution. Then, $\|\varsigma\|_{L^2}\le C\|g\|_{H^1}$, and we have
\begin{equation}\label{eq:bem-gal}
\|\varsigma-\varsigma_h\|_{H^{-1/2}}\le C\inf_{\chi}\|\varsigma-\chi\|_{H^{-1/2}}\le C\|\varsigma-\Pi_0\varsigma\|_{H^{-1/2}}\le Ch^{1/2}\|\varsigma\|_{L^2},
\end{equation}
where $\chi$ ranges over $\mathcal P_0$, or over its subspace of mean zero in the bordered case, and $\Pi_0$ is the $L^2$ projection onto $\mathcal P_0$. If moreover $\varsigma\in H^1$, then $\|\varsigma-\Pi_0\varsigma\|_{H^{-1/2}}\le Ch^{3/2}\|\varsigma\|_{H^1}$ and $\|\varsigma-\varsigma_h\|_{H^{-1/2}}\le Ch^{3/2}\|\varsigma\|_{H^1}$.
\item Let $X\in\Nd$, which uses \ref{K2} with $r=1$ through the constant $C_I$. In the setting of (c), we have $\|\varsigma_h\|_{L^2}\le C\|\varsigma\|_{L^2}\le C\|g\|_{H^1}$. In particular, the Galerkin densities $\varsigma_u=\tilde V^{-1}\tilde Q_{\rm int}u$ of a P1 trace $u$, and likewise the bordered exterior densities, satisfy
\begin{equation}\label{eq:bem-dens}
\|\varsigma_u\|_{L^2}\le C\|u\|_{H^1}.
\end{equation}
\end{enumerate}
\end{lemma}

\begin{proof}
For (a), the form $\langle\mathcal L_X\cdot,\cdot\rangle$ is coercive on $\mathcal P_0$ by \cref{lem:bem-Vreg}(a), so $\tilde V$ is positive definite, and its restriction to the densities of mean zero is positive definite as well, which makes the bordered system uniquely solvable. The exact inverses exist by \cref{lem:bem-Vreg}(b) and \cref{lem:bem-gauge}. By the definition of $\tilde V$ and $\tilde C$, the pairings $\langle\mathcal L_{X,h}^{-1}a,b\rangle$ and $\langle\mathcal L_{X,h,0}^{-1}a,b\rangle$ equal $\tilde b^\T\tilde V^{-1}\tilde a$ and $\tilde b^\T\tilde C\tilde a$, and since both matrices are symmetric, the Galerkin inverses are symmetric.

For (b), testing with $\varsigma$ gives $\langle\mathcal L_X\varsigma,\varsigma\rangle=\langle g,\varsigma\rangle$ in each case, because the constants and the multipliers drop against densities of mean zero. Hence, \cref{lem:bem-Vreg}(a) gives \eqref{eq:bem-inv}, and for the Galerkin inverses the same identity gives $\langle g,\varsigma\rangle\ge0$.

For (c), we have $\|\varsigma\|_{L^2}\le C\|g\|_{H^1}$ by \cref{lem:bem-Vreg,lem:bem-gauge}. We apply C\'ea's lemma \cite[Thm.~8.1]{Steinbach2008} on the space of $L^2$ densities with the norm of $H^{-1/2}$, with the Galerkin space $\mathcal P_0$, or its subspace of mean zero in the bordered case. On this space, the form $\langle\mathcal L_X\cdot,\cdot\rangle$ is elliptic by \cref{lem:bem-Vreg}(a) and bounded by the bound of $\mathcal L_X$ from $H^{-1/2}$ to $H^{1/2}$ on $L^2$ densities in \cref{lem:bem-Vreg}(b), and these two constants determine the constant of the lemma. The projection $\Pi_0$ preserves the mean, so it maps into the subspace of mean zero in the bordered case. With the approximation property of $\mathcal P_0$ \cite[Thm.~10.4]{Steinbach2008}, we obtain \eqref{eq:bem-gal}. If moreover $\varsigma\in H^1$, then the same approximation property gives $\|\varsigma-\Pi_0\varsigma\|_{H^{-1/2}}\le Ch^{3/2}\|\varsigma\|_{H^1}$, and the second inequality of \eqref{eq:bem-gal} gives the bound for $\varsigma-\varsigma_h$.

For (d), the inverse estimate $\|\chi\|_{L^2}\le Ch^{-1/2}\|\chi\|_{H^{-1/2}}$ on $\mathcal P_0$, applied to $\varsigma_h-\Pi_0\varsigma$, and \eqref{eq:bem-gal} give $\|\varsigma_h\|_{L^2}\le C\|\varsigma\|_{L^2}\le C\|g\|_{H^1}$. For a P1 trace $u$, the data $Q_{\rm int}u$ and $Q_{\rm ext}u$ lie in $H^1$ with norms bounded by $C\|u\|_{H^1}$ by \cref{lem:bem-Vreg}(c), and their Galerkin solutions are the Galerkin densities of $u$. This gives \eqref{eq:bem-dens}.
\end{proof}

\begin{definition}[Schur forms]\label{def:bem-schur}
Let $X\in\mathfrak C$, and let $Q_{\rm int}=\tfrac12+\mathcal K_X$ and $Q_{\rm ext}=\tfrac12-\mathcal K_X$ in the parameter form. For Lipschitz functions $v$ and $w$, put
\begin{equation}\label{eq:bem-ext}
\mathfrak s_{X,h}(v,w)=2\langle\mathcal L_X\partial_\xi v,\partial_\xi w\rangle+\langle\mathcal L_{X,h}^{-1}Q_{\rm int}v,Q_{\rm int}w\rangle+\langle\mathcal L_{X,h,0}^{-1}Q_{\rm ext}v,Q_{\rm ext}w\rangle ,
\end{equation}
and let $\langle\mathcal S_Xv,w\rangle$ be the same expression with $\mathcal L_X^{-1}$ and $\mathcal L_{X,0}^{-1}$ in place of the Galerkin inverses. The exact inverses apply because $Q_{\rm int}v$ and $Q_{\rm ext}v$ lie in $H^1$ by \cref{lem:bem-Vreg}(c). If $X$ is a polygon, the same holds for $v\in H^1$ by \cref{lem:bem-Vreg}(c), and we use both forms also for $v\in H^1$ and Lipschitz $w$.
\end{definition}

\begin{lemma}[Schur forms]\label{lem:bem-schur}
Assume \ref{K1} and \ref{K3}, and let $X\in\mathfrak C$. Then, the following hold with constants that do not depend on $X$ and $h$.
\begin{enumerate}[label=(\alph*)]
\item If $X$ is one of the curves of \cref{prop:bem-dtn}, and $g$ and $w$ are as in \cref{prop:bem-dtn}(ii), then $\langle\mathcal S_Xg,w\rangle=\int_\Tor\varsigma\,w\,\dd\xi$ with the density $\varsigma$ of \cref{prop:bem-dtn}(ii). If $X$ is moreover not a polygon, then $\langle\mathcal S_Xg,w\rangle=\int_\Tor(\Lam_Xg)\,w\,|X'|\,\dd\xi$.
\item Let $X$ be a polygon. For P1 functions $v$ and $w$, we have $\mathfrak s_{X,h}(\one,w)=\mathfrak s_{X,h}(w,\one)=0$ and $\mathfrak s_{X,h}(v,w)=w^\T\Lh(X)v$.
\item For Lipschitz functions $v$ and $w$, we have $|\mathfrak s_{X,h}(v,w)|\le C\|v\|_{H^{1/2}}\|w\|_{H^{1/2}}$.
\item Let $X$ be a polygon and let $w$ be a P1 function. Then, we have
\[
\mathfrak s_{X,h}(w,w)\ge2\langle\mathcal L_X\partial_\xi w,\partial_\xi w\rangle\ge c\|\partial_\xi w\|_{H^{-1/2}}^2 .
\]
If the vertex values of $w$ sum to zero, then $\mathfrak s_{X,h}(w,w)\ge c\|w\|_{H^{1/2}}^2$.
\item Let $u$ and $w$ be Lipschitz, or let $X$ be a polygon, $u\in H^1$ and $w$ Lipschitz. Let $\varsigma_{\rm int}=\mathcal L_X^{-1}Q_{\rm int}u$ and $\varsigma_{\rm ext}=\mathcal L_{X,0}^{-1}Q_{\rm ext}u$, and let $\varsigma_{{\rm int},h}$ and $\varsigma_{{\rm ext},h}$ be the corresponding Galerkin solutions. Then, $\|\varsigma_{\rm int}\|_{L^2}+\|\varsigma_{\rm ext}\|_{L^2}\le C\|u\|_{H^1}$, and
\[
\mathfrak s_{X,h}(u,w)-\langle\mathcal S_Xu,w\rangle=\langle\varsigma_{{\rm int},h}-\varsigma_{\rm int},Q_{\rm int}w\rangle+\langle\varsigma_{{\rm ext},h}-\varsigma_{\rm ext},Q_{\rm ext}w\rangle .
\]
Moreover, we have
\begin{equation}\label{eq:bem-sym}
|\mathfrak s_{X,h}(u,w)-\langle\mathcal S_Xu,w\rangle|\le C\big(\|\varsigma_{\rm int}-\Pi_0\varsigma_{\rm int}\|_{H^{-1/2}}+\|\varsigma_{\rm ext}-\Pi_0\varsigma_{\rm ext}\|_{H^{-1/2}}\big)\|w\|_{H^{1/2}}
\end{equation}
and
\begin{equation}\label{eq:bem-symest}
|\mathfrak s_{X,h}(u,w)-\langle\mathcal S_Xu,w\rangle|\le Ch^{1/2}\|u\|_{H^1}\|w\|_{H^{1/2}}.
\end{equation}
\end{enumerate}
\end{lemma}

\begin{proof}
Part (a) is \cref{prop:bem-dtn}(ii) with $Y=X$.

For (b), $\mathcal K_X\one=-\tfrac12$ off the vertices by \cref{prop:bem-gauss}, which gives $Q_{\rm int}\one=0$, and the bordered system absorbs the constant datum $Q_{\rm ext}\one=\one$, so $\mathfrak s_{X,h}(\one,w)=\mathfrak s_{X,h}(w,\one)=0$ for P1 functions $w$. The Galerkin solves have the consistent loads $\tilde Q_{\rm int}v$ and $\tilde Q_{\rm ext}v$ for P1 functions $v$, so \cref{lem:bem-reduction} gives $\mathfrak s_{X,h}(v,w)=w^\T\Lh(X)v$, because the deflation $P$ changes $v$ and $w$ only by multiples of $\one$.

Part (c) follows from \cref{lem:bem-Vreg}(b), (c) and \eqref{eq:bem-inv}.

For (d), the two Schur terms of $\mathfrak s_{X,h}(w,w)$ are nonnegative by \cref{lem:bem-galinv}(b), and $\partial_\xi w$ is a density of $\mathcal P_0$, so \cref{lem:bem-Vreg}(a) gives the first claim (cf.\ \cite[Lemma~12.11]{Steinbach2008}). If the vertex values of $w$ sum to zero, then $\int_\Tor w\,\dd\xi=0$, and comparing Fourier coefficients gives $\|w\|_{H^{1/2}}\le\|\partial_\xi w\|_{H^{-1/2}}$, which proves the second claim.

For (e), the two forms differ only in their Schur terms, which gives the identity. By \cref{lem:bem-Vreg}(c), we have $\|Q_{\rm int}w\|_{H^{1/2}}+\|Q_{\rm ext}w\|_{H^{1/2}}\le C\|w\|_{H^{1/2}}$, so the identity and \eqref{eq:bem-gal} give \eqref{eq:bem-sym}. By \cref{lem:bem-Vreg}(b), (c) and \cref{lem:bem-gauge}, we have $\|\varsigma_{\rm int}\|_{L^2}+\|\varsigma_{\rm ext}\|_{L^2}\le C\|u\|_{H^1}$, and \eqref{eq:bem-gal} gives \eqref{eq:bem-symest}.
\end{proof}

\begin{lemma}[{The kernel of $\Lh$}]\label{lem:bem-a}
For every polygon $X\in\mathfrak C$, we have $\ker\Lh(X)=\operatorname{span}\{\one\}$.
\end{lemma}

\begin{proof}
Let $X\in\mathfrak C$ be a polygon and let $\Lh(X)v=0$. By \cref{lem:bem-reduction}, $0=v^\T\Lh(X)v$ is the sum of $2(\mathsf D_\xi Pv)^\T\tilde V(\mathsf D_\xi Pv)$, $(\tilde Q_{\rm int}Pv)^\T\tilde V^{-1}(\tilde Q_{\rm int}Pv)$ and $(\tilde Q_{\rm ext}Pv)^\T\tilde C(\tilde Q_{\rm ext}Pv)$. The matrix $\tilde V$ is positive definite by \cref{lem:bem-Vreg}(a), since $a^\T\tilde Va=\langle\mathcal L_X\varsigma,\varsigma\rangle$ for the density $\varsigma\in\mathcal P_0$ with coefficients $a$. Hence, $\tilde V^{-1}$ is positive definite as well, and $\tilde C$ is positive semidefinite, since $a^\T\tilde Ca=\varsigma^\T\tilde V\varsigma$ for the solution $\varsigma$ of the bordered system with load $a$. Thus, all three terms vanish, and in particular $\mathsf D_\xi Pv=0$. Then, $Pv$ is constant, and since $\one^\T Pv=0$, we obtain $Pv=0$, that is, $v\in\operatorname{span}\{\one\}$. Conversely, $\Lh(X)\one=0$.
\end{proof}

\subsection{Higher regularity on the smooth curve}

The consistency estimates compare the polygon with the smooth curve $x$, and on $x$ they need the double layer operator and the inverses of the single layer on Sobolev spaces of higher order, with constants that are uniform on $\mathscr K$. A function $g\in H^2$ is identified with its continuous representative, which is Lipschitz, so that $\mathcal K_xg$ is defined by \eqref{eq:bem-pv}.

\begin{lemma}[Uniform higher mapping properties]\label{lem:bem-higher}
Assume \ref{K1}, \ref{K2} with $r=1$, and \ref{K3}.
\begin{enumerate}[label=(\alph*)]
\item For every $g\in H^2$, we have $\mathcal K_xg\in H^2$ and
\[
\|\mathcal K_xg\|_{H^2}\le C\|g\|_{H^2}.
\]
\item Every $\varsigma\in L^2$ with $\mathcal L_x\varsigma\in H^2$ lies in $H^1$. If \ref{K2} holds with $r=3$, then every $\varsigma\in L^2$ with $\mathcal L_x\varsigma\in H^3$ lies in $H^2$. In both cases, with $s=1$ and $s=2$, respectively, we have
\[
\|\varsigma\|_{H^s}\le C\|\mathcal L_x\varsigma\|_{H^{s+1}} .
\]
\item The bordered inverse satisfies $\|\mathcal L_{x,0}^{-1}g\|_{H^1}\le C\|g\|_{H^2}$ for $g\in H^2$, and, if \ref{K2} holds with $r=3$, $\|\mathcal L_{x,0}^{-1}g\|_{H^2}\le C\|g\|_{H^3}$ for $g\in H^3$. For $g\in H^2$, the densities $\mathcal L_x^{-1}(\tfrac12+\mathcal K_x)g$ and $\mathcal L_{x,0}^{-1}(\tfrac12-\mathcal K_x)g$ lie in $H^1$, and their norms in $H^1$ are bounded by $C\|g\|_{H^2}$.
\end{enumerate}
Here $C$ depends only on $M_*$, $c_*$, $\mu_1$, $B_1$, $C_1$ and the constant $C_r$ of \ref{K2} for the order named.
\end{lemma}

\begin{proof}
For $t\in\R$, let $\tau_tu=u(\cdot+t)$ and $x_t=x(\cdot+t)$. The substitution $\eta\mapsto\eta+t$ in the integrals that define $\mathcal L_x$ and $\mathcal K_x$ gives
\begin{equation}\label{eq:bem-tr}
\tau_t(\mathcal L_x\varsigma)=\mathcal L_{x_t}(\tau_t\varsigma),\qquad\tau_t(\mathcal K_xa)=\mathcal K_{x_t}(\tau_ta)
\end{equation}
for $\varsigma\in L^2$ and Lipschitz $a$, and the same identities hold with the translate $x_\lambda$, $\lambda\in\R$, in place of $x$. The curve $x_t$ satisfies \ref{K1} with the constants of $x$, and $\|x_t'-x'\|_{L^\infty}\le|t|\,\|x''\|_{L^\infty}$. We fix $\delta\in(0,1]$ with $C_1\delta\le2r_1$. For $|t|\le\delta$, the curve $x_t$ then lies in $\mathfrak C$, so the bounds of \cref{lem:bem-Vreg} hold for $x_t$, and \cref{lem:bem-cont} applies to $x$ and $x_t$ with the constant $C_L$ of the class $\mathfrak F(M_*+c_*/2,c_*/2)$.

By the definition of the norms in \cref{sec:setting}, we have $\|u\|_{H^s}^2=\sum_kw_k^{2s}|\hat u(k)|^2$ with $w_k=(1+k^2)^{1/2}$, and $\widehat{\tau_tu}(k)=\eu^{2\pi\iu kt}\hat u(k)$. Hence, $\tau_t$ preserves the norm of $H^s$, and $|\eu^{2\pi\iu kt}-1|\le2\pi w_k|t|$ gives $\|\tau_tu-u\|_{H^s}\le2\pi|t|\,\|u\|_{H^{s+1}}$ for $u\in H^{s+1}$. Conversely, let $s\ge0$ and $B\ge0$, and let $u\in H^s$ satisfy $\|\tau_tu-u\|_{H^s}\le Bt$ for $0<t\le\delta$. Letting $t\to0$ in each finite part of the series $t^{-2}\|\tau_tu-u\|_{H^s}^2=\sum_kw_k^{2s}|\eu^{2\pi\iu kt}-1|^2t^{-2}|\hat u(k)|^2\le B^2$, we obtain $\sum_{|k|\le K}w_k^{2s}(2\pi k)^2|\hat u(k)|^2\le B^2$, and since $w_k^{2s+2}=w_k^{2s}(1+k^2)$, we have
\begin{equation}\label{eq:bem-fq}
u\in H^{s+1},\qquad\|u\|_{H^{s+1}}\le\|u\|_{H^s}+\frac B{2\pi} .
\end{equation}

For (a), let $g\in H^2$ and $u=\mathcal K_xg$. By \cref{lem:bem-Vreg}(c), we have $u\in H^1$ and $\|u\|_{H^1}\le C\|g\|_{H^1}$. Let $0<t\le\delta$. By \eqref{eq:bem-tr}, we have
\[
\tau_tu-u=\underbrace{\mathcal K_{x_t}(\tau_tg-g)}_{=:\,R_g}+\underbrace{(\mathcal K_{x_t}-\mathcal K_x)g}_{=:\,R_x},
\]
where $g$ and $\tau_tg-g$ are Lipschitz. By \cref{lem:bem-Vreg}(c) for $x_t$, $R_g$ has norm in $H^1$ bounded by $C\|\tau_tg-g\|_{H^1}\le Ct\|g\|_{H^2}$. By \cref{lem:bem-cont}(ii) and $\|x_t'-x'\|_{L^\infty}\le C_1t$, $R_x$ has norm in $H^1$ bounded by $Ct\|g\|_{H^1}$. Hence, $\|\tau_tu-u\|_{H^1}\le Ct\|g\|_{H^2}$, and \eqref{eq:bem-fq} with $s=1$ gives $u\in H^2$ and $\|u\|_{H^2}\le C\|g\|_{H^2}$.

For (b) with $s=1$, let $\varsigma\in L^2$ with $G=\mathcal L_x\varsigma\in H^2$, and let $0<t\le\delta$. By \eqref{eq:bem-tr}, $\mathcal L_{x_t}(\tau_t\varsigma)=\tau_tG$, and therefore we have
\[
\mathcal L_x(\tau_t\varsigma-\varsigma)=\underbrace{(\tau_tG-G)}_{=:\,T_G}-\underbrace{(\mathcal L_{x_t}-\mathcal L_x)\tau_t\varsigma}_{=:\,T_x} .
\]
By the bound of translations above, $\|T_G\|_{H^1}\le2\pi t\|G\|_{H^2}$. By \cref{lem:bem-cont}(i), we have $\|T_x\|_{H^1}\le Ct\|\varsigma\|_{L^2}$, and $\|\varsigma\|_{L^2}\le B_1\|G\|_{H^1}$ by \ref{K3}. Hence, \ref{K3} applied to $\tau_t\varsigma-\varsigma$ gives
\[
\|\tau_t\varsigma-\varsigma\|_{L^2}\le B_1\|\mathcal L_x(\tau_t\varsigma-\varsigma)\|_{H^1}\le Ct\|G\|_{H^2}.
\]
Let $\varphi$ be a real periodic $C^1$ function, and let $E_t=(\tau_t\varsigma-\varsigma)/t$, so that $\|E_t\|_{L^2}\le C\|G\|_{H^2}$. The substitution $\xi\mapsto\xi+t$ gives $\int E_t\varphi\,\dd\xi=\int\varsigma(\xi)\,t^{-1}(\varphi(\xi-t)-\varphi(\xi))\,\dd\xi$, and as $t\to0$ the difference quotients of $\varphi$ tend to $-\varphi'$ and are bounded by $\sup|\varphi'|$. Hence, dominated convergence gives $\int E_t\varphi\,\dd\xi\to-\int\varsigma\varphi'\,\dd\xi$, and the criterion at the beginning of the proof of \cref{lem:bem-H1}, applied to $E_n=E_{t_n}$ for a sequence $t_n\to0$, gives $\varsigma\in H^1$ and $\|\varsigma'\|_{L^2}\le C\|G\|_{H^2}$. Together with $\|\varsigma\|_{L^2}\le B_1\|G\|_{H^1}$ and the inequality for $\|u\|_{H^1}$ in the proof of \cref{lem:bem-Vreg}(c), this gives $\|\varsigma\|_{H^1}\le C\|G\|_{H^2}$.

For (b) with $s=2$, let \ref{K2} hold with $r=3$, so that $|x^{(i)}|\le C_3$ for $1\le i\le4$. We first show that the second differences of the single layer in the translations of the curve satisfy
\begin{equation}\label{eq:bem-sd}
\|(\mathcal L_{x_{a+b}}-\mathcal L_{x_a}-\mathcal L_{x_b}+\mathcal L_x)\psi\|_{H^1}\le C|a|\,|b|\,\|\psi\|_{L^2}\qquad(a,b\in\R,\ \psi\in L^2).
\end{equation}
Shifting $\eta$ by an integer changes none of the kernels, so we may take $|\xi-\eta|\le\frac34$. For such $\xi,\eta$ and $c\in\R$, let
\[
P_c(\xi,\eta)=\int_0^1x'\big(\eta+\theta(\xi-\eta)+c\big)\,\dd\theta,
\]
so that $x_c(\xi)-x_c(\eta)=(\xi-\eta)P_c(\xi,\eta)$. Since $d(\xi,\eta)\ge|\xi-\eta|/3$, \ref{K1} gives $|P_c|\ge m$ with $m=c_*/3$. The terms $\log|\xi-\eta|$ cancel, and the kernel of the operator in \eqref{eq:bem-sd} is $-\frac1{2\pi}k$ with
\[
k=\log|P_{a+b}|-\log|P_a|-\log|P_b|+\log|P_0|=\log|p|-\log|q|,\qquad p=P_{a+b}P_0,\quad q=P_aP_b .
\]
The bounds of $x''$ and $x'''$ under the integral give $|P_a-P_0|\le C_3|a|$ and $|P_{a+b}-P_a-P_b+P_0|\le C_3|a|\,|b|$, and $|P_c|\le C_3$ for every $c$. Hence, the identity
\[
p-q=(P_{a+b}-P_a-P_b+P_0)P_0-(P_a-P_0)(P_b-P_0)
\]
gives $|p-q|\le C|a|\,|b|$. Since $|p|,|q|\ge m^2$ and $|\log y_1-\log y_2|\le|y_1-y_2|/m^2$ for $y_1,y_2\ge m^2$, we obtain $|k|\le C|a|\,|b|$. For the Lipschitz bound in $\xi$, fix $\eta$, and let $p_i$ and $q_i$ be the values of $p$ and $q$ at $(\xi_i,\eta)$, where $|\xi_i-\eta|\le\frac34$. The bounds of $x''$, $x'''$ and $x''''$ under the integrals show that $P_c$, $P_a-P_0$ and $P_{a+b}-P_a-P_b+P_0$ are Lipschitz in $\xi$ with the constants $C_3$, $C_3|a|$ and $C_3|a|\,|b|$. Hence, the identity above gives $|(p_1-q_1)-(p_2-q_2)|\le C|a|\,|b|\,|\xi_1-\xi_2|$, and we also have $|q_1-q_2|\le C|\xi_1-\xi_2|$. Since $|q_i|\ge m^2$ and $|p_2-q_2|\le C|a|\,|b|$, the identity
\[
\frac{p_1}{q_1}-\frac{p_2}{q_2}=\frac{(p_1-q_1)-(p_2-q_2)}{q_1}+\frac{(p_2-q_2)(q_2-q_1)}{q_1q_2}
\]
bounds the left side by $C|a|\,|b|\,|\xi_1-\xi_2|$. Since $k=\log|p/q|$ and $|p/q|\ge m^2/C_3^2$, the same bound holds for $|k(\xi_1,\eta)-k(\xi_2,\eta)|$. The ranges $|\xi-\eta|\le\frac34$ for neighboring integer shifts overlap, so $k(\cdot,\eta)$ is Lipschitz with constant $C|a|\,|b|$ for every $\eta$. Hence, the image of $\psi$ in \eqref{eq:bem-sd} is bounded by $C|a|\,|b|\,\|\psi\|_{L^1}$ and Lipschitz with the same constant, so its norm in $H^1$ is bounded by $C|a|\,|b|\,\|\psi\|_{L^1}\le C|a|\,|b|\,\|\psi\|_{L^2}$. This proves \eqref{eq:bem-sd}.

Next, for $\lambda\in\R$ and $u\in H^1$, let $D_\lambda u=(\mathcal L_{x_\lambda}-\mathcal L_x)u$. We show that
\begin{equation}\label{eq:bem-Dk}
\|D_\lambda u\|_{H^2}\le C|\lambda|\,\|u\|_{H^1}.
\end{equation}
The curve $x_\lambda$ satisfies \ref{K1} with the constants of $x$, so \cref{lem:bem-cont} applies to $x$ and $x_\lambda$ for every $\lambda\in\R$. By \cref{lem:bem-cont}(i) and $\|x_\lambda'-x'\|_{L^\infty}\le C_3|\lambda|$, we have $\|D_\lambda u\|_{H^1}\le C|\lambda|\,\|u\|_{L^2}$. For $0<t\le\delta$, \eqref{eq:bem-tr} and $(x_\lambda)_t=x_{\lambda+t}$ give
\[
\tau_t(D_\lambda u)-D_\lambda u=(\mathcal L_{x_{\lambda+t}}-\mathcal L_{x_\lambda}-\mathcal L_{x_t}+\mathcal L_x)\tau_tu+D_\lambda(\tau_tu-u).
\]
By \eqref{eq:bem-sd} and the bound of $D_\lambda$ in $H^1$, the right side has norm in $H^1$ bounded by $C|\lambda|t\|u\|_{L^2}+C|\lambda|\,\|\tau_tu-u\|_{L^2}\le C|\lambda|t\|u\|_{H^1}$. Thus, \eqref{eq:bem-fq} with $s=1$ gives \eqref{eq:bem-Dk}.

Now let $\varsigma\in L^2$ with $G=\mathcal L_x\varsigma\in H^3$. By the case $s=1$, we have $\varsigma\in H^1$ and $\|\varsigma\|_{H^1}\le C\|G\|_{H^2}$. For $0<t\le\delta$, as in the case $s=1$, we have
\[
\mathcal L_x(\tau_t\varsigma-\varsigma)=(\tau_tG-G)-D_t(\tau_t\varsigma),
\]
and by \eqref{eq:bem-Dk}, the right side lies in $H^2$ with norm bounded by $Ct(\|G\|_{H^3}+\|\varsigma\|_{H^1})$. The case $s=1$, applied to $\tau_t\varsigma-\varsigma$, gives $\|\tau_t\varsigma-\varsigma\|_{H^1}\le Ct\|G\|_{H^3}$. Thus, \eqref{eq:bem-fq} with $s=1$ gives $\varsigma\in H^2$ and $\|\varsigma\|_{H^2}\le C\|G\|_{H^3}$.

For (c), let $g\in H^{s+1}$, where $s=1$, or $s=2$ if \ref{K2} holds with $r=3$, and let $(\varsigma,c)$ be the pair of \cref{lem:bem-gauge} for $Y=x$, so that $\mathcal L_x\varsigma=g+c$, $\int\varsigma\,\dd\xi=0$ and $|c|\le C\|g\|_{H^1}$. Since $\|g+c\|_{H^{s+1}}\le\|g\|_{H^{s+1}}+|c|$, (b) gives $\|\mathcal L_{x,0}^{-1}g\|_{H^s}=\|\varsigma\|_{H^s}\le C\|g\|_{H^{s+1}}$. For $g\in H^2$, (a) gives $\|(\tfrac12\pm\mathcal K_x)g\|_{H^2}\le C\|g\|_{H^2}$. Hence, (b) with $s=1$ gives the bound for $\mathcal L_x^{-1}(\tfrac12+\mathcal K_x)g$, and the bound for the bordered inverse with $s=1$ gives that for $\mathcal L_{x,0}^{-1}(\tfrac12-\mathcal K_x)g$.
\end{proof}

\subsection{Discrete regularity}\label{sec:bem-a}

What is specific to the matrix in this subsection is the treatment of the constants, the upgrade from $H^{1/2}$ to $H^1$ by an inverse estimate, and the lumped loads.

Let $X\in\Nd$. The lumped pairing is bounded by the $L^2$ norms. Indeed, $\mu_j\le M_1h$ because $|\Gamma_p|\le M_1h$, and every P1 function $v$ satisfies $\|v\|_{L^2}^2\ge\tfrac h3\sum_jv_j^2$, because on a panel with end values $a$ and $b$ the integral of $v^2$ is $\tfrac h3(a^2+ab+b^2)\ge\tfrac h6(a^2+b^2)$. Hence, the Cauchy--Schwarz inequality gives, for P1 functions $v$ and $w$,
\begin{equation}\label{eq:bem-lump}
|w^\T M_vv|\le M_1h\Big(\sum_jv_j^2\Big)^{1/2}\Big(\sum_jw_j^2\Big)^{1/2}\le3M_1\|v\|_{L^2}\|w\|_{L^2}.
\end{equation}
The loads of the shape derivative of $\Om$ are lumped loads of this kind. Since $\Om(X)$ is linear in $X$, we have $D\Om[e]f=\Om(e)f$ for a piecewise linear direction $e$ and a P1 field $f$, and this is the lumped load $M_vg$ of the P1 function $g$ with vertex values $\omega_j(e)\cdot f_j/\mu_j$, where $\omega_j(e)=\tfrac12\Jm(e_{j+1}-e_{j-1})$. Let $e'_p$ be the value of $e'$ on panel $p$. Then, $|\omega_j(e)|\le\tfrac h2(|e'_{j-1}|+|e'_j|)$ and $\mu_j\ge c_{\rm ca}h$, and since $\|g\|_{L^2}^2\le h\sum_jg_j^2$ and $h\sum_j|f_j|^2\le3\|f\|_{L^2}^2$, we obtain
\begin{equation}\label{eq:bem-dOm}
\|g\|_{L^2}\le C\|e'\|_{L^2}\|f\|_{L^\infty},\qquad\|g\|_{L^2}\le C\|e'\|_{L^\infty}\|f\|_{L^2}.
\end{equation}
For $e=X$, we have $|\omega_j(X)|\le\tfrac12(|\Gamma_{j-1}|+|\Gamma_j|)=\mu_j$, so $\Om(X)f$ is the lumped load of a P1 function with $L^2$ norm bounded by $C\|f\|_{L^2}$. Hence, \cref{lem:bem-areg} below and \eqref{eq:bem-dens} show that $z_f=\Lhp(X)\Om(X)f$ and its Galerkin density $\varsigma_f=\tilde V^{-1}\tilde Q_{\rm int}z_f$ satisfy
\begin{equation}\label{eq:bem-zf}
\|z_f\|_{H^1}+\|\varsigma_f\|_{L^2}\le C\|f\|_{L^2},
\end{equation}
and so does the bordered exterior density of $z_f$.

A load that is bounded by the $L^2$ norm is a consistent load. Let $M$ be the parameter mass matrix, $M_{ij}=\int_\Tor\varphi_i\varphi_j\,\dd\xi$, so that $v^\T Mv=\|v\|_{L^2}^2$. Let $y\in\R^N$ with $\one^\T y=0$ and $|w^\T y|\le B\|w\|_{L^2}$ for all P1 functions $w$ whose vertex values sum to zero. The vertex values of $v=M^{-1}y$ sum to zero, because $\sum_i\varphi_i=1$ and $\int\varphi_j\,\dd\xi=h$ give $h\sum_jv_j=\one^\T Mv=\one^\T y=0$, so testing with $w=v$ gives $\|v\|_{L^2}^2=v^\T y\le B\|v\|_{L^2}$. Since $g=v/|X'|$ satisfies $\int g\varphi_j\,|X'|\,\dd\xi=(Mv)_j=y_j$ and $|X'|\ge c_{\rm ca}$, the vector $y$ is the consistent load of $g$, and we have
\begin{equation}\label{eq:bem-riesz}
\|g\|_{L^2}\le c_{\rm ca}^{-1}B .
\end{equation}

The discrete regularity rests on the following perturbed Galerkin estimate.

\begin{lemma}[Perturbed Galerkin estimate]\label{lem:bem-strang}
Let $E$ and $V$ be real vector spaces, let $\iota:V\to E$ be linear and let $W\subset V$ be a subspace. Let $p:E\to[0,\infty)$ and $q:V\to[0,\infty)$, let $a:E\times V\to\R$ be linear in its first argument, let $\ell:V\to\R$, and let $\alpha>0$, $\beta\ge0$ and $\eta\ge0$. Assume that, for all $x\in E$ and $w\in W$,
\[
\alpha\,q(w)^2\le a(\iota w,w),\qquad|a(x,w)|\le\beta\,p(x)\,q(w),
\]
and let $u\in E$ and $u_h\in W$ satisfy, for all $w\in W$,
\[
a(\iota u_h,w)=\ell(w),\qquad|\ell(w)-a(u,w)|\le\eta\,q(w).
\]
Then, for every $v\in W$, we have
\[
q(u_h-v)\le\frac{\eta+\beta\,p(u-\iota v)}{\alpha}.
\]
\end{lemma}

\begin{proof}
Let $v\in W$ and $w=u_h-v\in W$. Since $a$ is linear in its first argument and $a(\iota u_h,w)=\ell(w)$, we have
\[
a(\iota w,w)=\big(\ell(w)-a(u,w)\big)+a(u-\iota v,w).
\]
The ellipticity, the consistency bound and the bound of $a$ give $\alpha\,q(w)^2\le\eta\,q(w)+\beta\,p(u-\iota v)\,q(w)$. If $q(w)>0$, we divide by $\alpha\,q(w)$, and if $q(w)=0$, the claim holds because its right side is nonnegative.
\end{proof}

Unlike Strang's lemma \cite[Thm.~8.3]{Steinbach2008}, the estimate allows different sizes $p$ and $q$ on the trial and test spaces and assumes nothing about the unperturbed problem.

The upgrade from $H^{1/2}$ to $H^1$ uses two properties of $\mathcal P_1$ on the uniform grid. By the inequality for $\|u\|_{H^1}$ in the proof of \cref{lem:bem-Vreg}(c), the following lemma also bounds the norms of $H^1$.

\begin{lemma}[P1 functions on the uniform grid]\label{lem:bem-p1}
There is a constant $C$, independent of $N$, with the following properties.
\begin{enumerate}[label=(\alph*)]
\item Every P1 function $w$ satisfies $\|w\|_{L^2}+\|w'\|_{L^2}\le Ch^{-1/2}\|w\|_{H^{1/2}}$.
\item Let $u\in H^1$, let $m_p$ be the mean of $u$ on the panel $[\xi_p,\xi_{p+1}]$, and let $v$ be the P1 function with $v(\xi_p)=m_p$. Then, we have
\[
\|u-v\|_{H^{1/2}}\le Ch^{1/2}\|u\|_{H^1},\qquad\|v\|_{L^2}+\|v'\|_{L^2}\le C\|u\|_{H^1}.
\]
\end{enumerate}
\end{lemma}

The first bound in (b) is the approximation property of $\mathcal P_1$ in $H^{1/2}$ \cite[Thm.~10.9]{Steinbach2008}, realized by an approximant that is also bounded in $H^1$.

\begin{proof}
For (a), we have $\|w\|_{L^2}\le\|w\|_{H^{1/2}}$, and the identity $\widehat{w'}(n)=2\pi\iu n\,\hat w(n)$ gives $\|w'\|_{H^{-1/2}}\le2\pi\|w\|_{H^{1/2}}$. Since $w'$ lies in $\mathcal P_0$, the inverse estimate on $\mathcal P_0$ \cite[Lemma~10.10]{Steinbach2008} gives
\[
\|w'\|_{L^2}\le Ch^{-1/2}\|w'\|_{H^{-1/2}}\le Ch^{-1/2}\|w\|_{H^{1/2}} ,
\]
and $h\le1$ proves (a).

For (b), let $D=\sum_p(m_{p+1}-m_p)^2$, with $m_N=m_0$. The Cauchy--Schwarz inequality on each panel gives $hm_p^2\le\int_{\xi_p}^{\xi_{p+1}}u^2\,\dd\xi$. On a panel with end values $a$ and $b$, the integral of $v^2$ is bounded by $\tfrac h2(a^2+b^2)$, so $\|v\|_{L^2}^2\le h\sum_pm_p^2\le\|u\|_{L^2}^2$. To bound $D$, let $z$ be the P1 function with $z(\xi_p)=m_p-m_{p-1}$. On the panel $[\xi_p,\xi_{p+1}]$, we have $z'=(z_{p+1}-z_p)/h$, and the integral of $u$ is $hm_p$. Summation by parts therefore gives
\[
\int_\Tor u\,z'\,\dd\xi=\sum_p(z_{p+1}-z_p)\,m_p=-\sum_pz_p\,(m_p-m_{p-1})=-D .
\]
We integrate by parts and use $\|z\|_{L^2}^2\le h\sum_pz_p^2=hD$ and $\|u'\|_{L^2}\le2\pi\|u\|_{H^1}$. Then, we obtain $D\le\|u'\|_{L^2}\|z\|_{L^2}\le2\pi\|u\|_{H^1}(hD)^{1/2}$, so that $D\le Ch\|u\|_{H^1}^2$. Since $v'=(m_{p+1}-m_p)/h$ on the panel $[\xi_p,\xi_{p+1}]$, we have $\|v'\|_{L^2}^2=D/h\le C\|u\|_{H^1}^2$, which proves the second bound of (b).

For the first bound, let $\Pi_0u\in\mathcal P_0$ be the function with the value $m_p$ on the panel $[\xi_p,\xi_{p+1}]$. The approximation property of $\mathcal P_0$ \cite[Thm.~10.4]{Steinbach2008} gives $\|u-\Pi_0u\|_{L^2}\le Ch\|u\|_{H^1}$. On the panel $[\xi_p,\xi_{p+1}]$, the function $v-\Pi_0u$ is linear, vanishes at $\xi_p$ and equals $m_{p+1}-m_p$ at $\xi_{p+1}$, so that $\|v-\Pi_0u\|_{L^2}^2=\tfrac h3D\le Ch^2\|u\|_{H^1}^2$. Hence, $\varepsilon_u=u-v$ satisfies $\|\varepsilon_u\|_{L^2}\le Ch\|u\|_{H^1}$, and $\|\varepsilon_u\|_{H^1}\le\|u\|_{H^1}+\|v\|_{L^2}+\|v'\|_{L^2}\le C\|u\|_{H^1}$. Finally, $(1+n^2)^{1/2}\le\tfrac12\big(h^{-1}+h(1+n^2)\big)$, so every $g\in H^1$ satisfies
\begin{equation}\label{eq:bem-interp}
\|g\|_{H^{1/2}}^2\le\tfrac12\big(h^{-1}\|g\|_{L^2}^2+h\|g\|_{H^1}^2\big).
\end{equation}
With $g=\varepsilon_u$, we obtain $\|\varepsilon_u\|_{H^{1/2}}^2\le Ch\|u\|_{H^1}^2$.
\end{proof}

In the Strang argument, the first argument of the discrete Schur form is a function of $H^{1/2}$. Let $X\in\Nd$. For $v\in H^{1/2}$ and a P1 function $w$, we define $\mathfrak s_{X,h}(v,w)$ by \eqref{eq:bem-ext}, with the hypersingular term read as $2\langle\partial_\xi v,\mathcal L_X\partial_\xi w\rangle$, where $\langle\partial_\xi v,\psi\rangle=\sum_n2\pi\iu n\,\hat v(n)\overline{\hat\psi(n)}$. Since the kernel of $\mathcal L_X$ is symmetric, this agrees with \eqref{eq:bem-ext} for $v\in H^1$. The Galerkin inverses are applied to $Q_{\rm int}v$ and $Q_{\rm ext}v$, which lie in $H^{1/2}$ by \cref{lem:bem-Vreg}(c), and the form is linear in $v$. Let the vertex values of $w$ sum to zero, so that $\|w\|_{H^{1/2}}\le\|\partial_\xi w\|_{H^{-1/2}}$ by \cref{sec:bem-unif}. Since $|n|\le(1+n^2)^{1/2}$, the hypersingular term is bounded by $4\pi\|v\|_{H^{1/2}}\|\mathcal L_X\partial_\xi w\|_{H^{1/2}}$, and \cref{lem:bem-Vreg}(b) bounds this by $C\|v\|_{H^{1/2}}\|\partial_\xi w\|_{H^{-1/2}}$. By \eqref{eq:bem-inv} and \cref{lem:bem-Vreg}(c), each Schur term is bounded by $C\|v\|_{H^{1/2}}\|w\|_{H^{1/2}}$. Hence, we have
\begin{equation}\label{eq:bem-sbd}
|\mathfrak s_{X,h}(v,w)|\le C\|v\|_{H^{1/2}}\|\partial_\xi w\|_{H^{-1/2}}\qquad\Big(v\in H^{1/2},\ w\in\mathcal P_1,\ \sum_jw_j=0\Big).
\end{equation}

\begin{lemma}[Discrete regularity]\label{lem:bem-areg}
Let $X\in\Nd$ and let $y$ be the load of $f$, either consistent, $y_j=\int f\varphi_j\,|X'|\,\dd\xi$ with $f\in L^2$, or lumped, $y_j=\mu_jf(\xi_j)$ with $f$ a P1 function. Then, $\|\Lhp y\|_{H^1}\le C\|f\|_{L^2}$.
\end{lemma}

\begin{proof}
We first treat consistent loads. Let $m=\one^\T y=\int f\,|X'|\,\dd\xi$ and $\tilde f=f-m/|X'|$. Since $\int\varphi_j\,\dd\xi=h$, the consistent load of $\tilde f$ is $y-mh\one=Py$, and $\Lhp y=\Lhp Py$ by \eqref{eq:mp-formula}. The parameter density $\varsigma=\tilde f|X'|=f|X'|-m$ has mean zero, and $\|\varsigma\|_{L^2}\le C\|f\|_{L^2}$. Let $c_0$ be the constant for which $u=\mathcal L_X\varsigma+c_0$ has mean zero. Since $|c_0|\le\|\mathcal L_X\varsigma\|_{L^2}$, \cref{lem:bem-Vreg}(b) gives $\|u\|_{H^1}\le2\|\mathcal L_X\varsigma\|_{H^1}\le C\|f\|_{L^2}$. The pair $(\varsigma,-c_0)$ solves $\mathcal L_X\varsigma=u-c_0$ and $\int\varsigma\,\dd\xi=0$, so \cref{prop:bem-dtn}(ii), applied to the polygon $X$ and to $g=u\in H^1$, gives
\[
\langle\mathcal S_Xu,w\rangle=\int_\Tor\varsigma\,w\,\dd\xi=w^\T Py
\]
for every P1 function $w$. The discrete solution $u_h=\Lhp y$ satisfies $\one^\T u_h=0$, and $\Lh u_h=Py$ by \eqref{eq:mp-formula}, which applies because $\ker\Lh=\operatorname{span}\{\one\}$ by \cref{lem:bem-a}. Hence, for every P1 function $w$ whose vertex values sum to zero, \cref{sec:bem-unif} gives $\mathfrak s_{X,h}(u_h,w)=w^\T\Lh u_h=w^\T Py$.

We apply \cref{lem:bem-strang} with $E=H^{1/2}$, $V=\R^N$, the map $\iota$ that sends vertex values to the P1 function, the subspace $W$ of vectors whose entries sum to zero, $a=\mathfrak s_{X,h}$, $\ell(w)=w^\T Py$, $p=\|\cdot\|_{H^{1/2}}$ and $q(w)=\|\partial_\xi w\|_{H^{-1/2}}$. The ellipticity holds with $\alpha=c$ by \cref{sec:bem-unif}, and the bound of $a$ holds with $\beta=C$ by \eqref{eq:bem-sbd}. For the consistency error, \eqref{eq:bem-symest} and $\|w\|_{H^{1/2}}\le\|\partial_\xi w\|_{H^{-1/2}}$ give, for $w\in W$,
\[
|\ell(w)-\mathfrak s_{X,h}(u,w)|=\big|\langle\mathcal S_Xu,w\rangle-\mathfrak s_{X,h}(u,w)\big|\le Ch^{1/2}\|u\|_{H^1}\|\partial_\xi w\|_{H^{-1/2}} ,
\]
so that $\eta=Ch^{1/2}\|u\|_{H^1}$. Let $v$ be the approximant of $u$ of \cref{lem:bem-p1}(b), let $\bar v=h\sum_jv_j$ be its mean, and let $Pv=v-\bar v$. Since $u$ has mean zero, we have $|\bar v|=|\int(u-v)\,\dd\xi|\le\|u-v\|_{H^{1/2}}$, and hence $\|u-Pv\|_{H^{1/2}}\le2\|u-v\|_{H^{1/2}}\le Ch^{1/2}\|u\|_{H^1}$. Applied with $Pv\in W$ in place of $v$, \cref{lem:bem-strang} gives
\[
\|\partial_\xi(u_h-Pv)\|_{H^{-1/2}}\le\frac{\eta+C\|u-Pv\|_{H^{1/2}}}{c}\le Ch^{1/2}\|u\|_{H^1}.
\]
The vertex values of $u_h-Pv$ sum to zero, so $\|u_h-Pv\|_{H^{1/2}}\le\|\partial_\xi(u_h-Pv)\|_{H^{-1/2}}$, and \cref{lem:bem-p1}(a) gives $\|u_h-Pv\|_{H^1}\le Ch^{-1/2}\|u_h-Pv\|_{H^{1/2}}\le C\|u\|_{H^1}$. By \cref{lem:bem-p1}(b), we also have $\|Pv\|_{H^1}\le\|v\|_{H^1}+|\bar v|\le C\|u\|_{H^1}$. Hence, we obtain
\[
\|u_h\|_{H^1}\le\|u_h-Pv\|_{H^1}+\|Pv\|_{H^1}\le C\|u\|_{H^1}\le C\|f\|_{L^2}.
\]

For a lumped load $y=M_vf$ with $f$ a P1 function, we have $\one^\T Py=0$, and for P1 functions $w$ whose vertex values sum to zero, $Pw=w$ and the symmetry of $P$ give $w^\T Py=w^\T M_vf$, and \eqref{eq:bem-lump} gives $|w^\T Py|\le C\|f\|_{L^2}\|w\|_{L^2}$. By \eqref{eq:bem-riesz}, $Py$ is the consistent load of a function $g$ with $\|g\|_{L^2}\le C\|f\|_{L^2}$. Since $\Lhp y=\Lhp Py$, the case of consistent loads gives $\|\Lhp y\|_{H^1}\le C\|g\|_{L^2}\le C\|f\|_{L^2}$.
\end{proof}

\subsection{Shape derivatives of the discrete Steklov--Poincar\'e matrix}\label{sec:bem-shape}

\begin{lemma}[{The first variation as a bilinear form on $H^1$}]\label{cor:bem-c1p-ab}
For $X\in\Nd$, a piecewise linear direction $e$ and P1 functions $u$ and $w$ whose vertex values sum to zero, $\one^\T u=\one^\T w=0$, we have
\[
\big|w^\T D_X\Lh(X)[e]\,u\big|\le C\,\|e'\|_{L^2}\,\|u\|_{H^1}\,\|w\|_{H^1}.
\]
\end{lemma}

\begin{proof}
Let $\varsigma_u=\tilde V^{-1}\tilde Q_{\rm int}u$ and $\varsigma_w$ be the Galerkin densities of $u$ and $w$. By \cref{lem:bem-reduction}, we have $\Lh=P(\cdots)P$, and since $Pu=u$ and $Pw=w$, the number $w^\T D_X\Lh[e]u$ is the derivative of the Schur expression inside the deflations, evaluated at $u$ and $w$. It is a sum of terms of the forms $\langle D\mathcal L_X[e]\partial_\xi u,\partial_\xi w\rangle$, which comes from the hypersingular part by Maue's formula, $\langle D\mathcal K_X[e]u,\varsigma_w\rangle$, $\langle D\mathcal K_X[e]w,\varsigma_u\rangle$ and $\langle D\mathcal L_X[e]\varsigma_u,\varsigma_w\rangle$. The bordered exterior part has the same structure, with the bordered exterior densities. By \eqref{eq:bem-tri} and $\|v\|_{L^1}\le\|v\|_{L^2}$, each term is bounded by $C\|e'\|_{L^2}$ times a product of one of $\|u'\|_{L^2}$, $\|\varsigma_u\|_{L^2}$ and one of $\|w'\|_{L^2}$, $\|\varsigma_w\|_{L^2}$. Hence, we obtain
\[
|w^\T D_X\Lh[e]u|\le C\|e'\|_{L^2}\big(\|u\|_{H^1}+\|\varsigma_u\|_{L^2}\big)\big(\|w\|_{H^1}+\|\varsigma_w\|_{L^2}\big),
\]
and \eqref{eq:bem-dens}, which holds for the interior and for the bordered exterior densities, bounds the right side by $C\|e'\|_{L^2}\|u\|_{H^1}\|w\|_{H^1}$.
\end{proof}

\begin{proof}[Proof of \cref{thm:bem-c1p}]
The vectors $z_f$ and $z_\varphi$ lie in $\one^\perp$ by \eqref{eq:mp-formula}. By \cref{cor:bem-c1p-ab} with $u=z_f$ and $w=z_\varphi$, we have
\[
|z_\varphi^\T D_X\Lh[e]z_f|\le C\|e'\|_{L^2}\|z_f\|_{H^1}\|z_\varphi\|_{H^1}.
\]
By \eqref{eq:bem-zf}, the factors on the right are bounded by $C\|f\|_{L^2}$ and $C\|\varphi\|_{L^2}$, which proves \eqref{eq:C1}.

For the increment, let $Y$ be a polygon of the segment $[X,X+e]$. Since $Y\in\Nd$, \eqref{eq:mp-deriv} applies at $Y$. Since $\Om$ is linear in $X$ and $\Lhp$ is symmetric, we obtain
\[
D_Yb_{Y,h}(f,\varphi)[e]=\underbrace{-2z_\varphi^\T D_Y\Lh(Y)[e]\,z_f}_{=:\,T_\Lambda}+\underbrace{2\big(D\Om[e]\varphi\big)^\T z_f}_{=:\,T_\varphi}+\underbrace{2z_\varphi^\T D\Om[e]f}_{=:\,T_f},
\]
where $z_f$ and $z_\varphi$ are taken at $Y$. The term $T_\Lambda$ obeys \eqref{eq:C1} at $Y$, and $\|f\|_{L^2}\le\|f\|_{L^\infty}$. By \eqref{eq:bem-dOm}, $D\Om[e]\varphi$ is the lumped load $M_vg$ of a P1 function $g$ with $\|g\|_{L^2}\le C\|e'\|_{L^2}\|\varphi\|_{L^\infty}$. Hence, the lumped pairing \eqref{eq:bem-lump} and \eqref{eq:bem-zf} bound $T_\varphi$ by $C\|e'\|_{L^2}\|f\|_{L^2}\|\varphi\|_{L^\infty}$. In the same way, $T_f$ is bounded by $C\|e'\|_{L^2}\|f\|_{L^\infty}\|\varphi\|_{L^2}$. The mean value inequality, applied to $t\mapsto b_{X+te,h}(f,\varphi)$ on $[0,1]$, gives the bound for the increment.
\end{proof}

\begin{proof}[Proof of \cref{thm:bem-c1}]
Let $M=M_1$ and $c=c_{\rm ca}$, so that $\mathfrak C\subset\mathfrak F(M,c)$. Let $X\in\mathfrak C$, let $e$ and $z$ be Lipschitz directions, and let $X_\varepsilon=X+\varepsilon_1e+\varepsilon_2z$ for $\varepsilon=(\varepsilon_1,\varepsilon_2)\in\R^2$. Since $|\Delta e|\le\|e'\|_{L^\infty}d(\xi,\eta)$ and $|\Delta X|\ge c\,d(\xi,\eta)$, the quotients relative to $X$ satisfy $|q_e|\le\|e'\|_{L^\infty}/c$ and $|q_z|\le\|z'\|_{L^\infty}/c$.

We first bound the second differences of the single layer. Let $\varsigma\in L^2$ be real, and let $s$ and $t$ be real numbers with
\begin{equation}\label{eq:bem-small}
|s|\,\|e'\|_{L^\infty}\le c/4,\qquad|t|\,\big(\|z'\|_{L^\infty}+1\big)\le c/12 .
\end{equation}
Since $|s|\,\|e'\|_{L^\infty}+|t|\,\|z'\|_{L^\infty}\le c/3$, the curves $X+se+tz$, $X+tz$, $X+se$ and $X$ lie in $\mathfrak F(M+c/2,c/2)$ by the bounds for the curves $W$ below. We put
\[
\Theta_{s,t}=\mathcal L_{X+se+tz}-\mathcal L_{X+tz}-\mathcal L_{X+se}+\mathcal L_X
\]
and show that
\begin{equation}\label{eq:bem-SD}
\|\Theta_{s,t}\varsigma\|_{H^1}\le C|s|\,|t|\,\|e'\|_{L^\infty}\|z'\|_{L^\infty}\|\varsigma\|_{L^2}.
\end{equation}

Let $\theta_1=sq_e$ and $\theta_2=tq_z$, so that $|\theta_1|\le\frac14$, $|\theta_2|\le\frac1{12}$ and $\Delta(X+se+tz)=\Delta X\,(1+\theta_1+\theta_2)$. In the kernel of $\Theta_{s,t}$, the terms $\log|\Delta X|$ cancel, and since $(1+\theta_1+\theta_2)/\big((1+\theta_1)(1+\theta_2)\big)=1-\omega$ with $\omega=\theta_1\theta_2/\big((1+\theta_1)(1+\theta_2)\big)$, the kernel is $-\frac1{2\pi}\log|1-\omega|$. Since $|1+\theta_1|\ge\frac34$ and $|1+\theta_2|\ge\frac34$, we have $|\omega|\le\frac{16}9|\theta_1|\,|\theta_2|\le\frac12$, and since $|\log|1-\omega||\le2|\omega|$, the kernel is bounded by $C|\theta_1|\,|\theta_2|$. Hence, for every $\xi$, we have
\begin{equation}\label{eq:bem-SD0}
|\Theta_{s,t}\varsigma(\xi)|\le C|s|\,|t|\,\|e'\|_{L^\infty}\|z'\|_{L^\infty}\|\varsigma\|_{L^1},
\end{equation}
and the same bound holds for $\|\Theta_{s,t}\varsigma\|_{L^2}$, since $\|\varsigma\|_{L^1}\le\|\varsigma\|_{L^2}$.

For the derivative, we let the parameters be complex. For $\sigma,\tau\in\C$ with $|\sigma|\,\|e'\|_{L^\infty}+|\tau|\,\|z'\|_{L^\infty}\le c/2$, the curve $W=X+\sigma e+\tau z$ lies in $\mathfrak F(M+c/2,c/2)$, since $|W'|\le M+|\sigma|\,\|e'\|_{L^\infty}+|\tau|\,\|z'\|_{L^\infty}\le M+c/2$ almost everywhere and $|W(\xi)-W(\eta)|\ge\big(c-|\sigma|\,\|e'\|_{L^\infty}-|\tau|\,\|z'\|_{L^\infty}\big)d(\xi,\eta)\ge(c/2)\,d(\xi,\eta)$. For such $W$ and $0<\delta<\frac12$, we put
\[
P^\delta_W\varsigma=W'\,\mathcal C^\delta_W(\varsigma/W'),\qquad\text{so that}\qquad P^\delta_W\varsigma(\xi)=\frac{W'(\xi)}{\pi\iu}\int_{d(\xi,\eta)>\delta}\frac{\varsigma(\eta)}{W(\eta)-W(\xi)}\,\dd\eta,
\]
and $\mathcal D^\delta_W\varsigma=\frac12\operatorname{Re}\big(\iu P^\delta_W\varsigma\big)$ if $\sigma$ and $\tau$ are real. The decomposition displayed before \eqref{eq:bem-Dlip}, with \eqref{eq:bem-david} and \cref{lem:bem-holo}, gives, uniformly in $\delta$,
\begin{equation}\label{eq:bem-Plip}
\|P^\delta_{W_1}\varsigma-P^\delta_{W_2}\varsigma\|_{L^2}\le C\|W_1'-W_2'\|_{L^\infty}\|\varsigma\|_{L^2}\qquad\big(W_1,W_2\in\mathfrak F(M+c/2,c/2)\big).
\end{equation}

Fix $\delta$, let $s$ and $t$ satisfy \eqref{eq:bem-small}, let $0<\kappa\le1$, and put $R=c/\big(12(\|z'\|_{L^\infty}+\kappa)\big)$, so that $|t|\le R$. Write $P_{\sigma,\tau}=P^\delta_{X+\sigma e+\tau z}\varsigma$, let $u=P_{s,t}-P_{0,t}-P_{s,0}+P_{0,0}$, and let
\[
\Psi(\tau)=\int_\Tor\overline{u(\xi)}\,\big(P_{s,\tau}-P_{0,\tau}\big)(\xi)\,\dd\xi\qquad(\tau\in\C,\ |\tau|<3R).
\]
For $|\tau|<3R$, we have $|s|\,\|e'\|_{L^\infty}+|\tau|\,\|z'\|_{L^\infty}<c/2$, so the denominators $W(\eta)-W(\xi)$ of the curves $X+se+\tau z$ and $X+\tau z$ are at least $c\delta/2$ in absolute value on $\{d(\xi,\eta)>\delta\}$. Hence, $\Psi$ is holomorphic in $|\tau|<3R$, as $\Phi$ in the proof of \cref{lem:bem-holo}. For $|\tau|\le2R$, both curves lie in $\mathfrak F(M+c/2,c/2)$ and their derivatives differ by $se'$, so the Cauchy--Schwarz inequality and \eqref{eq:bem-Plip} give $|\Psi(\tau)|\le B$ with $B=C|s|\,\|e'\|_{L^\infty}\|u\|_{L^2}\|\varsigma\|_{L^2}$. Cauchy's estimate on the circle of radius $R$ around each $\zeta$ with $|\zeta|\le R$ gives $|\Psi'(\zeta)|\le B/R$. Since $\Psi(t)-\Psi(0)=\|u\|_{L^2}^2$ and $|t|\le R$, the mean value inequality gives $\|u\|_{L^2}^2\le B|t|/R$, that is, $\|u\|_{L^2}\le C|s|\,|t|\,\|e'\|_{L^\infty}(\|z'\|_{L^\infty}+\kappa)\|\varsigma\|_{L^2}$. Letting $\kappa\to0$, we obtain
\[
\|u\|_{L^2}\le C|s|\,|t|\,\|e'\|_{L^\infty}\|z'\|_{L^\infty}\|\varsigma\|_{L^2}
\]
uniformly in $\delta$. The second difference $E_\delta=\mathcal D^\delta_{X+se+tz}\varsigma-\mathcal D^\delta_{X+tz}\varsigma-\mathcal D^\delta_{X+se}\varsigma+\mathcal D^\delta_X\varsigma$ equals $\frac12\operatorname{Re}(\iu u)$, so $\|E_\delta\|_{L^2}\le\frac12\|u\|_{L^2}$. By \eqref{eq:bem-dsL} for the four curves, $\int E_\delta\varphi\,\dd\xi$ tends to $-\int \Theta_{s,t}\varsigma\,\varphi'\,\dd\xi$ for every real periodic $C^1$ function $\varphi$, and the criterion at the beginning of the proof of \cref{lem:bem-H1} gives $\|\partial_\xi \Theta_{s,t}\varsigma\|_{L^2}\le C|s|\,|t|\,\|e'\|_{L^\infty}\|z'\|_{L^\infty}\|\varsigma\|_{L^2}$. Together with \eqref{eq:bem-SD0}, this proves \eqref{eq:bem-SD}.

The kernel of $\Theta_{s,t}$ is real, symmetric and bounded, so $\Theta_{s,t}$ is symmetric on $L^2$, and \eqref{eq:bem-SD} and \cref{lem:bem-fs}(a) give
\begin{equation}\label{eq:bem-SDh}
\|\Theta_{s,t}\varsigma\|_{H^{1/2}}\le C|s|\,|t|\,\|e'\|_{L^\infty}\|z'\|_{L^\infty}\|\varsigma\|_{H^{-1/2}} .
\end{equation}
Let $\varsigma,\psi\in L^2$ be real, and let $F(\varepsilon)=\langle\mathcal L_{X_\varepsilon}\varsigma,\psi\rangle$ be as in \cref{lem:bem-cont}(iii). For $s$ and $t$ as in \eqref{eq:bem-small}, we have $F(s,t)-F(0,t)-F(s,0)+F(0,0)=\langle \Theta_{s,t}\varsigma,\psi\rangle$, so \eqref{eq:bem-SDh} gives
\[
\big|F(s,t)-F(0,t)-F(s,0)+F(0,0)\big|\le L|s|\,|t|,\qquad L=C\|e'\|_{L^\infty}\|z'\|_{L^\infty}\|\varsigma\|_{H^{-1/2}}\|\psi\|_{H^{-1/2}} .
\]
By \cref{lem:bem-cont}(iii), $F$ is twice continuously differentiable near $0$. Dividing by $|t|$ and letting $t\to0$, and then dividing by $|s|$ and letting $s\to0$, we obtain
\begin{equation}\label{eq:bem-mixF}
|\partial_{\varepsilon_1}\partial_{\varepsilon_2}F(0)|\le C\|e'\|_{L^\infty}\|z'\|_{L^\infty}\|\varsigma\|_{H^{-1/2}}\|\psi\|_{H^{-1/2}} .
\end{equation}
For the double layer, let $a$ be a real Lipschitz function and $G(\varepsilon)=\langle\mathcal K_{X_\varepsilon}a,\psi\rangle$. Let $\tilde F$ be the function $F$ for the directions $-\iu e$ and $z$ and the density $a'$. Since $q_{-\iu e}=-\iu q_e$ and $\operatorname{Re}(-\iu\lambda)=\operatorname{Im}\lambda$ for $\lambda\in\C$, the formulas of \cref{lem:bem-cont}(iii) give $\partial_{\varepsilon_1}\partial_{\varepsilon_2}G(0)=-\partial_{\varepsilon_1}\partial_{\varepsilon_2}\tilde F(0)$. Since $\|(-\iu e)'\|_{L^\infty}=\|e'\|_{L^\infty}$ and $\|a'\|_{H^{-1/2}}\le C\|a\|_{H^{1/2}}$, \eqref{eq:bem-mixF} gives
\begin{equation}\label{eq:bem-mixG}
|\partial_{\varepsilon_1}\partial_{\varepsilon_2}G(0)|\le C\|e'\|_{L^\infty}\|z'\|_{L^\infty}\|a\|_{H^{1/2}}\|\psi\|_{H^{-1/2}} .
\end{equation}

We now turn to the Schur form. Let $v$ and $w$ be real Lipschitz functions and $\mathfrak s(\varepsilon)=\mathfrak s_{X_\varepsilon,h}(v,w)$. For $\beta\in\R^N$, let $\varsigma_\beta=\sum_p\beta_p\iota_p\in\mathcal P_0$, and for $g\in L^2$, let $\ell(g)=(\langle g,\iota_p\rangle)_p$. Let $V(\varepsilon)=\tilde V_{X_\varepsilon}$, let $\tilde C(\varepsilon)$ be its bordered inverse, and let $Q^\varepsilon_\pm=\frac12\pm\mathcal K_{X_\varepsilon}$. By \eqref{eq:bem-ext} and the matrix form of the Galerkin inverses in \cref{lem:bem-galinv}(a), we have
\[
\mathfrak s(\varepsilon)=2\underbrace{\langle\mathcal L_{X_\varepsilon}\partial_\xi v,\partial_\xi w\rangle}_{=:\,H(\varepsilon)}+\underbrace{\ell(Q^\varepsilon_+w)^\T V(\varepsilon)^{-1}\ell(Q^\varepsilon_+v)}_{=:\,S_+(\varepsilon)}+\underbrace{\ell(Q^\varepsilon_-w)^\T\tilde C(\varepsilon)\,\ell(Q^\varepsilon_-v)}_{=:\,S_-(\varepsilon)} .
\]
The function $H$ is the function $F$ with $\varsigma=\partial_\xi v$ and $\psi=\partial_\xi w$, the numbers $\beta^\T V(\varepsilon)\alpha=\langle\mathcal L_{X_\varepsilon}\varsigma_\alpha,\varsigma_\beta\rangle$ are functions $F$, and the numbers $\beta^\T\ell(Q^\varepsilon_\pm v)=\frac12\langle v,\varsigma_\beta\rangle\pm\langle\mathcal K_{X_\varepsilon}v,\varsigma_\beta\rangle$ are functions $G$ up to a constant. By \cref{lem:bem-cont}(iii), $H$, $V$ and the vectors $\ell(Q^\varepsilon_\pm v)$, $\ell(Q^\varepsilon_\pm w)$ are twice continuously differentiable near $0$. The bordered inverse $\tilde C(\varepsilon)$ is the upper left $N\times N$ block of the inverse of the matrix $\mathsf B(V(\varepsilon))$, which borders $V(\varepsilon)$ by the column $h\one$ and the row $h\one^\T$ of the constraint $\int\varsigma_\beta\,\dd\xi=h\sum_p\beta_p=0$. The matrix $V(0)$ is positive definite by \cref{lem:bem-Vreg}(a), so $\mathsf B(V(0))$ is invertible, since the bordered system is uniquely solvable by \cref{lem:bem-galinv}(a). Since inversion is smooth on the open set of invertible matrices, and the derivative of $\mathsf B(V)$ in the direction of a matrix $A$ is $A$ in the upper left corner, $V^{-1}$ and $\tilde C$ are smooth functions of $V$ near $V(0)$, with the derivatives $-V^{-1}AV^{-1}$ and $-\tilde CA\tilde C$. Hence, $\mathfrak s$ is twice continuously differentiable near $0$.

For the bounds, we say that a vector $r\in\R^N$ has dual size $B$ if $|\beta^\T r|\le B\|\varsigma_\beta\|_{H^{-1/2}}$ for all $\beta$, and that a matrix $A$ has energy size $K$ if $|\beta^\T A\alpha|\le K\|\varsigma_\alpha\|_{H^{-1/2}}\|\varsigma_\beta\|_{H^{-1/2}}$ for all $\alpha$ and $\beta$. Let $c_0$ be the constant of \cref{lem:bem-Vreg}(a), and let $\mathsf J$ be $V(0)^{-1}$ or $\tilde C(0)$. If $r$ has dual size $B$ and $y=\mathsf J r$, then $y^\T V(0)y=y^\T r$, also in the bordered case, where $\sum_py_p=0$. Hence, $c_0\|\varsigma_y\|_{H^{-1/2}}^2\le y^\T r\le B\|\varsigma_y\|_{H^{-1/2}}$, and $\|\varsigma_y\|_{H^{-1/2}}\le B/c_0$. Since $\mathsf J$ is symmetric, it follows that, for vectors $r_w$, $r_v$ of dual sizes $B_w$, $B_v$ and matrices $A$, $A'$ of energy sizes $K$, $K'$, we have
\begin{equation}\label{eq:bem-types}
|r_w^\T\mathsf J r_v|\le\frac{B_wB_v}{c_0},\qquad|r_w^\T\mathsf J A\mathsf J r_v|\le\frac{B_wKB_v}{c_0^2},\qquad|r_w^\T\mathsf J A\mathsf J A'\mathsf J r_v|\le\frac{B_wKK'B_v}{c_0^3} .
\end{equation}
At $\varepsilon=0$, the vector $\ell(Q^0_\pm v)$ has dual size $C\|v\|_{H^{1/2}}$, since $|\beta^\T\ell(Q^0_\pm v)|=|\langle Q^0_\pm v,\varsigma_\beta\rangle|\le\|Q^0_\pm v\|_{H^{1/2}}\|\varsigma_\beta\|_{H^{-1/2}}$ and $\|Q^0_\pm v\|_{H^{1/2}}\le C\|v\|_{H^{1/2}}$ by \cref{lem:bem-Vreg}(c). By \cref{lem:bem-cont}(iii) and \eqref{eq:bem-mixG}, its first derivatives in $\varepsilon_1$ and $\varepsilon_2$ have dual sizes $C\|e'\|_{L^\infty}\|v\|_{H^{1/2}}$ and $C\|z'\|_{L^\infty}\|v\|_{H^{1/2}}$, and its mixed second derivative has dual size $C\|e'\|_{L^\infty}\|z'\|_{L^\infty}\|v\|_{H^{1/2}}$. The same holds for $w$. By \cref{lem:bem-cont}(iii) and \eqref{eq:bem-mixF}, the first derivatives of $V$ at $0$ have energy sizes $C\|e'\|_{L^\infty}$ and $C\|z'\|_{L^\infty}$, and the mixed second derivative has energy size $C\|e'\|_{L^\infty}\|z'\|_{L^\infty}$. With $V_i=\partial_{\varepsilon_i}V(0)$ and $V_{12}=\partial_{\varepsilon_1}\partial_{\varepsilon_2}V(0)$, the derivatives of $\mathsf J$ along $V(\varepsilon)$ at $0$ are $-\mathsf J V_i\mathsf J$ and $\mathsf J V_1\mathsf J V_2\mathsf J+\mathsf J V_2\mathsf J V_1\mathsf J-\mathsf J V_{12}\mathsf J$. By the product rule, $\partial_{\varepsilon_1}S_\pm(0)$ and $\partial_{\varepsilon_1}\partial_{\varepsilon_2}S_\pm(0)$ are sums of terms of the three forms in \eqref{eq:bem-types}, in which $r_w$ and $r_v$ are $\ell(Q^0_\pm w)$, $\ell(Q^0_\pm v)$ or their derivatives, and $A$, $A'$ are derivatives of $V$. Each term of $\partial_{\varepsilon_1}S_\pm(0)$ contains one derivative in $\varepsilon_1$, and each term of $\partial_{\varepsilon_1}\partial_{\varepsilon_2}S_\pm(0)$ contains one derivative in $\varepsilon_1$ and one in $\varepsilon_2$. Hence, \eqref{eq:bem-types} and the sizes above give
\begin{equation}\label{eq:bem-Sder}
|\partial_{\varepsilon_1}S_\pm(0)|\le C\|e'\|_{L^\infty}\|v\|_{H^{1/2}}\|w\|_{H^{1/2}},\qquad|\partial_{\varepsilon_1}\partial_{\varepsilon_2}S_\pm(0)|\le C\|e'\|_{L^\infty}\|z'\|_{L^\infty}\|v\|_{H^{1/2}}\|w\|_{H^{1/2}} .
\end{equation}
For $H$, \cref{lem:bem-cont}(iii), \eqref{eq:bem-mixF} and $\|\partial_\xi v\|_{H^{-1/2}}\le C\|v\|_{H^{1/2}}$ give
\begin{equation}\label{eq:bem-Hder}
|\partial_{\varepsilon_1}H(0)|\le C\|e'\|_{L^\infty}\|v\|_{H^{1/2}}\|w\|_{H^{1/2}},\qquad|\partial_{\varepsilon_1}\partial_{\varepsilon_2}H(0)|\le C\|e'\|_{L^\infty}\|z'\|_{L^\infty}\|v\|_{H^{1/2}}\|w\|_{H^{1/2}} .
\end{equation}
Adding \eqref{eq:bem-Hder}, twice, and \eqref{eq:bem-Sder} for both signs gives the bounds of the first and the mixed second derivative.

For the Lipschitz bound and the Taylor bound, let $Y\in\mathfrak C$, let $d=Y-X$, and let $\phi(t)=\mathfrak s_{X+td,h}(v,w)$. Since $\mathfrak C$ is convex, the curves $X+td$ lie in $\mathfrak C$ for $0\le t\le1$, and $\|d'\|_{L^\infty}=\|X'-Y'\|_{L^\infty}$. Since $\phi(t+\varepsilon_1+\varepsilon_2)=\mathfrak s_{X+td+\varepsilon_1d+\varepsilon_2d,h}(v,w)$, the above, applied at $X+td$ with both directions equal to $d$, shows that $\phi$ is twice continuously differentiable on $[0,1]$, with
\[
|\phi'(t)|\le C\|d'\|_{L^\infty}\|v\|_{H^{1/2}}\|w\|_{H^{1/2}},\qquad|\phi''(t)|\le C\|d'\|_{L^\infty}^2\|v\|_{H^{1/2}}\|w\|_{H^{1/2}} .
\]
The mean value inequality for $\phi$ gives the Lipschitz bound, and Taylor's formula gives $|\phi(1)-\phi(0)-\phi'(0)|\le\frac12\sup_{[0,1]}|\phi''|$, which is the Taylor bound, since $\phi'(0)=D\mathfrak s_{X,h}[Y-X](v,w)$.
\end{proof}

\begin{proof}[Proof of \cref{thm:bem-c3p}]
We first bound the second derivative of $Y\mapsto b_{Y,h}(f,\varphi)$ at $X$ in the direction $e$. By the paragraphs after \cref{def:nbhd}, \eqref{eq:mp-deriv} holds on a neighborhood of $X$, and $\Lh$ is twice differentiable at $X$. Write $R=\Lhp(X)$, $L_1=D\Lh(X)[e]$ and $L_2=D^2\Lh(X)[e,e]$. Differentiating \eqref{eq:mp-deriv} once more in the direction $e$, we obtain
\[
D^2\Lhp(X)[e,e]=R\,(2L_1RL_1-L_2)\,R .
\]
Since $\Om$ is linear in $X$, we have $D^2\Om=0$. Put $z_f=R\,\Om f$, $z_\varphi=R\,\Om\varphi$, $\zeta_f=R\,D\Om[e]f$ and $\zeta_\varphi=R\,D\Om[e]\varphi$. Differentiating $b_{X,h}(f,\varphi)=2(\Om\varphi)^\T R\,\Om f$ twice and using the symmetry of $R$, we obtain
\begin{align*}
D^2b_{X,h}(f,\varphi)[e,e]&=\underbrace{-2z_\varphi^\T L_2z_f}_{=:\,\mathrm{(A)}}+\underbrace{4z_\varphi^\T L_1RL_1z_f}_{=:\,\mathrm{(B)}}\\
&\quad+\underbrace{\underbrace{\big(-4\zeta_\varphi^\T L_1z_f\big)}_{=:\,\mathrm{(C_1)}}+\underbrace{\big(-4z_\varphi^\T L_1\zeta_f\big)}_{=:\,\mathrm{(C_2)}}+\underbrace{4\big(D\Om[e]\varphi\big)^\T \zeta_f}_{=:\,\mathrm{(C_3)}}}_{=:\,\mathrm{(C)}} .
\end{align*}
By \eqref{eq:bem-zf}, we have $\|z_f\|_{H^1}\le C\|f\|_{L^2}\le C\|f\|_{L^\infty}$ and $\|z_\varphi\|_{H^1}\le C\|\varphi\|_{L^2}$.

For (A), let $u$ and $w$ be P1 traces in $\one^\perp$, let $\varsigma_u=\tilde V^{-1}\tilde Q_{\rm int}u$ and $\tau_u=\tilde C\tilde Q_{\rm ext}u$ be the interior and the bordered exterior Galerkin densities of $u$, and define $\varsigma_w$ and $\tau_w$ in the same way. By \cref{lem:bem-reduction}, $\Lh$ is the function $P(\cdots)P$ of $\tilde V$ and $\tilde K$ near $X$. We differentiate it twice in the direction $e$, using $D(\tilde V^{-1})=-\tilde V^{-1}D\tilde V\,\tilde V^{-1}$, $D\tilde C=-\tilde C\,D\tilde V\,\tilde C$ and the symmetry of $\tilde V^{-1}$, $\tilde C$ and $D\tilde V[e]$. Then, we obtain
\[
w^\T L_2u=\underbrace{\mathcal A(u,w)}_{=:\,\mathrm{(A_1)}}+\underbrace{2\rho_w^\T\tilde V^{-1}\rho_u+2\hat\rho_w^\T\tilde C\hat\rho_u}_{=:\,\mathrm{(A_2)}},
\]
where $\rho_u=D\tilde K[e]u-D\tilde V[e]\varsigma_u$ and $\hat\rho_u=D\tilde K[e]u+D\tilde V[e]\tau_u$, and likewise for $w$. Here, $\mathcal A(u,w)$ is the sum of the terms of the forms listed in the proof of \cref{cor:bem-c1p-ab}, with $D^2\mathcal L_X[e,e]$ and $D^2\mathcal K_X[e,e]$ in place of the first derivatives.

In $\mathrm{(A_1)}$, the kernels of $D^2\mathcal L_X[e,e]$ and $D^2\mathcal K_X[e,e]$ are bounded by $C|q_e|^2$ by \cref{lem:bem-log}. Applying \eqref{eq:maximal} to both factors of $|q_e|^2$, and $Me'\le\|e'\|_{L^\infty}$ in one of them, we obtain $|q_e(\xi,\eta)|^2\le C\|e'\|_{L^\infty}Me'(\eta)$. Hence, for real $\varsigma\in L^2$ and every density $\psi$ of $\mathcal P_0$, the Cauchy--Schwarz inequality and $\|Me'\|_{L^2}\le C\|e'\|_{L^2}$ give
\[
|\langle D^2\mathcal L_X[e,e]\varsigma,\psi\rangle|\le C\|e'\|_{L^\infty}\int_\Tor|\psi|\,\dd\xi\int_\Tor Me'\,|\varsigma|\,\dd\eta\le C\|e'\|_{L^\infty}\|e'\|_{L^2}\|\varsigma\|_{L^2}\|\psi\|_{L^2}.
\]
The same bound holds for $\langle D^2\mathcal K_X[e,e]u,\psi\rangle$, with $\|u'\|_{L^2}$ in place of $\|\varsigma\|_{L^2}$. In each term of $\mathcal A(u,w)$, the second factor of the pairing is the derivative of a P1 function or a Galerkin density, hence a density of $\mathcal P_0$, so these bounds and \eqref{eq:bem-dens} give $|\mathrm{(A_1)}|\le C\|e'\|_{L^\infty}\|e'\|_{L^2}\|u\|_{H^1}\|w\|_{H^1}$.

For $\mathrm{(A_2)}$ and for (B), we use that the first derivatives of the layer operators map into $H^1$. For $\lambda\in\C$ with $|\lambda|\le1$ and real $\varsigma\in L^2$, let
\[
Q_\lambda\varsigma(\xi)=\int_\Tor\operatorname{Re}\big(\lambda\,q_e(\xi,\eta)\big)\,\varsigma(\eta)\,\dd\eta .
\]
By \cref{lem:bem-log}, we have $D\mathcal L_X[e]\varsigma=-\frac1{2\pi}Q_1\varsigma$, and $D\mathcal K_X[e]u=\frac1{2\pi}Q_{-\iu}u'$ for a P1 trace $u$. We show that
\begin{equation}\label{eq:c3-H1}
\|Q_\lambda\varsigma\|_{H^1}\le C\|e'\|_{L^\infty}\|\varsigma\|_{L^2}.
\end{equation}
Let $c=c_{\rm ca}$ and $s=\|e'\|_{L^\infty}$, so that $X\in\mathfrak F(M_1,c)$ and $|q_e|\le s/c$. If $s=0$, then $q_e=0$, and there is nothing to prove, so let $s>0$. For $0<t\le c/(2s)$, let $X_t=X+t\lambda e$ and
\[
u_t=\frac{(\mathcal L_{X_t}-\mathcal L_X)\varsigma}{t} .
\]
We have $|X_t'|\le M_1+c/2$ and $|\Delta X_t|\ge(c-ts)\,d(\xi,\eta)\ge\frac c2\,d(\xi,\eta)$, so $X$ and $X_t$ lie in $\mathfrak F(M_1+c/2,c/2)$. Since $\|X_t'-X'\|_{L^\infty}\le ts$, \cref{lem:bem-cont}(i) for this class gives $\|u_t\|_{H^1}\le Cs\|\varsigma\|_{L^2}$. As in the proof of \cref{lem:bem-cont}(iii), $\Delta X_t=\Delta X\,(1+t\lambda q_e)$ and $|t\lambda q_e|\le\frac12$ give $\log|\Delta X_t|=\log|\Delta X|+\operatorname{Re}\operatorname{Log}(1+t\lambda q_e)$, and hence
\[
u_t(\xi)=-\frac1{2\pi}\int_\Tor\frac{\operatorname{Re}\operatorname{Log}\big(1+t\lambda q_e(\xi,\eta)\big)}{t}\,\varsigma(\eta)\,\dd\eta
\]
for every $\xi$. Since $|\operatorname{Log}(1+z)|\le\frac32|z|$ for $|z|\le\frac12$, the integrand is bounded by $(3s/2c)\,|\varsigma(\eta)|$, and it tends to $\operatorname{Re}\big(\lambda q_e(\xi,\eta)\big)\varsigma(\eta)$ as $t\to0$. By dominated convergence, $u_t(\xi)\to-\frac1{2\pi}Q_\lambda\varsigma(\xi)$ for every $\xi$, and $|u_t(\xi)|\le(s/c)\|\varsigma\|_{L^1}$ for all $t$ and $\xi$. Hence, by dominated convergence on $\Tor$, the Fourier coefficients of $u_t$ converge to those of $-\frac1{2\pi}Q_\lambda\varsigma$. Letting $t\to0$ in each finite part of the series for $(2\pi)^2\|u_t\|_{H^1}^2\le(2\pi Cs\|\varsigma\|_{L^2})^2$, we obtain the same bound for each finite part of the series for $\|Q_\lambda\varsigma\|_{H^1}^2$. Thus, $Q_\lambda\varsigma$ lies in $H^1$, and \eqref{eq:c3-H1} holds.

By \eqref{eq:c3-H1}, we have $\|D\mathcal L_X[e]\varsigma\|_{H^1}\le C\|e'\|_{L^\infty}\|\varsigma\|_{L^2}$ and $\|D\mathcal K_X[e]u\|_{H^1}\le C\|e'\|_{L^\infty}\|u'\|_{L^2}$. In $\mathrm{(A_2)}$, the vector $\rho_u$ is, by \cref{lem:bem-log}, the load $(\langle g_u,\iota_p\rangle)_p$ of $g_u=D\mathcal K_X[e]u-D\mathcal L_X[e]\varsigma_u$, where $\varsigma_u$ is read as a density of $\mathcal P_0$, and \eqref{eq:bem-dens} gives $\|g_u\|_{H^1}\le C\|e'\|_{L^\infty}\|u\|_{H^1}$. The vector $\vartheta_u=\tilde V^{-1}\rho_u$ is the Galerkin density of $g_u$, so $\|\vartheta_u\|_{L^2}\le C\|e'\|_{L^\infty}\|u\|_{H^1}$ by \cref{lem:bem-galinv}(d). The remaining factor $\rho_w^\T\vartheta_u=\langle D\mathcal K_X[e]w-D\mathcal L_X[e]\varsigma_w,\vartheta_u\rangle$ is a first variation, which \eqref{eq:bem-tri}, $\|v\|_{L^1}\le\|v\|_{L^2}$ and \eqref{eq:bem-dens} bound by $C\|e'\|_{L^2}\|w\|_{H^1}\|\vartheta_u\|_{L^2}$. The bordered term is treated in the same way, with $\hat\rho_u$, $\tau_u$ and $\tilde C\hat\rho_u$. Hence, we have $|\mathrm{(A_2)}|\le C\|e'\|_{L^\infty}\|e'\|_{L^2}\|u\|_{H^1}\|w\|_{H^1}$. With $u=z_f$ and $w=z_\varphi$, which lie in $\one^\perp$ by \eqref{eq:mp-formula}, the bounds for $\mathrm{(A_1)}$, $\mathrm{(A_2)}$, $z_f$ and $z_\varphi$ give
\begin{equation}\label{eq:c3-A}
|\mathrm{(A)}|\le C\|e'\|_{L^\infty}\|e'\|_{L^2}\|f\|_{L^\infty}\|\varphi\|_{L^2}.
\end{equation}

For (B), we first show that $|w^\T L_1u|\le C\|e'\|_{L^\infty}\|u\|_{H^1}\|w\|_{L^2}$ for P1 traces $u,w\in\one^\perp$, with the densities of (A). By \cref{lem:bem-reduction}, $w^\T L_1u$ is the sum of a hypersingular term $2\langle D\mathcal L_X[e]\partial_\xi u,\partial_\xi w\rangle$, of terms in which the derivative falls on $u$ or on its densities and which are paired with the densities of $w$, and of the transposed double layer terms $\langle D\mathcal K_X[e]w,\varsigma_u\rangle$ and $\langle D\mathcal K_X[e]w,\tau_u\rangle$. Since $D\mathcal L_X[e]\partial_\xi u\in H^1$ by \eqref{eq:c3-H1}, integration by parts turns the hypersingular term into $-2\langle\partial_\xi D\mathcal L_X[e]\partial_\xi u,w\rangle$, which is bounded by $C\|e'\|_{L^\infty}\|u'\|_{L^2}\|w\|_{L^2}$. By the symmetry of the Galerkin inverses (\cref{lem:bem-galinv}(a)), the terms with the densities of $w$ are pairings $\langle(\tfrac12\pm\mathcal K_X)w,\rho\rangle$, where $\rho$ is the interior or bordered Galerkin density of one of $D\mathcal K_X[e]u$, $D\mathcal L_X[e]\varsigma_u$ and $D\mathcal L_X[e]\tau_u$. These functions have norm in $H^1$ bounded by $C\|e'\|_{L^\infty}\|u\|_{H^1}$ by \eqref{eq:c3-H1} and \eqref{eq:bem-dens}, so $\|\rho\|_{L^2}\le C\|e'\|_{L^\infty}\|u\|_{H^1}$ by \cref{lem:bem-galinv}(d). Since $\mathcal K_X$ is bounded on $L^2$ by \cref{lem:bem-Vreg}(c), these terms are bounded by $C\|e'\|_{L^\infty}\|u\|_{H^1}\|w\|_{L^2}$. For the transposed terms, let $\beta$ be $\varsigma_u$ or $\tau_u$. By \cref{lem:bem-log}, $D\mathcal K_X[e]$ acts on $w'$ with the bounded kernel $\frac1{2\pi}\operatorname{Im}q_e$. Since $q_e$ is symmetric in $\xi$ and $\eta$, Fubini's theorem gives $\langle D\mathcal K_X[e]w,\beta\rangle=\langle w',G\rangle$ with $G=\frac1{2\pi}Q_{-\iu}\beta$. Integrating by parts, \eqref{eq:c3-H1} and \eqref{eq:bem-dens} give $|\langle w',G\rangle|=|\langle w,G'\rangle|\le C\|e'\|_{L^\infty}\|w\|_{L^2}\|\beta\|_{L^2}\le C\|e'\|_{L^\infty}\|u\|_{H^1}\|w\|_{L^2}$. Adding the bounds of all terms proves the bound for $w^\T L_1u$.

Next, let $y=L_1z_f$. Since $L_1=PL_1P$ by \cref{lem:bem-reduction}, we have $\one^\T y=0$, and the bound above gives $|w^\T y|\le C\|e'\|_{L^\infty}\|z_f\|_{H^1}\|w\|_{L^2}$ for all P1 functions $w\in\one^\perp$. By \eqref{eq:bem-riesz} and the bound for $z_f$, $y$ is the consistent load of a function $F$ with $\|F\|_{L^2}\le C\|e'\|_{L^\infty}\|f\|_{L^\infty}$, and \cref{lem:bem-areg} gives $\|RL_1z_f\|_{H^1}\le C\|e'\|_{L^\infty}\|f\|_{L^\infty}$. Since $R=\Lhp$ maps into $\one^\perp$ by \eqref{eq:mp-formula}, \cref{cor:bem-c1p-ab} applies to $RL_1z_f$ and $z_\varphi$, and gives
\begin{equation}\label{eq:c3-B}
|\mathrm{(B)}|\le C\|e'\|_{L^\infty}\|e'\|_{L^2}\|f\|_{L^\infty}\|\varphi\|_{L^2}.
\end{equation}

For (C), \eqref{eq:bem-dOm} gives $D\Om[e]f=M_vg_f$ and $D\Om[e]\varphi=M_vg_\varphi$ with P1 functions $g_f$ and $g_\varphi$ such that
\[
\|g_f\|_{L^2}\le C\|e'\|_{L^\infty}\|f\|_{L^2},\qquad\|g_f\|_{L^2}\le C\|e'\|_{L^2}\|f\|_{L^\infty},\qquad\|g_\varphi\|_{L^2}\le C\|e'\|_{L^\infty}\|\varphi\|_{L^2}.
\]
By \cref{lem:bem-areg} for lumped loads and $\|f\|_{L^2}\le\|f\|_{L^\infty}$, we obtain $\|\zeta_\varphi\|_{H^1}\le C\|e'\|_{L^\infty}\|\varphi\|_{L^2}$ and the two bounds $\|\zeta_f\|_{H^1}\le C\|e'\|_{L^\infty}\|f\|_{L^\infty}$ and $\|\zeta_f\|_{H^1}\le C\|e'\|_{L^2}\|f\|_{L^\infty}$. The vectors $\zeta_f$ and $\zeta_\varphi$ lie in $\one^\perp$, since $R$ maps into $\one^\perp$. In $\mathrm{(C_1)}$, \cref{cor:bem-c1p-ab} with the bounds for $z_f$ and $\zeta_\varphi$, and in $\mathrm{(C_2)}$, \cref{cor:bem-c1p-ab} with the bound for $z_\varphi$ and the first bound for $\zeta_f$, give the bound $C\|e'\|_{L^\infty}\|e'\|_{L^2}\|f\|_{L^\infty}\|\varphi\|_{L^2}$. In $\mathrm{(C_3)}$, the lumped pairing \eqref{eq:bem-lump} gives $|(M_vg_\varphi)^\T \zeta_f|\le C\|g_\varphi\|_{L^2}\|\zeta_f\|_{L^2}$, and the second bound for $\zeta_f$ gives the same bound. Hence, we obtain
\begin{equation}\label{eq:c3-C}
|\mathrm{(C)}|\le C\|e'\|_{L^\infty}\|e'\|_{L^2}\|f\|_{L^\infty}\|\varphi\|_{L^2}.
\end{equation}

Adding \eqref{eq:c3-A}, \eqref{eq:c3-B} and \eqref{eq:c3-C}, we obtain
\[
|D^2b_{X,h}(f,\varphi)[e,e]|\le K,\qquad K=C\|e'\|_{L^\infty}\|e'\|_{L^2}\|f\|_{L^\infty}\|\varphi\|_{L^2} .
\]
For $0\le t<1$, the segment from $X+te$ in the direction $(1-t)e$ is $[X+te,X+e]\subset\Nd$, so this bound holds at $X+te$ for $(1-t)e$, and since both sides are homogeneous of degree two in the direction, it gives $|D^2b_{X+te,h}(f,\varphi)[e,e]|\le K$. For $t=1$, the bound at $X+e$ for the direction $-e$ gives the same. Let $\gamma(t)=b_{X+te,h}(f,\varphi)$ for $0\le t\le1$. Then, $\gamma'(t)=Db_{X+te,h}(f,\varphi)[e]$ is differentiable on $[0,1]$ with $|\gamma''|\le K$, and the mean value inequality, applied to $\gamma'$ and then to $t\mapsto\gamma(t)-t\gamma'(0)$, gives
\[
|\delta^2\bXh(f,\varphi)|=|\gamma(1)-\gamma(0)-\gamma'(0)|\le K .
\]
\end{proof}

\subsection{Consistency}\label{sec:bem-cons}

\begin{lemma}[The polygon and the curve]\label{lem:bem-geo}
Assume \ref{K1}, \ref{K2} with $r=1$, and \ref{K3}, and let $X=\Ih x$ with $h\le h_0$. Then, for real $\varsigma\in L^2$ and real Lipschitz functions $a$, we have
\begin{equation}\label{eq:bem-geo}
\|(\mathcal L_X-\mathcal L_x)\varsigma\|_{L^2}\le Ch^2\|\varsigma\|_{L^2},\qquad\|(\mathcal K_X-\mathcal K_x)a\|_{L^2}\le Ch^2\|a'\|_{L^2},
\end{equation}
and
\begin{equation}\label{eq:bweak-LK}
\|(\mathcal L_X-\mathcal L_x)\varsigma\|_{H^{1/2}}\le C\epsh\|\varsigma\|_{L^2},\qquad\|(\mathcal K_X-\mathcal K_x)a\|_{H^{1/2}}\le C\epsh\|a\|_{H^1} .
\end{equation}
Moreover, for real $\varsigma\in L^2$ and real Lipschitz functions $w$, we have
\begin{equation}\label{eq:bweak-pair}
|\langle\varsigma,(\mathcal K_X-\mathcal K_x)w\rangle|\le C\epsh\|\varsigma\|_{L^2}\|w\|_{H^{1/2}} .
\end{equation}
\end{lemma}

\begin{proof}
Let $\theta=X-x$. On an interval of length $h$, the difference between a real $W^{2,\infty}$ function and its linear interpolant is bounded by $\frac{h^2}8$ times the bound of the second derivative of the function, so \ref{K2} with $r=1$ gives $\|\theta\|_{L^\infty}\le Ch^2$. Since $X=\Ih x$ and $h\le h_0$, we also have $\|\theta'\|_{L^\infty}\le C_Ih\le r_1$, and $X\in\Nd\subset\mathfrak C$. Let $q=\Delta\theta/\Delta x$, so that $\Delta X/\Delta x=1+q$. Since $|\Delta\theta|\le\min(\|\theta'\|_{L^\infty}d(\xi,\eta),2\|\theta\|_{L^\infty})$ and $|\Delta x|\ge c_*d(\xi,\eta)$, we have $|q|\le C\min(h,h^2/d(\xi,\eta))$, and moreover $|q|\le\|\theta'\|_{L^\infty}/c_*\le r_1/c_*\le\tfrac14$. Let $\operatorname{Log}$ be the principal branch of the logarithm. By the definition of $\mathcal L_Y$, the kernel of $\mathcal L_X-\mathcal L_x$ is $-\frac1{2\pi}\big(\log|\Delta X|-\log|\Delta x|\big)=-\frac1{2\pi}\operatorname{Re}\operatorname{Log}(1+q)$. For the double layer, let $a$ be a real Lipschitz function. Since $\|\theta'\|_{L^\infty}\le r_1\le c_*/4<c_*/2$, the argument for $G$ in the proof of \cref{lem:bem-cont}(iii), applied to the curve $x\in\mathfrak F(M_*,c_*)$ in the direction $\theta$ at $\varepsilon=(1,0)$, gives, for almost every $\xi$,
\[
(\mathcal K_X-\mathcal K_x)a(\xi)=\frac1{2\pi}\int_\Tor\operatorname{Im}\operatorname{Log}\big(1+q(\xi,\eta)\big)\,a'(\eta)\,\dd\eta .
\]
Thus, both differences are integral operators with the bounded kernels $-\frac1{2\pi}\operatorname{Re}\operatorname{Log}(1+q)$ and $\frac1{2\pi}\operatorname{Im}\operatorname{Log}(1+q)$, applied to $f=\varsigma$ and to $f=a'$, respectively. We split $\operatorname{Log}(1+q)=q+\big(\operatorname{Log}(1+q)-q\big)$.

Since $|q|\le\tfrac14$, we have
\[
|\operatorname{Log}(1+q)-q|\le|q|^2\le C\min\Big(h^2,\frac{h^4}{d(\xi,\eta)^2}\Big)\le C\,(h^2)^{7/8}\Big(\frac{h^4}{d(\xi,\eta)^2}\Big)^{1/8}=C\,h^{9/4}d(\xi,\eta)^{-1/4} .
\]
Since $\int_\Tor d(\xi,\eta)^{-1/2}\,\dd\eta\le4$ for every $\xi$, the Cauchy--Schwarz inequality bounds the contribution of $\operatorname{Log}(1+q)-q$ at every $\xi$ by $Ch^{9/4}\|f\|_{L^2}$, and hence in $L^2$ by $Ch^2\|f\|_{L^2}$.

The part with $q$ is the limit of commutators with the truncated Cauchy integrals. For $0<\delta<\frac12$ and $f\in L^2(\Tor;\C)$, we have
\[
\int_{d(\xi,\eta)>\delta}q(\xi,\eta)f(\eta)\,\dd\eta=-\pi\iu\Big(\theta(\xi)\,\mathcal C^\delta_x(f/x')(\xi)-\mathcal C^\delta_x(\theta f/x')(\xi)\Big),
\]
and by \eqref{eq:bem-david} for $\mathfrak F(M_*,c_*)$ and $|x'|\ge c_*$, the right side has norm in $L^2$ bounded by $C\|\theta\|_{L^\infty}\|f\|_{L^2}$, uniformly in $\delta$. Since $|q|\le\frac14$, the truncated integrals converge to $\int_\Tor q(\xi,\eta)f(\eta)\,\dd\eta$ for every $\xi$ as $\delta\to0$ by dominated convergence, and Fatou's lemma gives the same bound for the limit. Hence, the part with $q$ of both kernels is bounded on $L^2$ by $C\|\theta\|_{L^\infty}\le Ch^2$. Together with the bound of the remainder, this proves \eqref{eq:bem-geo}.

By \cref{lem:bem-cont}(i) and (ii) for $\mathfrak F(M_1,c_{\rm ca})\supset\mathfrak C$, we have $\|(\mathcal L_X-\mathcal L_x)\varsigma\|_{H^1}\le C\|X'-x'\|_{L^\infty}\|\varsigma\|_{L^2}\le Ch\|\varsigma\|_{L^2}$ and $\|(\mathcal K_X-\mathcal K_x)a\|_{H^1}\le Ch\|a\|_{H^1}$. With $\|a'\|_{L^2}\le C\|a\|_{H^1}$, the interpolation inequality \eqref{eq:bem-interp} turns \eqref{eq:bem-geo} and the bounds in $H^1$ into \eqref{eq:bweak-LK}.

For \eqref{eq:bweak-pair}, let $v$ be the Fourier partial sum of $w$ of degree $1/h$, and split $w=v+(w-v)$ into Lipschitz functions. Comparing Fourier coefficients, we obtain
\[
\|v\|_{H^1}\le Ch^{-1/2}\|w\|_{H^{1/2}},\qquad\|w-v\|_{L^2}\le h^{1/2}\|w\|_{H^{1/2}}.
\]
For the low frequencies, the second bound of \eqref{eq:bem-geo} gives $\|(\mathcal K_X-\mathcal K_x)v\|_{L^2}\le Ch^2\|v\|_{H^1}\le Ch^{3/2}\|w\|_{H^{1/2}}$. For the high frequencies, the first bound of \cref{lem:bem-cont}(ii) for $\mathfrak F(M_1,c_{\rm ca})$ and $\|X'-x'\|_{L^\infty}\le Ch$ give $\|(\mathcal K_X-\mathcal K_x)(w-v)\|_{L^2}\le Ch\|w-v\|_{L^2}\le Ch^{3/2}\|w\|_{H^{1/2}}$. The Cauchy--Schwarz inequality then gives \eqref{eq:bweak-pair}.
\end{proof}

\begin{proof}[Proof of \cref{lem:bem-bweak}]
We assume \ref{K2} with $r=2$ throughout, which \ref{K2} with $r=3$ implies. Since $\mathsf r$ is linear in $g$, it suffices to treat real $g$. The polygon satisfies $X\in\Nd\subset\mathfrak C$ (proof of \cref{lem:bem-geo}), and $g$ and $\Ih g-g$ are real and Lipschitz, because $H^2\subset C^1$. Test with a P1 function $w$, with vertex values $w_j=w(\xi_j)$. Let $\mathcal S_X$ and $\mathcal S_x$ be the continuum forms of \eqref{eq:bem-ext} on the polygon and on the curve, so that $\langle\mathcal S_xg,w\rangle=\int_\Gamma(\Lam_\Gamma g)\,w\,\dd s$ by \cref{prop:bem-dtn}. By the identity $w^\T\Lh(X)\Ih g=\mathfrak s_{X,h}(\Ih g,w)$ of \cref{lem:bem-schur}(b) and the linearity of $\mathfrak s_{X,h}$ in its first argument, we obtain
\begin{multline*}
w^\T \mathsf r=\underbrace{\mathfrak s_{X,h}(\Ih g-g,w)}_{=:\,E_{\rm int}}+\underbrace{\mathfrak s_{X,h}(g,w)-\langle\mathcal S_Xg,w\rangle}_{=:\,E_{\rm Gal}}+\underbrace{\langle(\mathcal S_X-\mathcal S_x)g,w\rangle}_{=:\,E_{\rm geo}}\\
+\underbrace{\int_\Gamma(\Lam_\Gamma g)\,w\,\dd s-w^\T\Mv\Ih(\Lam_\Gamma g)}_{=:\,E_{\rm lump}} .
\end{multline*}

On each panel, $\Ih g-g$ vanishes at the end points and $(\Ih g)'$ is the mean of $g'$, so the Poincar\'e inequality on the panel gives $\|\Ih g-g\|_{L^2}\le h\|(\Ih g-g)'\|_{L^2}\le h^2\|g''\|_{L^2}$. Hence, $\|\Ih g-g\|_{H^1}\le Ch\|g\|_{H^2}$ by $\|b\|_{H^1}\le\|b\|_{L^2}+\|b'\|_{L^2}$, and \eqref{eq:bem-interp} gives $\|\Ih g-g\|_{H^{1/2}}\le Ch^{3/2}\|g\|_{H^2}$. With the bound for $\mathfrak s_{X,h}$ on Lipschitz functions in \cref{lem:bem-schur}(c), this gives
\begin{equation}\label{eq:bweak-int}
|E_{\rm int}|\le C\|\Ih g-g\|_{H^{1/2}}\|w\|_{H^{1/2}}\le C\epsh\|g\|_{H^2}\|w\|_{H^{1/2}} .
\end{equation}

Let $\varsigma_X=\mathcal L_X^{-1}Q_Xg$ and $\varsigma=\mathcal L_x^{-1}Q_xg$, with $Q_Y=\tfrac12+\mathcal K_Y$ for the interior Schur term, and with $Q_Y=\tfrac12-\mathcal K_Y$ and the bordered inverses for the exterior term. For the exterior term, \cref{lem:bem-gauge} replaces \cref{lem:bem-Vreg}(b) throughout. For the interior term, $\mathcal L_x\varsigma=Q_xg$ gives
\[
\mathcal L_X(\varsigma_X-\varsigma)=(\mathcal K_X-\mathcal K_x)g-(\mathcal L_X-\mathcal L_x)\varsigma .
\]
For the exterior term, the same identity holds with $-(\mathcal K_X-\mathcal K_x)g$ in place of $(\mathcal K_X-\mathcal K_x)g$, up to a constant. In both cases, we call the right side without the constant the datum. In the exterior case, $\varsigma_X-\varsigma$ has mean zero, so the uniqueness in \cref{lem:bem-gauge} shows that it is the bordered inverse of the datum. By \cref{lem:bem-schur}(e), we have $\|\varsigma\|_{L^2}\le C\|g\|_{H^1}$. As in the proof of \cref{lem:bem-geo}, \cref{lem:bem-cont}(i) and (ii) bound the datum in $H^1$ by $Ch\|g\|_{H^1}$, so \cref{lem:bem-Vreg}(b) gives $\|\varsigma_X-\varsigma\|_{L^2}\le Ch\|g\|_{H^1}$. By \eqref{eq:bem-sym}, $|E_{\rm Gal}|$ is bounded by $C\|\varsigma_X-\Pi_0\varsigma_X\|_{H^{-1/2}}\|w\|_{H^{1/2}}$, summed over the two Schur terms. Since $\|\varsigma\|_{H^1}\le C\|g\|_{H^2}$ by \cref{lem:bem-higher}(c), \cref{lem:bem-galinv}(c) give
\[
\begin{aligned}
\|\varsigma_X-\Pi_0\varsigma_X\|_{H^{-1/2}}&\le\|(I-\Pi_0)\varsigma\|_{H^{-1/2}}+\|(I-\Pi_0)(\varsigma_X-\varsigma)\|_{H^{-1/2}}\\
&\le C\big(h^{3/2}\|\varsigma\|_{H^1}+h^{1/2}\|\varsigma_X-\varsigma\|_{L^2}\big)\le C\epsh\|g\|_{H^2},
\end{aligned}
\]
and we obtain
\begin{equation}\label{eq:bweak-gal}
|E_{\rm Gal}|\le C\epsh\|g\|_{H^2}\|w\|_{H^{1/2}} .
\end{equation}

By \eqref{eq:bem-ext}, $E_{\rm geo}$ is the sum of the hypersingular part $2\langle(\mathcal L_X-\mathcal L_x)\partial_\xi g,\partial_\xi w\rangle$ and the differences of the two Schur terms. Since $\|\partial_\xi w\|_{H^{-1/2}}\le C\|w\|_{H^{1/2}}$, the first bound of \eqref{eq:bweak-LK} bounds the hypersingular part by $C\epsh\|g\|_{H^2}\|w\|_{H^{1/2}}$. With $\varsigma_X$, $\varsigma$ and $Q_Y$ as in $E_{\rm Gal}$, the difference of the interior Schur terms is
\[
\langle\varsigma_X,Q_Xw\rangle-\langle\varsigma,Q_xw\rangle=\underbrace{\langle\varsigma_X-\varsigma,Q_Xw\rangle}_{=:\,E_{\rm den}}+\underbrace{\langle\varsigma,(\mathcal K_X-\mathcal K_x)w\rangle}_{=:\,E_{\rm pair}} .
\]
In $E_{\rm den}$, the density $\varsigma_X-\varsigma$ is the exact inverse of the datum of $E_{\rm Gal}$. Hence, \eqref{eq:bem-inv}, \eqref{eq:bweak-LK} and the bound for $\|\varsigma\|_{L^2}$ give $\|\varsigma_X-\varsigma\|_{H^{-1/2}}\le C\epsh\|g\|_{H^1}$, and with the bound $\|Q_Xw\|_{H^{1/2}}\le C\|w\|_{H^{1/2}}$ of \cref{lem:bem-Vreg}(c), we obtain $|E_{\rm den}|\le C\epsh\|g\|_{H^2}\|w\|_{H^{1/2}}$. In $E_{\rm pair}$, \eqref{eq:bweak-pair} and the bound for $\|\varsigma\|_{L^2}$ give $|E_{\rm pair}|\le C\epsh\|g\|_{H^2}\|w\|_{H^{1/2}}$. The exterior Schur term is treated in the same way. Hence, we obtain
\begin{equation}\label{eq:bweak-geo}
|E_{\rm geo}|\le C\epsh\|g\|_{H^2}\|w\|_{H^{1/2}} .
\end{equation}

For $E_{\rm lump}$, let $\rho=\mathcal L_{x,0}^{-1}g$. Then, $\Lam_\Gamma g=\rho/|x'|$ by \cref{prop:bem-dtn}, and \cref{lem:bem-higher}(c) gives $\|\rho\|_{H^1}\le C\|g\|_{H^2}$, and $\|\rho\|_{H^2}\le C\|g\|_{H^3}$ under \ref{K2} with $r=3$. The interpolation error satisfies $\|\rho-\Ih\rho\|_{L^2}\le Ch\|\rho\|_{H^1}$ \cite[(10.11)]{Steinbach2008}, and $\|\rho-\Ih\rho\|_{L^2}\le Ch^2\|\rho\|_{H^2}$ under \ref{K2} with $r=3$ \cite[Lemma~10.8]{Steinbach2008}. Since $\rho\in H^1$ is absolutely continuous, we take for $\Lam_\Gamma g$ the continuous representative $\rho/|x'|$, so that $\big(\Mv\Ih(\Lam_\Gamma g)\big)_j=\mu_j\rho_j/|x'(\xi_j)|$ with $\rho_j=\rho(\xi_j)$. Since $\dd s=|x'|\,\dd\xi$, we put $\kappa_j=\mu_j/(h|x'(\xi_j)|)$ and obtain $E_{\rm lump}=E_0+E_\kappa$, where
\[
E_0=\int_\Tor\rho\,w\,\dd\xi-h\sum_jw_j\rho_j,\qquad E_\kappa=h\sum_jw_j\rho_j(1-\kappa_j).
\]

On a panel, the integral of the product of two linear functions with end values $a_0,a_1$ and $b_0,b_1$ differs from the trapezoidal value $\frac h2(a_0b_0+a_1b_1)$ by $-\frac h6(a_1-a_0)(b_1-b_0)$. Summing over the panels, with $\delta_p$ the difference of the end values on panel $p$, we obtain
\[
E_0=\int_\Tor(\rho-\Ih\rho)\,w\,\dd\xi-\frac h6\sum_p\delta_p\rho\,\delta_pw .
\]
Since $\delta_p\rho$ is the integral of $\rho'$ over panel $p$, the Cauchy--Schwarz inequality gives $\sum_p|\delta_p\rho|^2\le h\|\rho'\|_{L^2}^2$, and $\sum_p|\delta_pw|^2=h\|w'\|_{L^2}^2\le C\|w\|_{H^{1/2}}^2$ by \cref{lem:bem-p1}(a). With $\|w\|_{L^2}\le\|w\|_{H^{1/2}}$, we obtain
\begin{equation}\label{eq:bweak-E0}
|E_0|\le\big(\|\rho-\Ih\rho\|_{L^2}+Ch^{3/2}\|\rho\|_{H^1}\big)\|w\|_{H^{1/2}} .
\end{equation}

In $E_\kappa$, the vertices are $X_j=x(\xi_j)$, so $\mu_j/h$ is the mean of the two chord quotients $|x(\xi_j\pm h)-x(\xi_j)|/h$. By Taylor's formula, the difference between each of them and $|x'(\xi_j)|$ is bounded by $h\|x''\|_{L^\infty}$, and $|x'|\ge c_*$ gives $|1-\kappa_j|\le Ch$. Under \ref{K2} with $r=3$, let $u=x'(\xi_j)$ and $z=\frac h2x''(\xi_j)$. Taylor's formula gives $\big(x(\xi_j\pm h)-x(\xi_j)\big)/(\pm h)=u\pm z+A_\pm$ with $|A_\pm|\le Ch^2$, so that $\mu_j/h=\frac12\big(|u+z+A_+|+|u-z+A_-|\big)$, and the triangle inequality shows that the difference between $\mu_j/h$ and $\frac12(|u+z|+|u-z|)$ is bounded by $Ch^2$. The triangle and the Cauchy--Schwarz inequalities give
\[
2|u|\le|u+z|+|u-z|\le\big(2|u+z|^2+2|u-z|^2\big)^{1/2}=2\big(|u|^2+|z|^2\big)^{1/2}\le2|u|+\frac{|z|^2}{|u|} .
\]
The first-order terms thus cancel in the mean of the two chord quotients, and since $|z|^2/|u|\le Ch^2/c_*$, we obtain $|1-\kappa_j|\le Ch^2$. By the bound $h\sum_jv_j^2\le3\|v\|_{L^2}^2$ for P1 functions $v$ before \eqref{eq:bem-lump}, applied to $w$ and to $\Ih\rho$, and by $\|\Ih\rho\|_{L^2}\le\|\rho\|_{L^2}+\|\rho-\Ih\rho\|_{L^2}\le C\|\rho\|_{H^1}$, the Cauchy--Schwarz inequality gives
\begin{equation}\label{eq:bweak-Ek}
|E_\kappa|\le C\max_j|1-\kappa_j|\,\|\rho\|_{H^1}\|w\|_{H^{1/2}} .
\end{equation}

Under \ref{K2} with $r=2$, we add \eqref{eq:bweak-E0} and \eqref{eq:bweak-Ek} with $|1-\kappa_j|\le Ch$, $\|\rho-\Ih\rho\|_{L^2}\le Ch\|\rho\|_{H^1}$ and $\|\rho\|_{H^1}\le C\|g\|_{H^2}$. Under \ref{K2} with $r=3$, we add them with $|1-\kappa_j|\le Ch^2$, $\|\rho-\Ih\rho\|_{L^2}\le Ch^2\|\rho\|_{H^2}$ and $\|\rho\|_{H^1}\le\|\rho\|_{H^2}\le C\|g\|_{H^3}$. Since $h^2\le\epsh\le h$, we obtain, respectively,
\begin{equation}\label{eq:bweak-lump}
|E_{\rm lump}|\le Ch\|g\|_{H^2}\|w\|_{H^{1/2}},\qquad|E_{\rm lump}|\le C\epsh\|g\|_{H^3}\|w\|_{H^{1/2}} .
\end{equation}

Since $\epsh\le h$ and $\|g\|_{H^2}\le\|g\|_{H^3}$, adding \eqref{eq:bweak-int}--\eqref{eq:bweak-geo} to the first, respectively the second, bound of \eqref{eq:bweak-lump} and taking the supremum over $w$ gives the two bounds of the theorem in $\dnorm{\cdot}$.
\end{proof}

\begin{proof}[Proof of the bound in $(H^1)^*$ in \cref{lem:bem-bweak}]
We assume \ref{K2} with $r=3$, and, as in the proof of the bounds in $\dnorm{\cdot}$, we take $g$ real. We use the splitting $w^\T\mathsf r=E_{\rm int}+E_{\rm Gal}+E_{\rm geo}+E_{\rm lump}$ of that proof, together with its notation, and we show that each of the four terms is bounded by $Ch^2\|g\|_{H^3}\|w\|_{H^1}$. In the Schur terms of the bounds in $\dnorm{\cdot}$, the test function appears only through $Q_{\rm int}w$ and $Q_{\rm ext}w$ in $H^{1/2}$, and the densities that represent these functions are only bounded in $H^{-1/2}$. For a test function in $H^1$, these densities lie in $L^2$. In the Galerkin terms, the approximation of $L^2$ densities by $\mathcal P_0$ in $H^{-1/2}$ then gives the additional factor $h^{1/2}$. In the interpolation and geometry terms, the densities allow us to measure the errors in $L^2$, where they are of order $h^2$, instead of in $H^{1/2}$. In the lumping term, the norm of $w$ in $H^1$ replaces the inverse estimate of \cref{lem:bem-p1}(a). Throughout, we use $|\langle u,v\rangle|\le\|u\|_{H^{1/2}}\|v\|_{H^{-1/2}}$ for $u\in H^{1/2}$ and $v\in L^2$, which follows from the Cauchy--Schwarz inequality for the Fourier coefficients, and $\|u'\|_{L^2}\le2\pi\|u\|_{H^1}$ for $u\in H^1$, which follows from the definition of the norms.

We first introduce the densities of the test function. Since $X$ is a polygon of $\mathfrak C$ and $w$ is a P1 function, \cref{lem:bem-Vreg}(c) gives $Q_{\rm int}w,Q_{\rm ext}w\in H^1$ and $\|Q_{\rm int}w\|_{H^1}+\|Q_{\rm ext}w\|_{H^1}\le C\|w\|_{H^1}$. Let $\psi_{\rm int}=\mathcal L_X^{-1}Q_{\rm int}w$, and let $\psi_{\rm ext}=\mathcal L_{X,0}^{-1}Q_{\rm ext}w$ be the bordered inverse, so that $\mathcal L_X\psi_{\rm ext}=Q_{\rm ext}w+c_w$ with a constant $c_w$ and $\int_\Tor\psi_{\rm ext}\,\dd\xi=0$. By \cref{lem:bem-Vreg}(b) and \cref{lem:bem-gauge}, we have
\begin{equation}\label{eq:bweak1-dual}
\|\psi_{\rm int}\|_{L^2}+\|\psi_{\rm ext}\|_{L^2}\le C\|w\|_{H^1},
\end{equation}
and \eqref{eq:bem-gal} gives $\|\psi_\bullet-\Pi_0\psi_\bullet\|_{H^{-1/2}}\le Ch^{1/2}\|w\|_{H^1}$ for $\bullet\in\{{\rm int},{\rm ext}\}$. Since the kernel of $\mathcal L_X$ is symmetric, we obtain
\begin{equation}\label{eq:bweak1-rep}
\langle\varsigma,Q_\bullet w\rangle=\langle\varsigma,\mathcal L_X\psi_\bullet\rangle=\langle\mathcal L_X\varsigma,\psi_\bullet\rangle\qquad(\varsigma\in L^2),
\end{equation}
where, for $\bullet={\rm ext}$, we require that $\varsigma$ has mean zero, so that the constant $c_w$ drops against $\varsigma$. If $\mathcal L_X\varsigma$ is the sum of a function $f$ and a constant, then the constant drops against $\psi_{\rm ext}$, which has mean zero, and $\langle\mathcal L_X\varsigma,\psi_{\rm ext}\rangle=\langle f,\psi_{\rm ext}\rangle$.

For $E_{\rm Gal}$, let $\varsigma_X^{\rm int}=\mathcal L_X^{-1}Q_{\rm int}g$ and $\varsigma_X^{\rm ext}=\mathcal L_{X,0}^{-1}Q_{\rm ext}g$, let $\varsigma_{X,h}^\bullet$ be the corresponding Galerkin solutions, and let $\gamma_\bullet=\varsigma_{X,h}^\bullet-\varsigma_X^\bullet$. By the identity in \cref{lem:bem-schur}(e), we have $E_{\rm Gal}=\langle\gamma_{\rm int},Q_{\rm int}w\rangle+\langle\gamma_{\rm ext},Q_{\rm ext}w\rangle$. In the interior case, the Galerkin orthogonality $\langle\mathcal L_X\gamma_{\rm int},\chi\rangle=0$ holds for all $\chi\in\mathcal P_0$. In the exterior case, $\gamma_{\rm ext}$ has mean zero, and the orthogonality holds for the densities $\chi\in\mathcal P_0$ of mean zero, because the constant of the bordered inverse and the multiplier of the bordered system drop against them. Since $\Pi_0$ preserves the mean, these densities include $\Pi_0\psi_{\rm ext}$. Hence, \eqref{eq:bweak1-rep}, \cref{lem:bem-Vreg}(b) and \eqref{eq:bweak1-dual} give
\[
|\langle\gamma_\bullet,Q_\bullet w\rangle|=|\langle\mathcal L_X\gamma_\bullet,\psi_\bullet-\Pi_0\psi_\bullet\rangle|\le C\|\gamma_\bullet\|_{H^{-1/2}}\|\psi_\bullet-\Pi_0\psi_\bullet\|_{H^{-1/2}}\le C\|\gamma_\bullet\|_{H^{-1/2}}\,h^{1/2}\|w\|_{H^1}
\]
(cf.\ \cite[Thm.~4.2.17]{SauterSchwab2011}). By \eqref{eq:bem-gal}, $\|\gamma_\bullet\|_{H^{-1/2}}$ is bounded by $C\|\varsigma_X^\bullet-\Pi_0\varsigma_X^\bullet\|_{H^{-1/2}}$, and the proof of \eqref{eq:bweak-gal} bounds this norm by $C\epsh\|g\|_{H^2}$. Since $\epsh h^{1/2}=h^2$, we obtain
\begin{equation}\label{eq:bweak1-gal}
|E_{\rm Gal}|\le Ch^2\|g\|_{H^2}\|w\|_{H^1} .
\end{equation}

For $E_{\rm int}$, let $e=\Ih g-g$. The proof of \eqref{eq:bweak-int} gives $\|e\|_{L^2}\le h^2\|g''\|_{L^2}$ and $\|e\|_{H^{1/2}}\le C\epsh\|g\|_{H^2}$. By \eqref{eq:bem-ext}, $E_{\rm int}=\mathfrak s_{X,h}(e,w)$ is the sum of the hypersingular part $2\langle\mathcal L_X\partial_\xi e,\partial_\xi w\rangle$ and the two Schur terms $\langle\beta_\bullet,Q_\bullet w\rangle$, where $\beta_{\rm int}=\mathcal L_{X,h}^{-1}Q_{\rm int}e$ and $\beta_{\rm ext}=\mathcal L_{X,h,0}^{-1}Q_{\rm ext}e$. Since $\partial_\xi w\in\mathcal P_0$, the function $\phi=\mathcal L_X\partial_\xi w$ lies in $H^1$, and \cref{lem:bem-Vreg}(b) gives $\|\phi\|_{H^1}\le C\|\partial_\xi w\|_{L^2}\le C\|w\|_{H^1}$. The function $e$ is Lipschitz and $\phi$ is absolutely continuous, so the symmetry of $\mathcal L_X$ and an integration by parts on $\Tor$ give
\[
|2\langle\mathcal L_X\partial_\xi e,\partial_\xi w\rangle|=|2\langle\partial_\xi e,\phi\rangle|=|2\langle e,\partial_\xi\phi\rangle|\le4\pi\|e\|_{L^2}\|\phi\|_{H^1}\le Ch^2\|g\|_{H^2}\|w\|_{H^1} .
\]
For the Schur terms, the density $\beta_{\rm ext}$ has mean zero, so \eqref{eq:bweak1-rep} gives
\[
\langle\beta_\bullet,Q_\bullet w\rangle=\langle\mathcal L_X\beta_\bullet,\Pi_0\psi_\bullet\rangle+\langle\mathcal L_X\beta_\bullet,\psi_\bullet-\Pi_0\psi_\bullet\rangle .
\]
By the definition of the Galerkin inverses, the first term is $\langle Q_\bullet e,\Pi_0\psi_\bullet\rangle$, because, in the exterior case, the multiplier of the bordered system drops against $\Pi_0\psi_{\rm ext}$, which has mean zero. On the polygon, \cref{lem:bem-Vreg}(c) gives $\|Q_\bullet e\|_{L^2}\le C\|e\|_{L^2}$, and $\Pi_0$ does not increase the norm in $L^2$. Hence, by \eqref{eq:bweak1-dual}, the first term is bounded by $C\|e\|_{L^2}\|w\|_{H^1}\le Ch^2\|g\|_{H^2}\|w\|_{H^1}$. For the second term, \cref{lem:bem-Vreg}(b), \eqref{eq:bem-inv} and \cref{lem:bem-Vreg}(c) give $\|\mathcal L_X\beta_\bullet\|_{H^{1/2}}\le C\|\beta_\bullet\|_{H^{-1/2}}\le C\|Q_\bullet e\|_{H^{1/2}}\le C\|e\|_{H^{1/2}}$, so the second term is bounded by $C\epsh\|g\|_{H^2}\,h^{1/2}\|w\|_{H^1}=Ch^2\|g\|_{H^2}\|w\|_{H^1}$. Adding the three parts, we obtain
\begin{equation}\label{eq:bweak1-int}
|E_{\rm int}|\le Ch^2\|g\|_{H^2}\|w\|_{H^1} .
\end{equation}

For $E_{\rm geo}$, we use the decomposition of the proof of \eqref{eq:bweak-geo}, now with the bounds in $L^2$ of \eqref{eq:bem-geo}. The hypersingular part satisfies
\[
|2\langle(\mathcal L_X-\mathcal L_x)\partial_\xi g,\partial_\xi w\rangle|\le2\|(\mathcal L_X-\mathcal L_x)\partial_\xi g\|_{L^2}\|\partial_\xi w\|_{L^2}\le Ch^2\|g'\|_{L^2}\|w\|_{H^1} .
\]
In $E_{\rm den}=\langle\varsigma_X-\varsigma,Q_Xw\rangle$, let $R$ be the datum of the proof of \eqref{eq:bweak-gal}, that is, $R=(\mathcal K_X-\mathcal K_x)g-(\mathcal L_X-\mathcal L_x)\varsigma$ for the interior term and $R=-(\mathcal K_X-\mathcal K_x)g-(\mathcal L_X-\mathcal L_x)\varsigma$ for the exterior term. The image $\mathcal L_X(\varsigma_X-\varsigma)$ is $R$, up to a constant in the exterior case, in which $\varsigma_X-\varsigma$ has mean zero. Hence, \eqref{eq:bweak1-rep} and the remark after it give $E_{\rm den}=\langle R,\psi_\bullet\rangle$. By \eqref{eq:bem-geo}, we have $\|R\|_{L^2}\le Ch^2(\|g'\|_{L^2}+\|\varsigma\|_{L^2})\le Ch^2\|g\|_{H^1}$, and \eqref{eq:bweak1-dual} gives $|E_{\rm den}|\le Ch^2\|g\|_{H^1}\|w\|_{H^1}$. In $E_{\rm pair}=\langle\varsigma,(\mathcal K_X-\mathcal K_x)w\rangle$, the second bound of \eqref{eq:bem-geo} gives $|E_{\rm pair}|\le Ch^2\|\varsigma\|_{L^2}\|w'\|_{L^2}\le Ch^2\|g\|_{H^1}\|w\|_{H^1}$. The exterior Schur term is treated in the same way, and we obtain
\begin{equation}\label{eq:bweak1-geo}
|E_{\rm geo}|\le Ch^2\|g\|_{H^2}\|w\|_{H^1} .
\end{equation}

For $E_{\rm lump}$, we use the splitting $E_{\rm lump}=E_0+E_\kappa$ and the bounds for $\rho=\mathcal L_{x,0}^{-1}g$ of the proof of \eqref{eq:bweak-lump} under \ref{K2} with $r=3$, that is, $\|\rho\|_{H^1}\le\|\rho\|_{H^2}\le C\|g\|_{H^3}$, $\|\rho-\Ih\rho\|_{L^2}\le Ch^2\|\rho\|_{H^2}$ and $|1-\kappa_j|\le Ch^2$. In $E_0$, the first part is bounded by $\|\rho-\Ih\rho\|_{L^2}\|w\|_{L^2}$. In the second part, we use $\sum_p|\delta_pw|^2=h\|w'\|_{L^2}^2$ directly, instead of \cref{lem:bem-p1}(a), and the Cauchy--Schwarz inequality gives
\[
\frac h6\Big|\sum_p\delta_p\rho\,\delta_pw\Big|\le\frac h6\big(h\|\rho'\|_{L^2}^2\big)^{1/2}\big(h\|w'\|_{L^2}^2\big)^{1/2}\le Ch^2\|\rho\|_{H^1}\|w\|_{H^1} .
\]
The Cauchy--Schwarz inequality in the proof of \eqref{eq:bweak-Ek} gives $|E_\kappa|\le C\max_j|1-\kappa_j|\,\|\rho\|_{H^1}\|w\|_{L^2}$ before $\|w\|_{L^2}$ is bounded by $\|w\|_{H^{1/2}}$. With $\|w\|_{L^2}\le\|w\|_{H^1}$, we obtain
\begin{equation}\label{eq:bweak1-lump}
|E_{\rm lump}|\le Ch^2\|g\|_{H^3}\|w\|_{H^1} .
\end{equation}

Since $\|g\|_{H^1}\le\|g\|_{H^2}\le\|g\|_{H^3}$, adding \eqref{eq:bweak1-gal}--\eqref{eq:bweak1-lump} and taking the supremum over $w$ gives the bound in $(H^1)^*$.
\end{proof}

\begin{lemma}[Shape derivative of the continuum Schur form]\label{lem:bem-dS}
Assume \ref{K1} and \ref{K3}. Let $e$ be a P1 direction, let $x_\varepsilon=x+\varepsilon e$, and let $D\mathcal L_x[e]$ and $D\mathcal K_x[e]$ be the derivatives defined before \cref{lem:bem-c4}. Let $v$ and $w$ be real Lipschitz functions, and let
\[
\begin{aligned}
\varsigma_{\rm int}&=\mathcal L_x^{-1}(\tfrac12+\mathcal K_x)v,&\qquad\psi_{\rm int}&=\mathcal L_x^{-1}(\tfrac12+\mathcal K_x)w,\\
\varsigma_{\rm ext}&=\mathcal L_{x,0}^{-1}(\tfrac12-\mathcal K_x)v,&\qquad\psi_{\rm ext}&=\mathcal L_{x,0}^{-1}(\tfrac12-\mathcal K_x)w .
\end{aligned}
\]
Then, the function $\varepsilon\mapsto\langle\mathcal S_{x_\varepsilon}v,w\rangle$ is differentiable at $\varepsilon=0$, and its derivative $D\langle\mathcal S_xv,w\rangle[e]$ is
\begin{equation}\label{eq:bem-dS}
\begin{aligned}
D\langle\mathcal S_xv,w\rangle[e]={}&2D\mathcal L_x[e](\partial_\xi v,\partial_\xi w)+D\mathcal K_x[e](v,\psi_{\rm int})+D\mathcal K_x[e](w,\varsigma_{\rm int})-D\mathcal L_x[e](\varsigma_{\rm int},\psi_{\rm int})\\
&-D\mathcal K_x[e](v,\psi_{\rm ext})-D\mathcal K_x[e](w,\varsigma_{\rm ext})-D\mathcal L_x[e](\varsigma_{\rm ext},\psi_{\rm ext}).
\end{aligned}
\end{equation}
\end{lemma}
\begin{proof}
For $|\varepsilon|$ small, the curve $x_\varepsilon$ lies in $\mathfrak C$, and $\|x_\varepsilon'-x'\|_{L^\infty}=|\varepsilon|\,\|e'\|_{L^\infty}$. We apply \cref{lem:bem-cont} on the class $\mathfrak F(M_*+c_*/2,c_*/2)\supset\mathfrak C$. Since $\partial_\xi v,\partial_\xi w\in L^2$, the hypersingular term $2\langle\mathcal L_{x_\varepsilon}\partial_\xi v,\partial_\xi w\rangle$ has the derivative $2D\mathcal L_x[e](\partial_\xi v,\partial_\xi w)$ at $\varepsilon=0$.

Consider the interior Schur term $S(\varepsilon)=\langle\mathcal L_{x_\varepsilon}^{-1}Q_\varepsilon v,Q_\varepsilon w\rangle$, where $Q_\varepsilon=\frac12+\mathcal K_{x_\varepsilon}$. Write $\varsigma=\varsigma_{\rm int}$, $\psi=\psi_{\rm int}$ and $\Delta_\varepsilon a=(\mathcal K_{x_\varepsilon}-\mathcal K_x)a$ for Lipschitz $a$, so that $Q_\varepsilon a=Q_0a+\Delta_\varepsilon a$. By \cref{lem:bem-Vreg}(c), $Q_\varepsilon v$ lies in $H^1$, and by (b) there is $\varsigma_\varepsilon\in L^2$ with $\mathcal L_{x_\varepsilon}\varsigma_\varepsilon=Q_\varepsilon v$. Since $\mathcal L_x\psi=Q_0w$, $\mathcal L_x\varsigma=Q_0v$ and $\mathcal L_x$ is symmetric on $L^2$, we have
\[
\langle\varsigma_\varepsilon,Q_0w\rangle=\langle\mathcal L_x\varsigma_\varepsilon,\psi\rangle=\langle\mathcal L_{x_\varepsilon}\varsigma_\varepsilon,\psi\rangle-\langle(\mathcal L_{x_\varepsilon}-\mathcal L_x)\varsigma_\varepsilon,\psi\rangle,\qquad
\langle\mathcal L_{x_\varepsilon}\varsigma_\varepsilon,\psi\rangle=\langle Q_0v,\psi\rangle+\langle\Delta_\varepsilon v,\psi\rangle,
\]
and $\langle Q_0v,\psi\rangle=\langle\varsigma,Q_0w\rangle=S(0)$. We insert these identities into $S(\varepsilon)=\langle\varsigma_\varepsilon,Q_0w\rangle+\langle\varsigma_\varepsilon,\Delta_\varepsilon w\rangle$ and split $\varsigma_\varepsilon=\varsigma+(\varsigma_\varepsilon-\varsigma)$ in the terms that contain $\varsigma_\varepsilon$. This gives
\[
S(\varepsilon)=S(0)+\underbrace{\langle\Delta_\varepsilon v,\psi\rangle+\langle\Delta_\varepsilon w,\varsigma\rangle-\langle(\mathcal L_{x_\varepsilon}-\mathcal L_x)\varsigma,\psi\rangle}_{=:\,T(\varepsilon)}+\underbrace{\langle\Delta_\varepsilon w,\varsigma_\varepsilon-\varsigma\rangle-\langle(\mathcal L_{x_\varepsilon}-\mathcal L_x)(\varsigma_\varepsilon-\varsigma),\psi\rangle}_{=:\,R(\varepsilon)} .
\]
By the definition of $D\mathcal L_x[e]$ and $D\mathcal K_x[e]$, the function $T$ has the derivative $D\mathcal K_x[e](v,\psi)+D\mathcal K_x[e](w,\varsigma)-D\mathcal L_x[e](\varsigma,\psi)$ at $\varepsilon=0$. For $R$, the identities $\mathcal L_{x_\varepsilon}\varsigma_\varepsilon=Q_0v+\Delta_\varepsilon v$ and $\mathcal L_{x_\varepsilon}\varsigma=Q_0v+(\mathcal L_{x_\varepsilon}-\mathcal L_x)\varsigma$ give
\[
\mathcal L_{x_\varepsilon}(\varsigma_\varepsilon-\varsigma)=\Delta_\varepsilon v-(\mathcal L_{x_\varepsilon}-\mathcal L_x)\varsigma .
\]
By \cref{lem:bem-Vreg}(b) for $x_\varepsilon$ and \cref{lem:bem-cont}(i), (ii), we obtain
\[
\|\varsigma_\varepsilon-\varsigma\|_{L^2}\le C\big(\|\Delta_\varepsilon v\|_{H^1}+\|(\mathcal L_{x_\varepsilon}-\mathcal L_x)\varsigma\|_{H^1}\big)\le C|\varepsilon|\,\|e'\|_{L^\infty}\big(\|v\|_{H^1}+\|\varsigma\|_{L^2}\big).
\]
Since $\|\Delta_\varepsilon w\|_{H^1}$ and the norm of $\mathcal L_{x_\varepsilon}-\mathcal L_x$ from $L^2$ to $H^1$ are $O(|\varepsilon|)$ by \cref{lem:bem-cont}(i), (ii), this bound gives $|R(\varepsilon)|\le K\varepsilon^2$, where $K$ depends on $v$, $w$ and $e$. Hence, $R$ has the derivative $0$ at $\varepsilon=0$, so $S$ is differentiable at $\varepsilon=0$, and its derivative is the sum of the interior terms in \eqref{eq:bem-dS}.

The exterior Schur term is treated in the same way, with $Q_\varepsilon a=Q_0a-\Delta_\varepsilon a$ for $Q_\varepsilon=\frac12-\mathcal K_{x_\varepsilon}$, which changes the signs of the two terms of $T$ that contain $\Delta_\varepsilon$, and with the bordered inverses of \cref{lem:bem-gauge}. The identities above hold up to constants in the images under $\mathcal L_x$ and $\mathcal L_{x_\varepsilon}$, which drop against the densities of mean zero, and the difference $\varsigma_\varepsilon-\varsigma_{\rm ext}$, which has mean zero, is the bordered inverse on $x_\varepsilon$ of $-\Delta_\varepsilon v-(\mathcal L_{x_\varepsilon}-\mathcal L_x)\varsigma_{\rm ext}$, so the bound of \cref{lem:bem-gauge} replaces that of \cref{lem:bem-Vreg}(b). Adding the three terms gives \eqref{eq:bem-dS}.
\end{proof}

The second solves in the proof of \cref{lem:bem-c4} need the shape derivatives of \eqref{eq:bem-dS} as functions in $H^1$. For a Lipschitz direction $e$, the quotient $q_e$ of \cref{sec:shape} with $x$ in place of $X$ satisfies $|q_e|\le\|e'\|_{L^\infty}/c_*$.

\begin{lemma}[{Representing functions of the shape derivatives}]\label{lem:bem-rep}
Assume \ref{K1}, \ref{K2} with $r=1$, and \ref{K3}, and let $e$ be a Lipschitz direction.
\begin{enumerate}[label=(\roman*)]
\item For $\varsigma\in L^2$, the function $f_{\mathcal L}(\xi)=-\frac1{2\pi}\int_\Tor\operatorname{Re}q_e(\xi,\eta)\,\varsigma(\eta)\,\dd\eta$ lies in $H^1$, and
\[
\|f_{\mathcal L}\|_{H^1}\le C\|e'\|_{L^\infty}\|\varsigma\|_{L^2},\qquad D\mathcal L_x[e](\varsigma,\varphi)=\langle f_{\mathcal L},\varphi\rangle\quad(\varphi\in L^2).
\]
\item For a real Lipschitz function $a$, the function $f_{\mathcal K}(\xi)=\frac1{2\pi}\int_\Tor\operatorname{Im}q_e(\xi,\eta)\,a'(\eta)\,\dd\eta$ lies in $H^1$, and
\[
\|f_{\mathcal K}\|_{H^1}\le C\|e'\|_{L^\infty}\|a\|_{H^1},\qquad D\mathcal K_x[e](a,\varphi)=\langle f_{\mathcal K},\varphi\rangle\quad(\varphi\in L^2).
\]
\end{enumerate}
\end{lemma}

\begin{proof}
Let $s_0>0$ be so small that $x+se\in\mathfrak C$ and $s|q_e|\le\frac12$ for $0<s\le s_0$. For (i), let $u_s=(\mathcal L_{x+se}\varsigma-\mathcal L_x\varsigma)/s$. Since $x\in\mathfrak C\subset\mathfrak F(M_1,c_{\rm ca})$, the proof of \eqref{eq:c3-H1} applies with $X=x$ and $\lambda=1$, and for $t=s$ small its difference quotient $u_t$ is $u_s$. That proof shows that $u_s(\xi)\to f_{\mathcal L}(\xi)$ for every $\xi$, that $|u_s(\xi)|$ is bounded uniformly in $s$ and $\xi$, and that $f_{\mathcal L}=-\frac1{2\pi}Q_1\varsigma$ lies in $H^1$ with the bound in (i). Hence, by dominated convergence on $\Tor$, $u_s\to f_{\mathcal L}$ in $L^2$. Finally, $\langle u_s,\varphi\rangle$ is the difference quotient at $0$ of $\varepsilon\mapsto\langle\mathcal L_{x+\varepsilon e}\varsigma,\varphi\rangle$, so it tends to $D\mathcal L_x[e](\varsigma,\varphi)$, and it tends to $\langle f_{\mathcal L},\varphi\rangle$ by the convergence in $L^2$.

For (ii), let $u_s=(\mathcal K_{x+se}a-\mathcal K_xa)/s$. By the formula for $\mathcal K_{X_\varepsilon}a-\mathcal K_Xa$ in the proof of \cref{lem:bem-cont}(iii), with $X=x$ and $\varepsilon=(s,0)$, we have
\[
u_s(\xi)=\frac1{2\pi}\int_\Tor\frac{\operatorname{Im}\operatorname{Log}\big(1+sq_e(\xi,\eta)\big)}{s}\,a'(\eta)\,\dd\eta
\]
for almost every $\xi$. Since $|\operatorname{Im}\operatorname{Log}(1+z)|\le\frac32|z|$ for $|z|\le\frac12$, dominated convergence as in the proof of \eqref{eq:c3-H1}, with $a'$ in place of $\varsigma$, gives $u_s(\xi)\to f_{\mathcal K}(\xi)$ for almost every $\xi$, with $|u_s(\xi)|$ bounded uniformly in $s$ and $\xi$, and hence $u_s\to f_{\mathcal K}$ in $L^2$. The second bound of \cref{lem:bem-cont}(ii) gives $\|u_s\|_{H^1}\le C\|e'\|_{L^\infty}\|a\|_{H^1}$. Letting $s\to0$ in each finite part of the series for $\|u_s\|_{H^1}^2$, we obtain $f_{\mathcal K}\in H^1$ with the bound in (ii), and the identification of $D\mathcal K_x[e](a,\varphi)$ follows as in (i).
\end{proof}

\begin{proof}[Proof of \cref{lem:bem-c4}]
Let $X=\Ih x$, let $e$ be a P1 direction and $g\in H^2$, and let $w=\sum_jw_j\varphi_j$ be a P1 function. Since $g$ is continuously differentiable, $g$ and $\Ih g-g$ are Lipschitz. For a curve $Y\in\mathfrak C$, write $\mathfrak s_{Y,h}$ for the Galerkin form \eqref{eq:bem-ext} on the fixed spaces and $\mathcal S_Y$ for the continuum form with exact inverses. Since $\|X'-x'\|_{L^\infty}\le r_1$, the polygon $X+\varepsilon e$ lies in $\mathfrak C$ for $|\varepsilon|\,\|e'\|_{L^\infty}\le r_1$, and $w^\T\Lh(X+\varepsilon e)\Ih g=\mathfrak s_{X+\varepsilon e,h}(\Ih g,w)$ by \cref{lem:bem-schur}(b). Since $\Lh$ is differentiable at $X$ by the paragraph after \cref{def:nbhd}, we obtain
\[
w^\T D\Lh(X)[e]\,\Ih g=\frac{\dd}{\dd\varepsilon}\,\mathfrak s_{X+\varepsilon e,h}(\Ih g,w)\Big|_{\varepsilon=0}=D\mathfrak s_{X,h}[e](\Ih g,w),
\]
where $D\mathfrak s_{Y,h}[e](v,w)$ is the derivative of \cref{thm:bem-c1}. On the side of the continuum, \eqref{eq:bem-dtn-schur} holds for $|\varepsilon|$ small, and $\langle\mathcal S_{x_\varepsilon}g,w\rangle=\sum_jw_j\langle\mathcal S_{x_\varepsilon}g,\varphi_j\rangle$, since the form is linear in $w$. Hence, \eqref{eq:bem-dS} gives
\[
\sum_jw_j\,\frac{\dd}{\dd\varepsilon}\int_\Tor(\Lam_{x_\varepsilon}g)\,\varphi_j\,|x_\varepsilon'|\,\dd\xi\Big|_{\varepsilon=0}=D\langle\mathcal S_xg,w\rangle[e],
\]
and we obtain
\begin{multline*}
w^\T\mathsf r=\underbrace{D\mathfrak s_{X,h}[e](\Ih g,w)-D\mathfrak s_{X,h}[e](g,w)}_{=:\,E_{\rm int}}+\underbrace{D\mathfrak s_{X,h}[e](g,w)-D\mathfrak s_{x,h}[e](g,w)}_{=:\,E_{\rm geo}}\\
+\underbrace{D\mathfrak s_{x,h}[e](g,w)-D\langle\mathcal S_xg,w\rangle[e]}_{=:\,E_{\rm Gal}} .
\end{multline*}

Since the forms $\mathfrak s_{X+\varepsilon e,h}$ are linear in their first argument, we have $E_{\rm int}=D\mathfrak s_{X,h}[e](\Ih g-g,w)$. The first bound of \cref{thm:bem-c1}, with the Lipschitz function $v=\Ih g-g$, and the interpolation bound $\|\Ih g-g\|_{H^{1/2}}\le Ch^{3/2}\|g\|_{H^2}$ of \eqref{eq:bweak-int} give
\begin{equation}\label{eq:c4-int}
|E_{\rm int}|\le Ch^{3/2}\|e'\|_{L^\infty}\|g\|_{H^2}\|w\|_{H^{1/2}}.
\end{equation}

For $E_{\rm geo}$, let $d=\Ih x-x$ and $Y_t=x+td$ for $0\le t\le1$. The curves $Y_t$ lie in $\mathfrak C$, since $\mathfrak C$ is convex, and $E_{\rm geo}=\Phi(1)-\Phi(0)$ with $\Phi(t)=D\mathfrak s_{Y_t,h}[e](g,w)$. Let $F_t(\varepsilon)=\mathfrak s_{Y_t+\varepsilon_1e+\varepsilon_2d,h}(g,w)$. By \cref{thm:bem-c1}, $F_t$ is twice continuously differentiable near $\varepsilon=0$. Hence, its second derivative at $0$ is symmetric. Since $F_{t+\tau}(\varepsilon)=F_t(\varepsilon_1,\varepsilon_2+\tau)$, the function $\Phi$ is differentiable at $t$, and $\Phi'(t)=\partial_{\varepsilon_2}\partial_{\varepsilon_1}F_t(0)=\partial_{\varepsilon_1}\partial_{\varepsilon_2}F_t(0)$ is the mixed second derivative of \cref{thm:bem-c1} at $Y_t$ in the directions $e$ and $d$. By the mixed bound of \cref{thm:bem-c1}, the mean value theorem and $\|d'\|_{L^\infty}\le C_Ih$, we have
\begin{equation}\label{eq:c4-geo}
|E_{\rm geo}|\le C\|e'\|_{L^\infty}\|(\Ih x-x)'\|_{L^\infty}\|g\|_{H^{1/2}}\|w\|_{H^{1/2}}\le Ch\|e'\|_{L^\infty}\|g\|_{H^{1/2}}\|w\|_{H^{1/2}}.
\end{equation}

In $E_{\rm Gal}$, the two forms have the same hypersingular term $2\langle\mathcal L_{x_\varepsilon}\partial_\xi g,\partial_\xi w\rangle$, so only the Schur terms contribute. Write $\mathcal L=\mathcal L_x$, $\mathcal K=\mathcal K_x$, $D\mathcal L=D\mathcal L_x[e]$ and $D\mathcal K=D\mathcal K_x[e]$. Consider the interior term, and let $Q=\frac12+\mathcal K$. Let $\varsigma=\mathcal L^{-1}Qg$ and $\psi=\mathcal L^{-1}Qw$, and let $\varsigma_h=\mathcal L_{x,h}^{-1}Qg$ and $\psi_h=\mathcal L_{x,h}^{-1}Qw$ be their Galerkin approximations. In the notation of the proof of \cref{thm:bem-c1}, with $X=x$, $z=0$ and $\mathsf J=V(0)^{-1}$, the interior Galerkin term is $S_+(\varepsilon)=\ell(Q^\varepsilon_+w)^\T V(\varepsilon)^{-1}\ell(Q^\varepsilon_+g)$. The entries of $V$ and of the loads are functions $F$ and $G$ of \cref{lem:bem-cont}(iii), with the derivatives $D\mathcal L(\iota_{p'},\iota_p)$, $D\mathcal K(g,\iota_p)$ and $D\mathcal K(w,\iota_p)$, and the derivative of $V^{-1}$ is $-\mathsf J V_1\mathsf J$. The vectors $\mathsf J\ell(Qg)$ and $\mathsf J\ell(Qw)$ are the coefficients of $\varsigma_h$ and $\psi_h$. By the product rule, the symmetry of $\mathsf J$ and the linearity of $D\mathcal L$ and $D\mathcal K$ in the densities, we obtain
\[
\partial_{\varepsilon_1}S_+(0)=D\mathcal K(g,\psi_h)+D\mathcal K(w,\varsigma_h)-D\mathcal L(\varsigma_h,\psi_h).
\]
By \eqref{eq:bem-dS}, the derivative of the continuum term is the same expression with $\varsigma$ and $\psi$ in place of $\varsigma_h$ and $\psi_h$. Hence, the interior part of $E_{\rm Gal}$ is
\[
\underbrace{D\mathcal K(w,\varsigma_h-\varsigma)-D\mathcal L(\varsigma_h-\varsigma,\psi_h)}_{=:\,E_1\ \text{(first solve)}}+\underbrace{D\mathcal K(g,\psi_h-\psi)-D\mathcal L(\varsigma,\psi_h-\psi)}_{=:\,E_2\ \text{(second solve)}} .
\]
The exterior part is treated in the same way, with $Q=\frac12-\mathcal K$, the bordered inverses and $\mathsf J=\tilde C(0)$, and the derivative of $\tilde C$ at $0$ is $-\mathsf J V_1\mathsf J$ as well. Its terms with $D\mathcal K$ have the opposite sign, and it splits into $E_1$ and $E_2$ in the same way. We bound the interior part, and the bounds hold for the exterior part with the changes noted.

To bound $E_1$, note that $\|\varsigma\|_{H^1}\le C\|g\|_{H^2}$ by the last statement of \cref{lem:bem-higher}(c). Thus, \cref{lem:bem-galinv}(c) gives $\|\varsigma_h-\varsigma\|_{H^{-1/2}}\le Ch^{3/2}\|g\|_{H^2}$. By \eqref{eq:bem-inv} and \cref{lem:bem-Vreg}(c), we have $\|\psi_h\|_{H^{-1/2}}\le C\|Qw\|_{H^{1/2}}\le C\|w\|_{H^{1/2}}$. By the bounds of \cref{lem:bem-cont}(iii), we have $|D\mathcal K(w,\chi)|\le C\|e'\|_{L^\infty}\|w\|_{H^{1/2}}\|\chi\|_{H^{-1/2}}$ and $|D\mathcal L(\chi,\psi_h)|\le C\|e'\|_{L^\infty}\|\chi\|_{H^{-1/2}}\|\psi_h\|_{H^{-1/2}}$ for $\chi\in L^2$. With $\chi=\varsigma_h-\varsigma$, we obtain
\begin{equation}\label{eq:c4-E1}
|E_1|\le Ch^{3/2}\|e'\|_{L^\infty}\|g\|_{H^2}\|w\|_{H^{1/2}}.
\end{equation}

To bound $E_2$, we apply \cref{lem:bem-rep}(ii) with $a=g$ and \cref{lem:bem-rep}(i) with the density $\varsigma$. This gives $f=f_{\mathcal K}-f_{\mathcal L}\in H^1$ with $D\mathcal K(g,\chi)-D\mathcal L(\varsigma,\chi)=\langle f,\chi\rangle$ for all $\chi\in L^2$. Since $\|\varsigma\|_{L^2}\le\|\varsigma\|_{H^1}\le C\|g\|_{H^2}$, we have
\[
\|f\|_{H^1}\le C\|e'\|_{L^\infty}\big(\|g\|_{H^1}+\|\varsigma\|_{L^2}\big)\le C\|e'\|_{L^\infty}\|g\|_{H^2},
\]
and $E_2=\langle f,\psi_h-\psi\rangle$. Let $\varsigma_2=\mathcal L^{-1}f$, which exists by \cref{lem:bem-Vreg}(b), and let $\varsigma_{2,h}=\mathcal L_{x,h}^{-1}f$ be its Galerkin approximation. By the symmetry of the Galerkin inverse (\cref{lem:bem-galinv}(a)), we have $\langle f,\psi_h\rangle=\langle Qw,\varsigma_{2,h}\rangle$, and by the symmetry of $\mathcal L$, we have $\langle f,\psi\rangle=\langle\mathcal L\varsigma_2,\psi\rangle=\langle\varsigma_2,Qw\rangle$. Hence, we obtain
\[
E_2=\langle Qw,\varsigma_{2,h}-\varsigma_2\rangle .
\]
In the exterior part, $f$ represents $-D\mathcal K(g,\cdot)-D\mathcal L(\varsigma,\cdot)$, the density $\varsigma_2$ is the bordered inverse of $f$, and the same identity holds, because the constants drop against the densities $\psi$, $\psi_h$ and $\varsigma_2$ of mean zero. We have $\|\varsigma_2\|_{L^2}\le C\|f\|_{H^1}$ by \cref{lem:bem-Vreg}(b), or by \cref{lem:bem-gauge} for the bordered inverse, and \eqref{eq:bem-gal} gives $\|\varsigma_{2,h}-\varsigma_2\|_{H^{-1/2}}\le Ch^{1/2}\|f\|_{H^1}$. Since $\|Qw\|_{H^{1/2}}\le C\|w\|_{H^{1/2}}$ by \cref{lem:bem-Vreg}(c), we obtain
\begin{equation}\label{eq:c4-E2-half}
|E_2|\le Ch^{1/2}\|e'\|_{L^\infty}\|g\|_{H^2}\|w\|_{H^{1/2}}.
\end{equation}
For the bound in $(H^1)^*$, let $\vartheta=\mathcal L^{-1}(Qw)$, or the bordered inverse of $Qw$ in the exterior part. By \cref{lem:bem-Vreg}(b), (c), or by \cref{lem:bem-gauge}, we have $\|\vartheta\|_{L^2}\le C\|Qw\|_{H^1}\le C\|w\|_{H^1}$. Let $\rho=\varsigma_{2,h}-\varsigma_2$. By the symmetry of $\mathcal L$ and the Galerkin orthogonality $\langle\mathcal L\rho,\chi\rangle=0$ for $\chi\in\mathcal P_0$, we have
\[
\begin{aligned}
\langle Qw,\rho\rangle=\langle\mathcal L\vartheta,\rho\rangle=\langle\mathcal L\rho,\vartheta-\Pi_0\vartheta\rangle&\le\|\mathcal L\rho\|_{H^{1/2}}\|\vartheta-\Pi_0\vartheta\|_{H^{-1/2}}\\
&\le C\|\rho\|_{H^{-1/2}}\,h^{1/2}\|\vartheta\|_{L^2}\le Ch\|f\|_{H^1}\|w\|_{H^1},
\end{aligned}
\]
by \cref{lem:bem-Vreg}(b) and \eqref{eq:bem-gal} (cf.\ \cite[Thm.~4.2.17]{SauterSchwab2011}). In the exterior part, $\mathcal L\vartheta$ equals $Qw$ up to a constant, which drops against $\rho$, since $\rho$ has mean zero. There, the orthogonality holds for the densities $\chi\in\mathcal P_0$ of mean zero, which include $\Pi_0\vartheta$, since $\Pi_0$ preserves the mean. Hence, we obtain
\begin{equation}\label{eq:c4-E2-one}
|E_2|\le Ch\|e'\|_{L^\infty}\|g\|_{H^2}\|w\|_{H^1}.
\end{equation}

We use $h^{3/2}\le h\le h^{1/2}$, $\|g\|_{H^{1/2}}\le\|g\|_{H^2}$ and $\|w\|_{H^{1/2}}\le\|w\|_{H^1}$. Adding \eqref{eq:c4-int} and \eqref{eq:c4-geo} to \eqref{eq:c4-E1} and \eqref{eq:c4-E2-one} for the interior and the exterior part gives the bound in $(H^1)^*$. The bound in $\dnorm{\cdot}$ follows in the same way, with \eqref{eq:c4-E2-half} in place of \eqref{eq:c4-E2-one}.
\end{proof}

\section*{Formal verification}
Formal counterparts of \Cref{thm:bem-c1p,thm:bem-c3p,thm:bem-c1,lem:bem-bweak,lem:bem-c4} have been verified in Lean 4 \cite{MU21} (version 4.28.0-rc1) using Mathlib \cite{MC20} (revision 789108e1).
The formalization takes \cref{prop:bem-hl,prop:bem-david} and the part (ii) of \cref{prop:bem-dtn}, in the parameter forms used in the proofs, as hypotheses. The discrete operator is defined through the Schur form \eqref{eq:bem-ext}; its identification with the Dirichlet-to-Neumann map of harmonic extensions is not formalized. Subject to these explicitly stated hypotheses and definitions, the formal proofs are kernel-checked and use no additional axioms beyond Lean's standard axioms.
The Lean formalization is publicly available at \url{https://doi.org/10.5281/zenodo.24103886}.
\section*{Use of Generative-AI Tools Declaration}
The author used Claude Code (Anthropic), Codex (OpenAI), and ChatGPT (OpenAI) in the preparation of this manuscript.
These tools assisted in drafting and editing the mathematical exposition and proofs, as well as in developing the Lean 4 formalization mentioned above.
AI assistance was used throughout the manuscript, including its mathematical arguments and presentation.
The author reviewed the resulting mathematical statements, proofs, references, and formalization, and assumes full responsibility for the accuracy and integrity of the manuscript.

\bibliographystyle{abbrvnat}
\bibliography{references}

\begin{thebibliography}{24}
\providecommand{\natexlab}[1]{#1}
\providecommand{\url}[1]{\texttt{#1}}
\expandafter\ifx\csname urlstyle\endcsname\relax
  \providecommand{\doi}[1]{doi: #1}\else
  \providecommand{\doi}{doi: \begingroup \urlstyle{rm}\Url}\fi

\bibitem[Beretta et~al.(2017)Beretta, Francini, and Vessella]{BeFrVe2017}
E.~Beretta, E.~Francini, and S.~Vessella.
\newblock Differentiability of the {Dirichlet} to {Neumann} map under movements
  of polygonal inclusions with an application to shape optimization.
\newblock \emph{SIAM J. Math. Anal.}, 49\penalty0 (2):\penalty0 756--776, 2017.
\newblock ISSN 0036-1410.
\newblock \doi{10.1137/16M1082160}.

\bibitem[Bonnet(1997)]{Bonnet1997}
M.~Bonnet.
\newblock Differentiability of strongly singular and hypersingular boundary
  integral formulations with respect to boundary perturbations.
\newblock \emph{Comput. Mech.}, 19\penalty0 (3):\penalty0 240--246, 1997.
\newblock ISSN 0178-7675.
\newblock \doi{10.1007/s004660050172}.

\bibitem[Bonnet(2001)]{Bonnet2001}
M.~Bonnet.
\newblock Boundary element based formulations for crack shape sensitivity
  analysis.
\newblock \emph{Eng. Anal. Bound. Elem.}, 25\penalty0 (4-5):\penalty0 347--362,
  2001.
\newblock ISSN 0955-7997.
\newblock \doi{10.1016/S0955-7997(01)00025-X}.

\bibitem[Calder{\'o}n(1965)]{Calderon1965}
A.~P. Calder{\'o}n.
\newblock Commutators of singular integral operators.
\newblock \emph{Proc. Natl. Acad. Sci. USA}, 53:\penalty0 1092--1099, 1965.
\newblock ISSN 0027-8424.
\newblock \doi{10.1073/pnas.53.5.1092}.

\bibitem[Chandler-Wilde et~al.(2022)Chandler-Wilde, Hewett, and
  Moiola]{ChandlerWildeHewettMoiola2022}
S.~Chandler-Wilde, D.~Hewett, and A.~Moiola.
\newblock Corrigendum to: ``{Interpolation} of {Hilbert} and {Sobolev} spaces:
  quantitative estimates and counterexamples''.
\newblock \emph{Mathematika}, 68\penalty0 (4):\penalty0 1393--1400, 2022.
\newblock ISSN 0025-5793.
\newblock \doi{10.1112/mtk.12155}.

\bibitem[Chandler-Wilde et~al.(2015)Chandler-Wilde, Hewett, and
  Moiola]{ChandlerWildeHewettMoiola2015}
S.~N. Chandler-Wilde, D.~P. Hewett, and A.~Moiola.
\newblock Interpolation of {Hilbert} and {Sobolev} spaces: quantitative
  estimates and counterexamples.
\newblock \emph{Mathematika}, 61\penalty0 (2):\penalty0 414--443, 2015.
\newblock ISSN 0025-5793.
\newblock \doi{10.1112/S0025579314000278}.

\bibitem[Coifman and Meyer(1978)]{CoifmanMeyer1978}
R.~R. Coifman and Y.~Meyer.
\newblock Commutateurs d'int{\'e}grales singuli{\`e}res et op{\'e}rateurs
  multilin{\'e}aires.
\newblock \emph{Ann. Inst. Fourier}, 28\penalty0 (3):\penalty0 177--202, 1978.
\newblock ISSN 0373-0956.
\newblock \doi{10.5802/aif.708}.

\bibitem[Coifman et~al.(1982)Coifman, McIntosh, and Meyer]{CMM1982}
R.~R. Coifman, A.~McIntosh, and Y.~Meyer.
\newblock L'int{\'e}grale de {Cauchy} d{\'e}finit un op{\'e}rateur born{\'e}
  sur {{\(L^ 2 \)}}pour les courbes lipschitziennes.
\newblock \emph{Ann. Math. (2)}, 116:\penalty0 361--387, 1982.
\newblock ISSN 0003-486X.
\newblock \doi{10.2307/2007065}.

\bibitem[Costabel(1988)]{Costabel1988}
M.~Costabel.
\newblock Boundary integral operators on {Lipschitz} domains: {Elementary}
  results.
\newblock \emph{SIAM J. Math. Anal.}, 19\penalty0 (3):\penalty0 613--626, 1988.
\newblock ISSN 0036-1410.
\newblock \doi{10.1137/0519043}.

\bibitem[Costabel and Le~Lou{\"e}r(2012)]{CostabelLeLouer2012}
M.~Costabel and F.~Le~Lou{\"e}r.
\newblock Shape derivatives of boundary integral operators in electromagnetic
  scattering. {I}: {Shape} differentiability of pseudo-homogeneous boundary
  integral operators.
\newblock \emph{Integral Equations Oper. Theory}, 72\penalty0 (4):\penalty0
  509--535, 2012.
\newblock ISSN 0378-620X.
\newblock \doi{10.1007/s00020-012-1954-z}.

\bibitem[David(1984)]{David1984}
G.~David.
\newblock Op{\'e}rateurs int{\'e}graux singuliers sur certaines courbes du plan
  complexe.
\newblock \emph{Ann. Sci. {\'E}c. Norm. Sup{\'e}r. (4)}, 17:\penalty0 157--189,
  1984.
\newblock ISSN 0012-9593.
\newblock \doi{10.24033/asens.1469}.

\bibitem[de~Moura and Ullrich(2021)]{MU21}
L.~de~Moura and S.~Ullrich.
\newblock The {Lean} 4 theorem prover and programming language.
\newblock In \emph{Automated deduction -- CADE 28. 28th international
  conference on automated deduction, virtual event, July 12--15, 2021.
  Proceedings}, pages 625--635. Cham: Springer, 2021.
\newblock ISBN 978-3-030-79875-8; 978-3-030-79876-5.
\newblock \doi{10.1007/978-3-030-79876-5_37}.

\bibitem[D{\"o}lz and Henr{\'{\i}}quez(2024)]{DolzHenriquez2024}
J.~D{\"o}lz and F.~Henr{\'{\i}}quez.
\newblock Parametric shape holomorphy of boundary integral operators with
  applications.
\newblock \emph{SIAM J. Math. Anal.}, 56\penalty0 (5):\penalty0 6731--6767,
  2024.
\newblock ISSN 0036-1410.
\newblock \doi{10.1137/23M1576451}.

\bibitem[Eto(2026)]{E26-SGBEM}
T.~Eto.
\newblock A structure-preserving symmetric {G}alerkin boundary element method
  for anisotropic multi-phase {M}ullins--{S}ekerka flow.
\newblock Preprint, {arXiv}:2609.14742v2, 2026.
\newblock \doi{10.48550/arXiv.2609.14742}.

\bibitem[Henr{\'{\i}}quez and Schwab(2021)]{HenriquezSchwab2021}
F.~Henr{\'{\i}}quez and C.~Schwab.
\newblock Shape holomorphy of the {Calder{\'o}n} projector for the {Laplacian}
  in {{\(\mathbb{R}^2\)}}.
\newblock \emph{Integral Equations Oper. Theory}, 93\penalty0 (43):\penalty0
  40, 2021.
\newblock ISSN 0378-620X.
\newblock \doi{10.1007/s00020-021-02653-5}.

\bibitem[Henr{\'i}quez et~al.(2026)Henr{\'i}quez, D{\"o}lz, Pinto, and
  Schwab]{HDPS2026SAM}
F.~Henr{\'i}quez, J.~D{\"o}lz, J.~Pinto, and C.~Schwab.
\newblock Shape holomorphy of the {C}alder{\'o}n projector on {L}ipschitz
  surfaces.
\newblock SAM Research Report 2026-42, Seminar for Applied Mathematics, ETH
  Z{\"u}rich, 2026.
\newblock URL
  \url{https://math.ethz.ch/sam/research/reports.leftnav.html?id=1197}.

\bibitem[Mullins and Sekerka(1963)]{MS63}
W.~W. Mullins and R.~F. Sekerka.
\newblock Morphological {S}tability of a {P}article {G}rowing by {D}iffusion or
  {H}eat {F}low.
\newblock \emph{Journal of Applied Physics}, 34\penalty0 (2):\penalty0
  323--329, 1963.
\newblock \doi{10.1063/1.1702607}.

\bibitem[N{\'e}d{\'e}lec(1976)]{Nedelec1976}
J.-C. N{\'e}d{\'e}lec.
\newblock Curved finite element methods for the solution of singular integral
  equations on surfaces in {$\mathbb R^3$}.
\newblock \emph{Comput. Methods Appl. Mech. Engrg.}, 8\penalty0 (1):\penalty0
  61--80, 1976.
\newblock \doi{10.1016/0045-7825(76)90053-0}.

\bibitem[Pinto et~al.(2024)Pinto, Henr{\'{\i}}quez, and
  Jerez-Hanckes]{PHJH2024}
J.~Pinto, F.~Henr{\'{\i}}quez, and C.~Jerez-Hanckes.
\newblock Shape holomorphy of boundary integral operators on multiple open
  arcs.
\newblock \emph{J. Fourier Anal. Appl.}, 30\penalty0 (2):\penalty0 51, 2024.
\newblock ISSN 1069-5869.
\newblock \doi{10.1007/s00041-024-10071-5}.

\bibitem[Potthast(1994)]{Potthast1994}
R.~Potthast.
\newblock Fr{\'e}chet differentiability of boundary integral operators in
  inverse acoustic scattering.
\newblock \emph{Inverse Probl.}, 10\penalty0 (2):\penalty0 431--447, 1994.
\newblock ISSN 0266-5611.
\newblock \doi{10.1088/0266-5611/10/2/016}.

\bibitem[Sauter and Schwab(2011)]{SauterSchwab2011}
S.~A. Sauter and C.~Schwab.
\newblock \emph{Boundary element methods}, volume~39 of \emph{Springer Ser.
  Comput. Math.}
\newblock Berlin: Springer, 2011.
\newblock ISBN 978-3-540-68092-5; 978-3-540-68093-2.
\newblock \doi{10.1007/978-3-540-68093-2}.

\bibitem[Stein(1970)]{Stein1970}
E.~M. Stein.
\newblock \emph{Singular integrals and differentiability properties of
  functions}, volume~30 of \emph{Princeton Math. Ser.}
\newblock Princeton University Press, Princeton, NJ, 1970.
\newblock \doi{10.1515/9781400883882}.

\bibitem[Steinbach(2008)]{Steinbach2008}
O.~Steinbach.
\newblock \emph{Numerical approximation methods for elliptic boundary value
  problems. {Finite} and boundary elements. {Transl}. from the 2003 {German}
  original}.
\newblock New York, NY: Springer, 2008.
\newblock ISBN 978-0-387-31312-2.
\newblock \doi{10.1007/978-0-387-68805-3}.

\bibitem[{The mathlib Community}(2020)]{MC20}
{The mathlib Community}.
\newblock The {Lean} mathematical library.
\newblock In \emph{Proceedings of the 9th ACM SIGPLAN International Conference
  on Certified Programs and Proofs}, CPP 2020, pages 367--381, New York, NY,
  USA, 2020. Association for Computing Machinery.
\newblock ISBN 9781450370974.
\newblock \doi{10.1145/3372885.3373824}.

\end{thebibliography}

\end{document}